\documentclass[11pt,letterpaper]{amsart}
\usepackage[T1]{fontenc}
\usepackage{lmodern}
\usepackage{amssymb,mathtools,mathrsfs}
\usepackage{enumitem}
\usepackage[unicode,hidelinks,linktoc=page,bookmarksnumbered]{hyperref}
\hypersetup{pdftitle={Whitney-Type Real-Analytic Approximation on Banach Spaces I: The C0 and C1 Cases},pdfauthor={Zhouzhe Wang},pdfsubject={Fine real-analytic approximation, Whitney extension, and analytic geometry}}
\allowdisplaybreaks[2]
\numberwithin{equation}{section}
\newcounter{analyticcondition}

\newcounter{proofseries}
\newcounter{proofstep}[proofseries]

\newcommand{\proofstep}[1]{\refstepcounter{proofstep}\textbf{Step \theproofstep. #1}}
\newcommand{\firstproofstep}[1]{\stepcounter{proofseries}\proofstep{#1}}
\theoremstyle{plain}
\newtheorem{theorem}{Theorem}[section]
\newtheorem{lemma}[theorem]{Lemma}
\newtheorem{proposition}[theorem]{Proposition}
\newtheorem{corollary}[theorem]{Corollary}
\theoremstyle{definition}
\newtheorem{definition}[theorem]{Definition}
\theoremstyle{remark}
\newtheorem{remark}[theorem]{Remark}

\begin{document}
\title[Real-Analytic Approximation on Banach Spaces]{Whitney-Type Real-Analytic Approximation on Banach Spaces I: The $C^0$ and $C^1$ Cases}
\author{Zhouzhe Wang}
\address{Beijing International Center for Mathematical Research, Peking University, Beijing 100871, China}
\email{wangzhouzhe@pku.edu.cn}
\date{}
\begin{abstract}
The classical Whitney theorem on real-analytic approximation is a fundamental result in finite-dimensional analysis and geometry. Extending this theory to general infinite-dimensional Banach spaces involves longstanding open problems. This paper is the first part of our study of Whitney-type real-analytic approximation on Banach spaces, focusing on the $C^0$ and $C^1$ cases.

We prove that a real Banach space admits a separating polynomial if and only if it admits a polynomial topological embedding into a real Hilbert space. We also introduce classes of mappings that locally factor through polynomial maps into finite-dimensional spaces and develop methods of real-analytic approximation for these classes.  Using these tools, we establish fine $C^0$ and $C^1$ real-analytic approximation theorems for broad classes of Banach spaces, including nonseparable ones. 

As applications, we resolve several open problems, including those concerning Whitney-type extensions with real analyticity off closed sets, real-analytic approximation preserving Lipschitz constants, and the real-analytic equivalence of diffeomorphic open subsets of Hilbert spaces.

\end{abstract}
\maketitle
\noindent\textit{Declaration on the use of AI assistance.}
We used ChatGPT 6 Astra to assist with the literature review,
identify potential errors in the manuscript,
polish the exposition of several proofs,
and improve the language of the manuscript.
All research ideas, proof strategies, and mathematical
arguments were developed independently by the author.
\tableofcontents

\section[Introduction]{Introduction}\label{sec:main-results}
\subsection{Background and Previous Work}

Whitney approximation is a basic tool in finite-dimensional differential topology.
Whitney established fine $C^r$ real-analytic approximation of $C^r$ maps
between real-analytic submanifolds of Euclidean spaces, for every finite
integer $r\geqslant0$
\cite[Theorem~2, p.~654]{Whitney1936}.
His embedding and real-analytic approximation theorems also imply that
every differentiable manifold admits a compatible real-analytic structure
\cite[Theorems~I--III, pp.~462--463]{Whitney1935},
\cite[Theorem~1, p.~654]{Whitney1936}.

For an open set $U\subset\mathbb R^n$, an integer $r\geqslant0$,
$f\in C^r(U,\mathbb R^m)$, and $\varepsilon\in C(U,(0,\infty))$,
Whitney's approximation theorem
\cite[Lemma~6, pp.~76--78]{Whitney1934},
\cite[Lemma~9, p.~659]{Whitney1936}
gives $g\in C^\omega(U,\mathbb R^m)$ such that, for every $x\in U$,
\begin{equation*}
	\max\limits_{0\leqslant j\leqslant r}
	\|D^jg(x)-D^jf(x)\|<\varepsilon(x).
\end{equation*}

It is therefore natural to seek infinite-dimensional versions of Whitney
approximation. Removing separability assumptions and obtaining simultaneous
fine approximation of derivatives have led to longstanding open problems
\cite[Question~3, p.~150]{GonzaloJaramillo1997},
\cite[second question, p.~125]{AFK2012}.
A basic difficulty in extending the classical convolution argument is the
absence of a nonzero locally finite translation-invariant Borel measure,
already on an infinite-dimensional separable Banach space
\cite[Section~1, p.~218]{HuntSauerYorke1992}.
Thus regularization by convolution with respect to Lebesgue measure does not
transfer directly to the ambient infinite-dimensional space.

In 1954, Kurzweil proved fine $C^0$ real-analytic approximation for continuous
Banach-valued mappings on open subsets of separable real Banach spaces
admitting a separating polynomial
\cite[Theorem~2, pp.~221--222]{Kurzweil1954}.
Gonzalo and Jaramillo asked whether the corresponding whole-space uniform
approximation theorem remains valid without separability
\cite[Question~3, p.~150]{GonzaloJaramillo1997}.

In 1971, Moulis proved fine $C^1$ approximation of Banach-valued $C^1$ mappings
by $C^\infty$ mappings on open subsets of $c_0$ and $\ell^{2m}$,
where $m$ is a positive integer
\cite[Section~1.A, Theorem~1, p.~295]{Moulis1971}.
On open subsets of a separable Hilbert space she also proved fine $C^k$ approximation
of $C^{2k-1}$ mappings by $C^\infty$ mappings for every positive integer $k$
\cite[Section~1.B, Theorem~2, p.~306]{Moulis1971}.

In 1998, Fry obtained uniform real-analytic approximation of uniformly continuous
real-valued functions on open subsets of $c_0$ \cite{Fry1998}.
In 2001, Cepedello Boiso and H\'ajek obtained uniform real-analytic approximation
of uniformly continuous mappings from open subsets of separable Banach spaces
with property $(K)$ into closed convex subsets of arbitrary Banach spaces,
with the approximants taking values in the same closed convex set
\cite[Theorem~1, p.~83]{CepedelloHajek2001}.
The space $c_0$ has property $(K)$ but no separating polynomial
\cite[Proposition~3 and the following remarks, p.~93]{CepedelloHajek2001}.
Whether every continuous real-valued function on $c_0$ admits uniform
real-analytic approximation was asked in
\cite[Question~2, p.~150]{GonzaloJaramillo1997}
and recalled in \cite[first question, p.~125]{AFK2012}.

In 2005, Azagra et al.\ proved fine $C^1$ approximation of $C^1$ mappings
$f:X\to Y$ by $C^p$ mappings, $p\in\mathbb N\cup\{\infty\}$,
when $X$ has an unconditional Schauder basis and a $C^p$ smooth Lipschitz bump,
with $Y$ an arbitrary Banach space
\cite{AzagraEtAl2005}.
In 2010, H\'ajek and Johanis established fine $C^1$ approximation of $C^1$
mappings on open subsets of a separable Banach space by $C^k$ mappings,
provided that the domain admits a $C^k$ smooth Lipschitz bump
and the target has an unconditional Schauder basis and separable dual
\cite{HajekJohanis2010}.

In 2011, Azagra, Fry, and Keener proved uniform real-analytic approximation of both
a bounded Lipschitz $C^1$ function $f:X\to\mathbb R$ and its derivative,
assuming that $Df$ is uniformly continuous and $X$ is a separable Banach space
admitting a separating polynomial
\cite[Theorem~1, p.~954]{AFK2011}.
In their 2012 paper, they also asked whether a $C^k$ function on a Banach space
admitting a separating polynomial can be approximated in the fine topology
together with its derivatives through order $k$ by a real-analytic function
\cite[second question, p.~125]{AFK2012}.

In 2012, for Lipschitz real-valued functions on a separable Banach space with a
separating polynomial, Azagra, Fry, and Keener obtained uniform
real-analytic approximants satisfying
$\operatorname{Lip}(g)\leqslant C_X\operatorname{Lip}(f)$, where $C_X$ depends
only on the domain space
\cite[Theorem~1, p.~126]{AFK2012}.
On a separable Hilbert space they obtained
$\operatorname{Lip}(g)\leqslant\operatorname{Lip}(f)+\eta$ for every $\eta>0$
\cite[Theorem~4, p.~127]{AFK2012}.
Johanis gave a shorter proof of the former theorem \cite{Johanis2012}.

In 2020, Mytrofanov and Ravsky proved uniform real-analytic approximation of continuous
Banach-valued mappings on arbitrary subsets of separable real normed spaces
admitting a separating polynomial, with the approximating mappings defined
on open neighborhoods of those subsets
\cite[Theorem~3, p.~109]{MytrofanovRavsky2020}.
In 2024, Johanis proved that real-valued Lipschitz functions on superreflexive
Banach spaces admit uniform $C^{1,\alpha}$ approximation, for a suitable
$0<\alpha\leqslant1$, without increasing their Lipschitz constants in the given norm
\cite[Theorem~9]{Johanis2024}.

In 2025, Johanis and Zaj\'i\v{c}ek developed $C^1$ Whitney extension theorems under
suitable Lipschitz extension and smooth approximation hypotheses, and
left open basic and Lipschitz $C^1$ extension for $(H,H)$ when $H$ is a
nonseparable Hilbert space \cite[Section~2.6]{JohanisZajicek2025}.

\subsection{Main Results}\label{subsec:intro-main-results}

Let $X$ and $Y$ be real Banach spaces, and let $U\subset X$ be open.
All derivatives in this paper are Fr\'echet derivatives, and
$C^\omega(U,Y)$ denotes the space of real-analytic maps.

For $k\in\{0,1\}$, $f\in C^k(U,Y)$, and
$\varepsilon\in C(U,(0,\infty))$, define
\begin{equation}\label{eq:fine-approximation-definition}
	\mathcal N_k(f,\varepsilon)
	=
	\left\{
	\begin{gathered}
	g\in C^k(U,Y):
	\max\limits_{0\leqslant j\leqslant k}
	\|D^jg(x)-D^jf(x)\|<\varepsilon(x),\\
	\text{for every }x\in U
	\end{gathered}
	\right\},
\end{equation}
where $D^0f=f$ and first derivatives are measured in the operator norm.
These sets form a neighborhood base at $f$ for the
\emph{fine $C^k$ topology} on $C^k(U,Y)$.

We say that $f\in C^k(U,Y)$ admits
\emph{fine $C^k$ real-analytic approximation} if, for every
$\varepsilon\in C(U,(0,\infty))$,
\begin{equation*}
\mathcal N_k(f,\varepsilon)\cap C^\omega(U,Y)
\ne\varnothing.
\end{equation*}

An arbitrary subset $A\subset X$ is endowed with the relative topology,
and $C(A,Y)$ denotes the space of continuous maps from $A$ to $Y$.
We say that $f\in C(A,Y)$ admits \emph{fine $C^0$ real-analytic approximation}
if, for every $\varepsilon\in C(A,(0,\infty))$, there are an open neighborhood
$U\subset X$ of $A$ and a map $g\in C^\omega(U,Y)$ such that, for every $x\in A$,
\begin{equation*}
 \|g(x)-f(x)\|<\varepsilon(x).
\end{equation*}
A continuous real polynomial $P:X\to\mathbb R$ is called a
\emph{separating polynomial}
\cite[Section~1, p.~146]{GonzaloJaramillo1997} if
\begin{equation}\label{eq:separating-polynomial}
 P(0)=0,\qquad\inf\limits_{\|x\|=1}|P(x)|>0.
\end{equation}

\begin{theorem}\label{thm:main-nonseparable}
Let $X$ be a real Banach space admitting a separating polynomial.
For every real Banach space $Y$, subset $A\subset X$, map $f\in C(A,Y)$,
and function $\varepsilon\in C(A,(0,\infty))$, there are an open neighborhood
$U\subset X$ of $A$ and a map $g\in C^\omega(U,Y)$ such that, for every $x\in A$,
\begin{equation}\label{eq:main-nonseparable-approximation}
 \|g(x)-f(x)\|<\varepsilon(x).
\end{equation}
\end{theorem}
Thus the theorem removes the separability assumption from the preceding approximation results and strengthens uniform approximation to fine approximation, even for mappings defined on arbitrary subsets. However, the existence of a separating polynomial is sufficient but not necessary for fine $C^0$ real-analytic approximation. We next establish the same approximation property for the broader class of separable Banach spaces with property $(K)$, which includes $c_0$ despite its failure to admit a separating polynomial.

A separable real Banach space $X$ is said to have \emph{property $(K)$}
\cite[Definition, p.~83]{CepedelloHajek2001}
if there are $d\in C^\omega(X,\mathbb R)$, $R_d>0$, and $\alpha\in\mathbb R$
such that the Taylor series of $d$ at every point has radius of convergence
at least $R_d$ and
\begin{equation}\label{eq:property-K-definition}
 \varnothing\ne\{x\in X:d(x)<\alpha\}\quad\text{is bounded}.
\end{equation}

\begin{theorem}\label{thm:main-K}
Let $X$ be a separable real Banach space with property $(K)$.
For every real Banach space $Y$, subset $A\subset X$, map $f\in C(A,Y)$,
and function $\varepsilon\in C(A,(0,\infty))$, there are an open neighborhood
$U\subset X$ of $A$ and a map $g\in C^\omega(U,Y)$ such that, for every $x\in A$,
\begin{equation}\label{eq:main-K-approximation}
 \|g(x)-f(x)\|<\varepsilon(x).
\end{equation}
\end{theorem}

We next state the first-order approximation results.

Let $I$ be an arbitrary index set, and write
$\operatorname{Fin}(I)=\{F\subset I:\#F<\infty\}$, directed by inclusion.
A family $(E_i)_{i\in I}$ of closed subspaces of $X$ is called an
\emph{unconditional decomposition} of $X$ if every $x\in X$ has a unique family
$(x_i)_{i\in I}\in\prod\limits_{i\in I}E_i$ such that
\begin{equation}\label{eq:unconditional-decomposition-definition}
 \lim\limits_{F\uparrow I}\left\|x-\sum\limits_{i\in F}x_i\right\|=0,
\end{equation}
where $F$ ranges over $\operatorname{Fin}(I)$.
If every $E_i$ is finite-dimensional, this is called an
\emph{unconditional finite-dimensional decomposition}. 
The subspaces $E_i$ are called the blocks of the decomposition.
If $I=\mathbb N$ and every $E_i$ is one-dimensional,
then choosing a nonzero vector $e_i\in E_i$ for each $i\in \mathbb{N}$
gives an unconditional Schauder basis $(e_i)_{i\in\mathbb N}$ of $X$.

\begin{theorem}\label{thm:main-complete}
Let $X$ be a closed complemented subspace of a real Banach space admitting
a separating polynomial and an unconditional finite-dimensional decomposition.
Given a real Banach space $Y$, an open set $U\subset X$, a map $f\in C^1(U,Y)$,
and a function $\varepsilon\in C(U,(0,\infty))$, there is a map
$g\in C^\omega(U,Y)$ such that, for every $x\in U$,
\begin{equation}\label{eq:main-approximation}
 \|g(x)-f(x)\|<\varepsilon(x),\qquad
 \|Dg(x)-Df(x)\|_{\mathcal L(X,Y)}<\varepsilon(x).
\end{equation}
\end{theorem}
\begin{corollary}\label{cor:classical-domain-spaces}
For every positive integer $m$ and every index set $\Gamma$,
the conclusion of Theorem~\ref{thm:main-complete} holds for
$X=\ell^{2m}(\Gamma;\mathbb R)$, $X=L^{2m}(0,1;\mathbb R)$,
and every real Hilbert space $X$.
\end{corollary}
\begin{proof}
The first two classes admit the separating polynomials
\begin{equation*}
 Q_\Gamma(x)=\sum\limits_{\gamma\in\Gamma}x_\gamma^{2m},
 \qquad
 Q_L(f)=\int_0^1 f(t)^{2m}\,\mathrm{d}t,
\end{equation*}
and have, respectively, the standard unconditional basis and the unconditional
Haar basis \cite[Theorem~2, p.~110]{Gaposhkin}.
Every real Hilbert space has an orthonormal basis and the separating polynomial
$x\mapsto\|x\|^2$.
Apply Theorem~\ref{thm:main-complete}.
\end{proof}
A complementary first-order approximation theorem is obtained by imposing an unconditional finite-dimensional decomposition on the target space rather than on the domain, at the expense of requiring the domain to be separable.
\begin{theorem}\label{thm:main-target}
Let $X$ be a separable real Banach space admitting a separating polynomial,
and let $Y$ be a real Banach space admitting an unconditional
finite-dimensional decomposition and an equivalent norm that is $C^1$ smooth
away from the origin.
For every open set $U\subset X$, map $f\in C^1(U,Y)$, and function
$\varepsilon\in C(U,(0,\infty))$, there is a map $g\in C^\omega(U,Y)$ such that, for every $x\in U$,
\begin{equation}\label{eq:main-target-approximation}
 \|g(x)-f(x)\|<\varepsilon(x),\qquad
 \|Dg(x)-Df(x)\|_{\mathcal L(X,Y)}<\varepsilon(x).
\end{equation}
\end{theorem}
Every real Hilbert space admits an unconditional finite-dimensional
decomposition into one-dimensional subspaces induced by an orthonormal
basis, and its Hilbert norm is $C^\infty$ smooth away from the origin.
Thus Theorem~\ref{thm:main-target} yields the following corollary.
\begin{corollary}\label{cor:hilbert-target-space}
Every $C^1$ map from an open subset of a separable real Banach space
admitting a separating polynomial into a real Hilbert space admits
fine $C^1$ real-analytic approximation.
\end{corollary}
\subsection{Outline of the Proofs}
Our proofs begin with a class of $C^1$ mappings on a Hilbert space
$H=\ell^2(\Gamma)$ that locally have the form
$u(x)=b(x_S,\|x\|^2)$, where $S$ is finite.
We obtain real-analytic approximation for these mappings by
Gaussian regularization in finitely many variables, together with
suitable summation estimates. We then obtain fine approximation from uniform approximation by introducing a suitable positive smooth function and approximating the corresponding rescaled mappings.

For Banach spaces admitting separating polynomials, we use projectional resolutions of the identity to reduce the construction to smaller subspaces, and then apply transfinite induction to obtain polynomial topological embeddings into Hilbert spaces. The Hilbert-space approximation results
can then be pulled back to establish fine $C^0$ approximation.
For fine $C^1$ approximation, unconditional finite-dimensional
decompositions allow us to construct smooth approximations
of the identity with locally finite-dimensional range
and controlled derivatives. Combined with local affine
approximation and partitions of unity, these constructions
reduce general $C^1$ mappings to the preceding approximation
framework. When the unconditional finite-dimensional decomposition is imposed on the target space, a similar reduction is achieved through mappings with locally finite-dimensional ranges.

The rest of this paper is organized as follows:
In Section~\ref{ns:sec:notation}, we collect notation and preliminary
results. Section~\ref{sec:hilbert-finite-approximation} introduces locally
finite-dimensional composition representations and establishes the
corresponding real-analytic approximation results on Hilbert spaces.
Section~\ref{sec:nonseparable-c0} characterizes spaces admitting separating
polynomials by their embeddings into Hilbert spaces and proves the fine
$C^0$ approximation theorem.
Sections~\ref{sec:ufdd} and~\ref{sec:target-geometry} establish the fine
$C^1$ approximation results using decompositions of the domain and the
target, respectively. In Section~\ref{sec:lipschitz-upgrade}, we pass from
uniform approximation of Lipschitz mappings to fine $C^0$ approximation
and apply this result to spaces with property $(K)$.
Finally, Section~\ref{sec:applications} presents applications to
approximation preserving Lipschitz constants, Whitney extension,
approximation with surjective derivatives, and real-analytic equivalence
and isotopy of open sets.

\section[Preliminaries]{Preliminaries}\label{ns:sec:notation}

\subsection[Analytic maps and summation]{Analytic maps and summation}
For real or complex Banach spaces $E$ and $F$, let $\mathcal L(E,F)$ denote
the space of bounded linear operators over the corresponding scalar field,
and write $\mathcal L(E)=\mathcal L(E,E)$.
For a set $M$, let $I_M$ denote its identity map.
We use the operator norm topology on
\begin{equation*}
 GL(E)=\{T\in\mathcal L(E):\text{there exists }S\in\mathcal L(E)
              \text{ with }TS=ST=I_E\}.
\end{equation*}
Finite direct sums $\oplus$ carry the maximum norm, while  direct sums $\oplus_2$ carry the square-sum norm.
We use the conventions $\operatorname{dist}(x,\varnothing)=+\infty$ and
$\inf\limits\varnothing=+\infty$.
The open and closed balls with center $a$ and radius $r>0$ are denoted by
$B_E(a,r)$ and $\overline B_E(a,r)$, respectively, and $B_E=B_E(0,1)$.
We write $\|u\|_{\infty,V}=\sup\limits_{x\in V}\|u(x)\|$, or simply
$\|u\|_\infty$ when the domain is clear, and set
$\|u\|_{\infty,\varnothing}=0$.
Supports are taken relative to the domain of the function.
The spaces $\mathbb R^d$ and $\mathbb C^d$ carry their Euclidean norms.
For $D\subset E$, $f:D\to F$, and $A\subset D$, set
\begin{equation*}
 \operatorname{Lip}(f,A)
 =\sup\limits_{\substack{x,y\in A\\x\ne y}}
   \frac{\|f(x)-f(y)\|_F}{\|x-y\|_E},
 \qquad
 \operatorname{Lip}(f)=\operatorname{Lip}(f,D).
\end{equation*}
This quantity is defined to be zero when $A$ has at most one point.

Let $E$ be a real Banach space and let
$p\in\mathbb N\cup\{\infty\}$.
A \emph{$C^p$ smooth Lipschitz bump} on $E$ is a function
$b\in C^p(E,\mathbb R)$ such that $b\not\equiv0$,
$\operatorname{Lip}(b)<\infty$, and there exists $R>0$ with
$\operatorname{supp}b\subset\overline B_E(0,R)$.

For a continuous symmetric $m$-linear map $A:E^m\to F$, write
$A(x^m)=A(x,\ldots,x)$.
For $m\geqslant1$, a continuous $m$-homogeneous polynomial
is a mapping $x\mapsto A(x^m)$ with $A$ as above;
a continuous $0$-homogeneous polynomial is a constant mapping.
A continuous polynomial is a finite sum
$P(x)=\sum\limits_{m=0}^dA_m(x^m)$.
In estimates of polynomial degrees, the zero polynomial is assigned degree $0$.

For real Banach spaces $E,F$ and an open set $U\subset E$,
a mapping $f:U\to F$ belongs to $C^\omega(U,F)$ if, for every $a\in U$,
there exist $r>0$ and continuous $m$-homogeneous polynomials
$P_m:E\to F$, $m\geqslant0$, such that
\begin{equation*}
 \begin{gathered}
 B_E(a,r)\subset U,\qquad
 f(a+h)=\sum\limits_{m=0}^{\infty}P_m(h),\\
 \sum\limits_{m=0}^{\infty}
 \sup\limits_{\|h\|\leqslant r}\|P_m(h)\|<\infty.
 \end{gathered}
\end{equation*}
The expansion holds for every $h\in E$ with $\|h\|<r$.

For a real Banach space $E$, its complexification $E_{\mathbb C}=E+iE$
is equipped with the Taylor complexification norm
\cite[Chapter~1, Definition~57 and equation~(11), p.~23]{HajekJohanis2014},
\begin{equation*}
 \|x+iy\|_{E_{\mathbb C}}
 =\sup\limits_{\theta\in[0,2\pi]}\|x\cos\theta-y\sin\theta\|_E.
\end{equation*}
For Hilbert spaces, we use the usual complex Hilbert norm.
Every continuous polynomial $P:E\to F$ between real Banach spaces
has a unique complex polynomial extension
$P^{\mathbb C}:E_{\mathbb C}\to F_{\mathbb C}$;
see \cite[Chapter~1, Section~2, p.~25]{HajekJohanis2014}.
The complexification of a bounded real linear operator $T:E\to F$ is likewise
denoted by $T_{\mathbb C}=T^{\mathbb C}$.

For complex Banach spaces $E,F$ and an open set $V\subset E$,
a mapping $f:V\to F$ is \emph{Fr\'echet-holomorphic} if, for every $a\in V$,
there exists a bounded complex-linear map $L_a:E\to F$ such that
\begin{equation*}
 \lim\limits_{\substack{h\to0\\h\ne0}}
 \frac{\|f(a+h)-f(a)-L_ah\|}{\|h\|}=0.
\end{equation*}
We write $Df(a)=L_a$ and denote the class of these mappings by
$\mathcal H(V,F)$.
The mapping $f$ is \emph{G\^ateaux holomorphic} if, for every $a\in V$
and $h\in E$, the mapping $\lambda\mapsto f(a+\lambda h)$ is holomorphic
on $\{\lambda\in\mathbb C:a+\lambda h\in V\}$.
For real Banach spaces $X,Y$, a map $f:X\to Y$ is called
\emph{entire real-analytic} if it has a holomorphic extension to all of
$X_{\mathbb C}$ with values in $Y_{\mathbb C}$. For background on holomorphic mappings between complex
Banach spaces, see \cite{Mujica1986}.
All integrals of Banach-valued functions are Bochner integrals.

For an arbitrary index set $I$, the sums of a family $(a_i)_{i\in I}$ of
nonnegative real numbers and a Banach-valued family satisfying
$\sum\limits_{i\in I}\|z_i\|<\infty$ are defined, respectively, by
\begin{equation*}
 \sum\limits_{i\in I}a_i
 =\sup\limits_{F\in\operatorname{Fin}(I)}\sum\limits_{i\in F}a_i,
 \qquad
 \sum\limits_{i\in I}z_i
 =\lim\limits_{F\uparrow I}\sum\limits_{i\in F}z_i.
\end{equation*}
An absolutely summable family $(z_i)_{i\in I}$ has at most countably many
nonzero terms; the second limit is that of the net of finite partial sums.

\subsection[Cardinals and superreflexivity]{Cardinals, well-orderings, and superreflexivity}\label{subsec:density-superreflexivity}

The cardinality of a set $A$ is denoted by $|A|$, and the cardinality of a
countably infinite set by $\aleph_0$.
A total order is a \emph{well-order} if every nonempty subset has a least element.
By the well-ordering theorem and
\cite[Definition~9.13 and Lemma~9.14(1), pp.~135--136]{Moschovakis1994},
every infinite set $I$ admits a well-order $\leqslant$ for which each proper
initial segment has cardinality less than $|I|$; that is, for every $i\in I$,
\begin{equation*}
 I_{<i}=\{j\in I:j<i\},\qquad |I_{<i}|<|I|.
\end{equation*}
Adding one element to a finite initial segment leaves it finite, and adding
one element to an infinite initial segment does not change its cardinality.
Consequently, for every $i\in I$,
\begin{equation}\label{eq:well-order-small-parts}
 I_{\leqslant i}=\{j\in I:j\leqslant i\},\qquad
 |I_{\leqslant i}|<|I|.
\end{equation}
In particular, if $|I|>\aleph_0$, then, for every $i\in I$,
\begin{equation}\label{eq:uncountable-cardinal-segments}
 \max\limits\{\aleph_0,|I_{\leqslant i}|\}<|I|.
\end{equation}
The \emph{norm density} of a normed space $E$, denoted by
$\operatorname{dens}E$, is the least cardinality of a norm-dense subset:
\begin{equation}\label{eq:norm-density-definition}
 \operatorname{dens}E
 =\min\limits\bigl\{|D|:D\subset E,\ \overline D^{\|\cdot\|}=E\bigr\}.
\end{equation}
If $E\ne\{0\}$, then $E$ is separable if and only if
$\operatorname{dens}E=\aleph_0$.
Equivalent norms leave the norm density unchanged.

Let $E$ and $F$ be real Banach spaces.
We say that $F$ is \emph{finitely representable} in $E$ if, for every
finite-dimensional linear subspace $F_0\subset F$ and every $\eta>0$,
there are a finite-dimensional linear subspace $E_0\subset E$ and a linear
isomorphism $T:F_0\to E_0$ such that
\begin{equation}\label{eq:finite-representability-definition}
 \|T\|\,\|T^{-1}\|\leqslant1+\eta.
\end{equation}

A Banach space $E$ is \emph{superreflexive} if every Banach space finitely
representable in $E$ is reflexive; see
\cite[Section~11.1, pp.~428--429]{Pisier2016}.
It is \emph{weakly compactly generated} if there exists a weakly compact
set $K\subset E$ such that $E=\overline{\operatorname{span}}K$.

For a nonzero real Banach space $(E,\|\cdot\|)$, its modulus of convexity
and modulus of smoothness are, respectively,
\begin{equation*}
 \begin{aligned}
 \delta_E(t)
 &=\inf\limits_{\substack{\|x\|=\|y\|=1\\\|x-y\|\geqslant t}}
       \left(1-\left\|\frac{x+y}{2}\right\|\right),
       &&\text{for every }0<t\leqslant2,\\
 \rho_E(t)
 &=\sup\limits_{\|x\|=\|y\|=1}
       \left(\frac{\|x+ty\|+\|x-ty\|}{2}-1\right),
       &&\text{for every }t\geqslant0.
 \end{aligned}
\end{equation*}
The norm is \emph{uniformly convex} if $\delta_E(t)>0$ for every
$0<t\leqslant2$.
It is \emph{uniformly convex of power type $q$}, where $2\leqslant q<\infty$,
if there exists $c>0$ such that
$\delta_E(t)\geqslant ct^q$ for every $0<t\leqslant2$.
Its modulus of smoothness has \emph{power type $p$}, where $1<p\leqslant2$,
if there exists $C>0$ such that
$\rho_E(t)\leqslant Ct^p$ for every $t\geqslant0$.

\subsection[Approximation on arbitrary subsets]{Continuous approximation on arbitrary subsets}\label{subsec:arbitrary-subset-neighborhood}

\begin{lemma}\label{lem:arbitrary-subset-neighborhood}
Let $X$ and $Y$ be real Banach spaces, let $A\subset X$ be an arbitrary subset,
and let $f\in C(A,Y)$ and $\varepsilon\in C(A,(0,\infty))$.
There are an open neighborhood $U\subset X$ of $A$ and locally Lipschitz maps
$F:U\to Y$ and $\delta:U\to(0,\infty)$ such that, for every $x\in A$,
\begin{equation}\label{eq:arbitrary-subset-neighborhood}
 \|F(x)-f(x)\|<\delta(x)<\frac{\varepsilon(x)}{4}.
\end{equation}
Consequently, if every continuous $Y$-valued map on every open subset of $X$
admits fine $C^0$ real-analytic approximation, then so does every continuous
$Y$-valued map on an arbitrary subset of $X$.
\end{lemma}
\begin{proof}
If $A=\varnothing$, take $U=X$, $F=0$, and $\delta=1$.
Assume henceforth that $A\ne\varnothing$.
For each $a\in A$, continuity of $f$ and $\varepsilon$ yields an open
neighborhood $B_a\ni a$ in $X$ such that, for every $x\in A\cap B_a$,
\begin{equation}\label{eq:arbitrary-subset-local-control}
 \|f(x)-f(a)\|<\frac{\varepsilon(a)}{8},
 \qquad \varepsilon(x)>\frac{\varepsilon(a)}{2}.
\end{equation}
Set $U=\bigcup\limits_{a\in A}B_a$.
By paracompactness of metric spaces \cite[Corollary~1, p.~979]{Stone1948},
the open cover $(B_a)_{a\in A}$ has an open refinement $(V_i)_{i\in I}$
that is locally finite in $U$.
For each $i\in I$, choose $a_i\in A$ with $V_i\subset B_{a_i}$.
Inspired by \cite[Theorem~2.31 and Section~4.3, pp.~355, 358--359]{Dugundji1951},
set, for every $x\in U$,
\begin{equation*}
 \begin{gathered}
 b_i(x)=\min\limits\{1,\operatorname{dist}(x,X\setminus V_i)\},\qquad
 b(x)=\sum\limits_{i\in I}b_i(x),\\
 \varphi_i(x)=\frac{b_i(x)}{b(x)}.
 \end{gathered}
\end{equation*}
Each $b_i$ is $1$-Lipschitz. Local finiteness and the covering property imply
that $b$ is continuous, $b>0$, $\varphi_i\geqslant0$, and
$\sum\limits_i\varphi_i=1$.
Fix $a\in U$. Choose an open neighborhood $W\ni a$, a finite set $I_a$,
and a constant $c_a>0$ such that $b_i|_W=0$ for $i\notin I_a$ and
$b|_W\geqslant c_a$.
Then, for $x,y\in W$ and $i\in I_a$,
\begin{equation*}
 |\varphi_i(x)-\varphi_i(y)|
 \leqslant\frac{|b_i(x)-b_i(y)|}{c_a}
       +\frac{|b(x)-b(y)|}{c_a^2}
 \leqslant\left(\frac{1}{c_a}+\frac{\#I_a}{c_a^2}\right)\|x-y\|.
\end{equation*}
Define
\begin{equation*}
 F(x)=\sum\limits_{i\in I}\varphi_i(x)f(a_i),\qquad
 \delta(x)=\frac{1}{8}\sum\limits_{i\in I}\varphi_i(x)\varepsilon(a_i).
\end{equation*}
Set $L_a=\frac{1}{c_a}+\frac{\#I_a}{c_a^2}$. For $x,y\in W$,
\begin{equation*}
 \begin{aligned}
 \|F(x)-F(y)\|
 &\leqslant L_a\left(\sum\limits_{i\in I_a}\|f(a_i)\|\right)\|x-y\|,\\
 |\delta(x)-\delta(y)|
 &\leqslant\frac{L_a}{8}
       \left(\sum\limits_{i\in I_a}\varepsilon(a_i)\right)\|x-y\|.
 \end{aligned}
\end{equation*}
Both sums are finite, so $F$ and $\delta$ are locally Lipschitz.
Moreover, $\delta>0$ because the coefficients $\varphi_i$ are nonnegative
and sum to one.
If $x\in A$ and $\varphi_i(x)>0$, then $x\in V_i\subset B_{a_i}$.
Taking finite convex combinations of the two inequalities in
\eqref{eq:arbitrary-subset-local-control} gives
\begin{equation*}
 \begin{aligned}
 \|F(x)-f(x)\|
 &\leqslant\sum\limits_{i\in I}\varphi_i(x)\|f(a_i)-f(x)\|\\
 &<\frac{1}{8}\sum\limits_{i\in I}\varphi_i(x)\varepsilon(a_i)
 =\delta(x)<\frac{\varepsilon(x)}{4}.
 \end{aligned}
\end{equation*}
If fine $C^0$ real-analytic approximation holds on open sets, choose
$g\in C^\omega(U,Y)$ with $\|g(x)-F(x)\|<\delta(x)$ on $U$.
For $x\in A$, we then have
\begin{equation*}
\|g(x)-f(x)\|
 \leqslant\|g(x)-F(x)\|+\|F(x)-f(x)\|
 <2\delta(x)<\varepsilon(x).\qedhere
\end{equation*}
\end{proof}

\subsection[Partitions of unity and affine terms]{Partitions of unity and differentiation of affine terms}

\begin{lemma}\label{lem:successive-partition}
Let $U\subset X$ be open, let $k\in\{1,2,\ldots\}\cup\{\infty\}$,
and let $(b_n)_{n\geqslant1}\subset C^k(U,[0,1])$.
Suppose that every $x\in U$ has an open neighborhood $V\subset U$
and  $m\in \mathbb{N}$ such that $b_m|_V=1$. Then the family defined by
\begin{equation}\label{eq:successive-partition}
 \psi_n=b_n\prod\limits_{j<n}(1-b_j),\qquad\text{for every }n\geqslant1,
\end{equation}
is a nonnegative locally finite $C^k$ partition of unity on $U$, and
$\operatorname{supp}_U\psi_n\subset\operatorname{supp}_U b_n$.
For $N\geqslant1$ and $x\in U$,
\begin{equation}\label{eq:successive-partition-sum}
 \sum\limits_{n=1}^{N}\psi_n
 =1-\prod\limits_{n=1}^{N}(1-b_n),
\end{equation}
\begin{equation}\label{eq:successive-partition-derivative}
 \|D\psi_n(x)\|\leqslant\sum\limits_{j\leqslant n}\|Db_j(x)\|.
\end{equation}
\end{lemma}
\begin{proof}
The construction \eqref{eq:successive-partition} is inspired by
\cite[Section~1.A, part~A2 of the proof of Theorem~1, pp.~302--303]{Moulis1971}.
Set $p_0=1$ and $p_n=\prod\limits_{j\leqslant n}(1-b_j)$.
Then $\psi_n=p_{n-1}-p_n$, and summing gives
\eqref{eq:successive-partition-sum}.
If $b_m|_V=1$, then $\psi_n|_V=0$ for $n>m$, and
\eqref{eq:successive-partition-sum} with $N\geqslant m$ gives
$\sum\limits_n\psi_n=1$ on $V$.
Thus the family is locally finite; nonnegativity and the support inclusion
follow from \eqref{eq:successive-partition}.
The product rule and $0\leqslant b_j\leqslant1$ yield
\begin{equation*}
 D\psi_n
 =\left(\prod\limits_{j<n}(1-b_j)\right)Db_n
 -b_n\sum\limits_{j<n}\left(\prod\limits_{\substack{l<n\\l\ne j}}(1-b_l)\right)Db_j.
\end{equation*}
Taking operator norms gives \eqref{eq:successive-partition-derivative}.
\end{proof}

\begin{lemma}\label{lem:partition-summation}
Let $U\subset X$ be open, let $k\in\{1,2,\ldots\}\cup\{\infty\}$,
and let $(\varphi_i)_{i\in I}$ be a nonnegative locally finite $C^k$
partition of unity on $U$.
For each $i\in I$, let $V_i\subset U$ be open, with
$\operatorname{supp}_U\varphi_i\subset V_i$, and let $u_i\in C^k(V_i,Y)$.
Then each product $\varphi_i u_i$, extended by zero on $U\setminus V_i$,
belongs to $C^k(U,Y)$, and
\begin{equation*}
 u=\sum\limits_{i\in I}\varphi_i u_i\in C^k(U,Y).
\end{equation*}
For any map $b:U\to Y$, set $I_x=\{i\in I:\varphi_i(x)>0\}$.
Then
\begin{equation}\label{eq:partition-derivative-comparison}
 Du(x)[h]
 =\sum\limits_{i\in I_x}\varphi_i(x)Du_i(x)[h]
  +\sum\limits_{i\in I_x}D\varphi_i(x)[h]\bigl(u_i(x)-b(x)\bigr).
\end{equation}
In particular,
\begin{equation}\label{eq:partition-summation-bound}
 \|Du(x)\|
 \leqslant\sum\limits_{i\in I_x}\varphi_i(x)\|Du_i(x)\|
  +\sum\limits_{i\in I_x}\|D\varphi_i(x)\|\,\|u_i(x)-b(x)\|.
\end{equation}
If $f\in C^1(U,Y)$, setting $e_i=u_i-f|_{V_i}$ gives
\begin{equation}\label{eq:partition-difference-derivative}
 \begin{aligned}
 u(x)-f(x)&=\sum\limits_{i\in I_x}\varphi_i(x)e_i(x),\\
 D(u-f)(x)[h]
 &=\sum\limits_{i\in I_x}\varphi_i(x)De_i(x)[h]
   +\sum\limits_{i\in I_x}D\varphi_i(x)[h]e_i(x).
 \end{aligned}
\end{equation}
\end{lemma}
\begin{proof}
Define, for $x\in U$,
\begin{equation*}
 w_i(x)=
 \begin{cases}
  \varphi_i(x)u_i(x),&x\in V_i,\\
  0,&x\in U\setminus V_i.
 \end{cases}
\end{equation*}
The open sets $V_i$ and $U\setminus\operatorname{supp}_U\varphi_i$ cover $U$.
On their intersection, $\varphi_i u_i=0$, so $w_i\in C^k(U,Y)$ and
$\operatorname{supp}_U w_i\subset\operatorname{supp}_U\varphi_i$.
For every $a\in U$, choose an open neighborhood $W\subset U$ and a finite
set $F_a\subset I$ such that $\varphi_i|_W=0$ for $i\notin F_a$.
Then
\begin{equation*}
 u|_W=\sum\limits_{i\in F_a}w_i|_W\in C^k(W,Y),\qquad
 Du|_W=\sum\limits_{i\in F_a}Dw_i|_W.
\end{equation*}
If $\varphi_i(x)=0$, then $x$ is a minimum of the nonnegative differentiable
function $\varphi_i$, so $D\varphi_i(x)=0$. Therefore, for $h\in X$,
\begin{equation*}
 \begin{aligned}
 Du(x)[h]
 &=\sum\limits_{i\in I_x}
       \bigl(\varphi_i(x)Du_i(x)[h]+D\varphi_i(x)[h]u_i(x)\bigr),\\
 \sum\limits_{i\in I_x}D\varphi_i(x)[h]
 &=D\left(\sum\limits_{i\in I}\varphi_i\right)(x)[h]=D1(x)[h]=0.
 \end{aligned}
\end{equation*}
For any mapping $b:U\to Y$, the latter identity gives
\begin{equation*}
 \sum\limits_{i\in I_x}D\varphi_i(x)[h]u_i(x)
 =\sum\limits_{i\in I_x}D\varphi_i(x)[h](u_i(x)-b(x)).
\end{equation*}
This proves \eqref{eq:partition-derivative-comparison}. Taking the supremum over $\|h\|\leqslant1$ gives
\eqref{eq:partition-summation-bound}.
If $f\in C^1(U,Y)$ and $e_i=u_i-f|_{V_i}$, then
\begin{equation*}
 u-f=\sum\limits_i\varphi_i e_i,\qquad
 D(u-f)(x)[h]
 =\sum\limits_{i\in I_x}
   \bigl(\varphi_i(x)De_i(x)[h]+D\varphi_i(x)[h]e_i(x)\bigr),
\end{equation*}
which is \eqref{eq:partition-difference-derivative}.
\end{proof}

\begin{lemma}\label{lem:affine-difference-gluing}
Let $U\subset X$ be open, let $f\in C^1(U,Y)$, and let
$(\varphi_i)_{i\in I}$ be a nonnegative locally finite $C^1$ partition
of unity on $U$.
For each $i$, let $a_i\in U$, let $V_i\subset U$ be open, with
$\operatorname{supp}_U\varphi_i\subset V_i$, and let $v_i\in C^1(V_i,Y)$.
Extend each product by zero outside $V_i$ and define
\begin{equation}\label{eq:affine-difference-gluing-definition}
 \begin{gathered}
 T_i(x)=f(a_i)+Df(a_i)(x-a_i),\qquad
 r_i=f|_{V_i}-T_i|_{V_i},\qquad e_i=v_i-r_i,\\
 A(x)=\sum\limits_i\varphi_i(x)Df(a_i),\qquad
 b(x)=\sum\limits_i\varphi_i(x)\bigl(f(a_i)-Df(a_i)a_i+v_i(x)\bigr),\\
 h(x)=b(x)+A(x)x.
 \end{gathered}
\end{equation}
Then $A$, $b$, and $h$ are $C^1$ maps, and
\eqref{eq:partition-difference-derivative} holds with $u=h$ and $e_i=v_i-r_i$.
\end{lemma}
\begin{proof}
The construction \eqref{eq:affine-difference-gluing-definition} is inspired by
\cite[Section~1.A, part~A2 of the proof of Theorem~1, p.~303]{Moulis1971}
and \cite[proof of Proposition~17, pp.~578--580]{HajekJohanis2010}.
By Lemma~\ref{lem:partition-summation}, $A$ and $b$ are $C^1$ maps.
The continuous bilinear map $(T,x)\mapsto Tx$ gives $h\in C^1(U,Y)$ and
\begin{equation*}
 Dh(x)[v]=Db(x)[v]+DA(x)[v]x+A(x)v.
\end{equation*}
The definitions and $\sum\limits_i\varphi_i=1$ give
\begin{equation*}
 \begin{aligned}
 h(x)&=\sum\limits_i\varphi_i(x)(T_i(x)+v_i(x)),\\
 h(x)-f(x)&=\sum\limits_i\varphi_i(x)(v_i(x)-r_i(x))
          =\sum\limits_i\varphi_i(x)e_i(x).
 \end{aligned}
\end{equation*}
Applying the product rule to this locally finite sum gives
\begin{equation*}
 D(h-f)(x)[v]
 =\sum\limits_{i\in I_x}\varphi_i(x)De_i(x)[v]
  +\sum\limits_{i\in I_x}D\varphi_i(x)[v]e_i(x),
\end{equation*}
where $I_x=\{i:\varphi_i(x)>0\}$.
\end{proof}

\begin{lemma}\label{lem:affine-analytic-approximation}
Let $U\subset X$ be open, let $b\in C^1(U,Y)$, and let
$A\in C^1(U,\mathcal L(X,Y))$.
If the map $(b,A):U\to Y\oplus\mathcal L(X,Y)$ admits fine $C^1$
real-analytic approximation, then so does $h(x)=b(x)+A(x)x$.
\end{lemma}
\begin{proof}
Given $\eta\in C(U,(0,\infty))$, set
$\delta(x)=\frac{\eta(x)}{2(2+\|x\|)}$.
Using the maximum norm, choose a real-analytic map $(\widehat b,\widehat A)$
such that, for every $x\in U$,
\begin{equation*}
 \max\limits_{j=0,1}
 \max\limits\{\|D^j\widehat b(x)-D^jb(x)\|,
              \|D^j\widehat A(x)-D^jA(x)\|\}<\delta(x).
\end{equation*}
Since the continuous bilinear map $(T,x)\mapsto Tx$ is analytic,
the map $g(x)=\widehat b(x)+\widehat A(x)x$ is real-analytic.
Write $B=\widehat b-b$ and $C=\widehat A-A$.
Then $(g-h)(x)=B(x)+C(x)x$ and
\begin{equation*}
 D(g-h)(x)[v]=DB(x)[v]+DC(x)[v]x+C(x)v.
\end{equation*}
Consequently,
\begin{equation*}
 \begin{aligned}
 \|g(x)-h(x)\|&<(1+\|x\|)\delta(x)<\eta(x),\\
 \|Dg(x)-Dh(x)\|&<(2+\|x\|)\delta(x)<\eta(x).
 \end{aligned}
\end{equation*}
This gives the required fine $C^1$ approximation.
\end{proof}

\subsection[Closed complemented subspaces]{Closed complemented subspaces}

\begin{proposition}\label{prop:permanence}
For $k\in\{0,1\}$, the property of fine $C^k$ real-analytic approximation
for maps with values in arbitrary Banach spaces passes to closed
complemented subspaces.
\end{proposition}
\begin{proof}
Suppose that $X$ has the stated approximation property, and let $E\subset X$
be the range of a bounded projection $\Pi:X\to E$.
Let $Y$ be a real Banach space, let $U\subset E$ be open, and let
$f\in C^k(U,Y)$ and $\varepsilon\in C(U,(0,\infty))$.
Then $\Pi^{-1}(U)$ is open in $X$.
Choose $G\in C^\omega(\Pi^{-1}(U),Y)$ such that, for every
$x\in\Pi^{-1}(U)$,
\begin{equation*}
 \max\limits_{0\leqslant j\leqslant k}
 \|D^jG(x)-D^j(f\circ\Pi)(x)\|<\varepsilon(\Pi x).
\end{equation*}
Set $g=G|_U$. Then $g-f=(G-f\circ\Pi)|_U$, and, when $k=1$,
\begin{equation*}
 Dg(x)-Df(x)=\bigl(DG(x)-D(f\circ\Pi)(x)\bigr)|_E.
\end{equation*}
Since $\|T|_E\|\leqslant\|T\|$ for $T\in\mathcal L(X,Y)$ and
$\Pi x=x$ for $x\in U$, we have
$\max\limits_{0\leqslant j\leqslant k}\|D^jg(x)-D^jf(x)\|<\varepsilon(x)$
for every $x\in U$.
For $k=0$, the last assertion of Lemma~\ref{lem:arbitrary-subset-neighborhood}
gives the version for arbitrary subsets of $E$.
\end{proof}
\section[Analytic Approximation on Hilbert Spaces]{Locally Finite-Dimensional Composition Representations and Analytic Approximation on Hilbert Spaces}\label{sec:hilbert-finite-approximation}

Fix a real Hilbert space $H=\ell^2(\Gamma)$, an open set $U\subset H$, and a real Banach space $Z$.
It suffices to consider the case $H\ne\{0\}$ and $U\ne\varnothing$.

\subsection[Local representations and partitions of unity]{Locally finite-dimensional composition representations and partitions of unity}\label{ns:sec:partition}

Fix a well-ordering of the index set $\Gamma$; finite collections of coordinates are always listed in this order.
For $S\in\operatorname{Fin}(\Gamma)$, write
\begin{equation*}
 \begin{gathered}
 x_S=(x_\gamma)_{\gamma\in S},\qquad
 \Pi_Sx=\sum\limits_{\gamma\in S}x_\gamma e_\gamma,\\
 Q_S=I_H-\Pi_S,\qquad
 \operatorname{supp}x=\{\gamma:x_\gamma\ne0\}.
 \end{gathered}
\end{equation*}
On $H_{\mathbb C}$, define
\begin{equation}\label{ns:eq:qPS}
 q(z)=\sum\limits_{\gamma\in\Gamma}z_\gamma^2,\qquad
 P_S(z)=(z_S,q(z)),\qquad P_\varnothing=q.
\end{equation}
We use the same notation for the restrictions to the real space.
The Cauchy--Schwarz inequality gives
\begin{equation}\label{ns:eq:DP}
 \begin{gathered}
 |q(z)|\leqslant\|z\|^2,\qquad q(x)=\|x\|^2,\qquad
 Dq(x)[v]=2\langle x,v\rangle,\\
 DP_S(x)[v]=(v_S,2\langle x,v\rangle),\qquad
 \|DP_S(x)\|\leqslant1+2\|x\|.
 \end{gathered}
\end{equation}

Define the smooth step function
\begin{equation}\label{ns:eq:smooth-step-function}
 \beta_0(t)=
 \begin{cases}e^{-\frac{1}{t}},&t>0,\\0,&t\leqslant0,\end{cases}
 \qquad
 \sigma(t)=\frac{\beta_0(t)}{\beta_0(t)+\beta_0(1-t)}.
\end{equation}
Then $0\leqslant\sigma\leqslant1$, $\sigma'\geqslant0$,
$\sigma|_{(-\infty,0]}=0$, and $\sigma|_{[1,\infty)}=1$.
Set
\begin{equation}\label{ns:eq:chi-c}
 \chi(t)=\sigma(4t-1),\qquad K_\chi=\|\chi'\|_\infty,\qquad
 c(t)=\begin{cases}\frac{\chi(t)}{t},&t>\frac{1}{4},\\0,&t\leqslant\frac{1}{4}.
 \end{cases}
\end{equation}
Then $c\in C^\infty(\mathbb R)$, $tc(t)=\chi(t)$,
$\|c\|_\infty+\|c'\|_\infty<\infty$,
$\chi=0$ on $(-\infty,\frac{1}{4}]$, and $\chi=1$ on $[\frac{1}{2},\infty)$.
For $u\in C^1(U,Z)$, write
\begin{equation}\label{eq:jetnorm}
 j_1(u;x)=\|u(x)\|+\|Du(x)\|.
\end{equation}
If $B:E_1\times E_2\to F$ is continuous and bilinear with $\|B\|\leqslant1$, then
\begin{equation*}
 D(B(u,v))(x)[h]
 =B(Du(x)[h],v(x))+B(u(x),Dv(x)[h]).
\end{equation*}
Consequently,
\begin{equation}\label{eq:jetproduct}
 j_1(B(u,v);x)\leqslant j_1(u;x)j_1(v;x).
\end{equation}

\begin{definition}\label{ns:def:class}
Let $k\in\mathbb N_0\cup\{\infty\}$ and $u:U\to Z$. We write $u\in\mathscr R_H^k(U,Z)$ if, for every $a\in U$, there are
$S\in\operatorname{Fin}(\operatorname{supp}a)$, an open neighborhood $a\in V\subset U$,
an open set $O\subset\mathbb R^{|S|+1}$, and $b\in C^k(O,Z)$ such that
\begin{equation}\label{ns:eq:local-class}
 P_S(V)\subset O,\qquad u|_V=b\circ P_S|_V,
 \qquad P_S(x)=(x_S,\|x\|^2).
\end{equation}

\end{definition}

\begin{lemma}\label{ns:lem:class-closure}
For $k\in\mathbb N_0\cup\{\infty\}$, one has $\mathscr R_H^k(U,Z)\subset C^k(U,Z)$.
This class is closed under finite products, composition with continuous multilinear mappings, and sums with locally finite families of supports.
Moreover, if $E$ is a real Banach space, then
\begin{equation}\label{ns:eq:class-composition}
 \begin{gathered}
 u\in\mathscr R_H^k(U,E),\quad u(U)\subset W\subset E,\quad W\text{ open},\quad
 \Phi\in C^k(W,Z)\\
 \Longrightarrow\quad \Phi\circ u\in\mathscr R_H^k(U,Z).
 \end{gathered}
\end{equation}
If $V_0\subset U$ is open, $\theta\in\mathscr R_H^k(U,\mathbb R)$,
$u\in\mathscr R_H^k(V_0,Z)$, and $\operatorname{supp}_U\theta\subset V_0$,
then the extension of $\theta u$ by zero to $U\setminus V_0$ belongs to $\mathscr R_H^k(U,Z)$.
\end{lemma}
\begin{proof}
Since $P_S$ is a polynomial, the local representation~\eqref{ns:eq:local-class} gives $\mathscr R_H^k(U,Z)\subset C^k(U,Z)$.
Given a finite family $u_j\in\mathscr R_H^k(U,Z_j)$, choose representations near $a\in U$ of the form
$u_j|_{V_j}=b_j\circ P_{S_j}|_{V_j}$, where $P_{S_j}(V_j)\subset O_j$.
Set
\begin{equation*}
 S=\bigcup\limits_{j=1}^NS_j,\qquad V=\bigcap\limits_{j=1}^NV_j,\qquad
 \pi_j(s,t)=(s_{S_j},t),\qquad
 O=\bigcap\limits_{j=1}^N\pi_j^{-1}(O_j).
\end{equation*}
Then $S\subset\operatorname{supp}a$, $P_S(V)\subset O$, and
\begin{equation}\label{ns:eq:class-common-projection}
 (u_1,\ldots,u_N)|_V
 =(b_1\circ\pi_1,\ldots,b_N\circ\pi_N)\circ P_S|_V.
\end{equation}
This proves the assertion about finite products.
If $u|_V=b\circ P_S|_V$, then $O'=O\cap b^{-1}(W)$ is an open set containing $P_S(V)$, and
$(\Phi\circ u)|_V=(\Phi\circ b|_{O'})\circ P_S|_V$, which proves~\eqref{ns:eq:class-composition}.
For a continuous multilinear mapping $M:Z_1\times\cdots\times Z_N\to Z$,
define $B_M\in C^k(O,Z)$ by
\begin{equation*}
 B_M(w)=M(b_1(\pi_1w),\ldots,b_N(\pi_Nw)).
\end{equation*}
Then $M(u_1,\ldots,u_N)|_V=B_M\circ P_S|_V$.
If $(\operatorname{supp}_Uu_i)_{i\in I}$ is locally finite, choose an open
neighborhood $V_0$ of $a$ and a finite $I_a\subset I$ such that
$u_i|_{V_0}=0$ for $i\notin I_a$. After taking the common representation
\eqref{ns:eq:class-common-projection} on $V\subset V_0$, set
\begin{equation*}
 B_\Sigma(w)=\sum\limits_{i\in I_a}b_i(\pi_iw),\qquad
 \left(\sum\limits_{i\in I}u_i\right)|_V=B_\Sigma\circ P_S|_V.
\end{equation*}
For the last assertion, define
\begin{equation*}
 w(x)=\begin{cases}\theta(x)u(x),&x\in V_0,\\0,&x\in U\setminus V_0.\end{cases}
\end{equation*}
On $V_0$, a common local representation $\theta=b_0\circ P_S$,
$u=b_1\circ P_S$ gives $w=(b_0b_1)\circ P_S$.
On the open set $U\setminus\operatorname{supp}_U\theta$, we have $w=0$.
These open sets cover $U$, and hence $w\in\mathscr R_H^k(U,Z)$.
\end{proof}

\begin{lemma}\label{ns:lem:base}
Let $a\in H$ and let $V$ be an open neighborhood of $a$. There are a finite set
$S\subset\operatorname{supp}a$ and a bounded open rectangle $B\subset\mathbb R^{|S|+1}$ such that
\begin{equation*}
 P_S(a)\in B,\qquad P_S^{-1}(\overline B)\subset V,\qquad
 \overline B\cap\{s_\gamma=0\}=\varnothing.
\end{equation*}
The last equality holds for every $\gamma\in S$. If a finite set $S_0\subset\operatorname{supp}a$ and an open set $O_0$ containing $P_{S_0}(a)$ are also prescribed,
one may additionally require $S\supset S_0$ and $\overline B\subset\pi^{-1}(O_0)$,
where $\pi(s,t)=(s_{S_0},t)$.
\end{lemma}
\begin{proof}
Choose $r>0$ and a finite set $S\subset\operatorname{supp}a$ such that
\begin{equation*}
 B(a,r)\subset V,\qquad \|Q_Sa\|<\frac{r}{4},
 \qquad S\supset S_0\ \text{if $S_0$ is prescribed}.
\end{equation*}
Set
\begin{equation*}
 O=\left\{(s,t):\|s-a_S\|<\frac{r}{4},\quad
                  t-\|s\|^2<\frac{r^2}{16}\right\}.
\end{equation*}
Then $P_S(a)\in O$; if $P_S(x)\in O$, then
\begin{equation*}
 \begin{gathered}
 \|Q_Sx\|^2=\|x\|^2-\|x_S\|^2<\frac{r^2}{16},\\
 \|x-a\|\leqslant\|x_S-a_S\|+\|Q_Sx\|+\|Q_Sa\|<\frac{3r}{4}.
 \end{gathered}
\end{equation*}
Since $a_\gamma\ne0$ for every $\gamma\in S$, there is a bounded open rectangle $B$ such that
\begin{equation*}
 P_S(a)\in B,\qquad
 \overline B\subset O\cap\bigcap\limits_{\gamma\in S}\{s_\gamma\ne0\}.
\end{equation*}
If $O_0$ is prescribed, then $P_S(a)\in\pi^{-1}(O_0)$, so we may also require
$\overline B\subset\pi^{-1}(O_0)$. Thus $P_S^{-1}(\overline B)\subset B(a,r)\subset V$.
\end{proof}
We write $B'\Subset B$ when $\overline{B'}$ is a compact subset of $B$.

Fix an open cover $\mathcal V$ of $U$, a mapping $u\in\mathscr R_H^1(U,Z)$,
and $\tau\in C(U,(0,\infty))$.
Inspired by \cite[proof of Theorem~2, pp.~696--698]{Johanis2017}, enumerate, with indices $j\geqslant1$, all triples
\begin{equation*}
 (m_j,B'_j,B_j)
\end{equation*}
such that $m_j\in\mathbb N_0$, the sets $B'_j,B_j\subset\mathbb R^{m_j+1}$
are open rectangles with rational endpoints, $B'_j\Subset B_j$, and
$\overline B_j\cap\{s_l=0\}=\varnothing$ for every $1\leqslant l\leqslant m_j$.
Write $B_j=\prod\limits_{\ell=1}^{m_j+1}(a_\ell,b_\ell)$ and
$B'_j=\prod\limits_{\ell=1}^{m_j+1}(a'_\ell,b'_\ell)$.
For every $1\leqslant\ell\leqslant m_j+1$, choose numbers satisfying
\begin{equation*}
 a_\ell<c_\ell<d_\ell<a'_\ell<b'_\ell<e_\ell<f_\ell<b_\ell
\end{equation*}
and define, using \eqref{ns:eq:smooth-step-function},
\begin{equation*}
 \beta_j(y)=\prod\limits_{\ell=1}^{m_j+1}
 \sigma\left(\frac{y_\ell-c_\ell}{d_\ell-c_\ell}\right)
 \sigma\left(\frac{f_\ell-y_\ell}{f_\ell-e_\ell}\right).
\end{equation*}
Then $\beta_j\in C_c^\infty(B_j,[0,1])$,
$\operatorname{supp}\beta_j\subset\prod\limits_\ell[c_\ell,f_\ell]\Subset B_j$,
and $\beta_j=1$ on $\prod\limits_\ell(d_\ell,e_\ell)\supset\overline B'_j$.
For $i=(S,j)$ with $|S|=m_j$, write
\begin{equation}\label{ns:eq:catalogue-atoms}
 S_i=S,\quad j(i)=j,\quad m_i=|S_i|,\quad
 w_i=\beta_j\circ P_{S_i},\quad
 K_i=P_{S_i}^{-1}(\operatorname{supp}\beta_j).
\end{equation}
If $S_i\subset\operatorname{supp}a$, then~\eqref{ns:eq:catalogue-atoms} gives a factorization of $w_i$ through $P_{S_i}$ near $a$.
If $S_i\not\subset\operatorname{supp}a$, then $P_{S_i}(a)\notin\operatorname{supp}\beta_j$.
The continuity of $P_{S_i}$ and the closedness of $\operatorname{supp}\beta_j$ imply that $w_i$ vanishes on a neighborhood of $a$; in this case, take $S=\varnothing$ and $b=0$ in Definition~\ref{ns:def:class}.
Hence $w_i\in\mathscr R_H^\infty(H,\mathbb R)$.
Let $I$ be the set of indices for which there are
$V_i\in\mathcal V$, an open set $V_i^{(0)}\subset U$, an open set
$O_i\subset\mathbb R^{m_i+1}$, and $b_i\in C^1(O_i,Z)$ satisfying
\begin{equation}\label{ns:eq:catalogue-representation}
 \begin{gathered}
 K_i\subset V_i\cap V_i^{(0)},\qquad
 \operatorname{supp}\beta_{j(i)}\Subset O_i,\qquad
 P_{S_i}(V_i^{(0)})\subset O_i,\\
 u|_{V_i^{(0)}}=b_i\circ P_{S_i}|_{V_i^{(0)}}.
 \end{gathered}
\end{equation}
For every $i\in I$, fix $V_i,V_i^{(0)},O_i,b_i$ satisfying~\eqref{ns:eq:catalogue-representation}, and write
\begin{equation*}
 N_i=\sup\limits_{w\in\operatorname{supp}\beta_{j(i)}}
       \bigl(\|b_i(w)\|+\|Db_i(w)\|\bigr)<\infty.
\end{equation*}

\begin{lemma}\label{ns:lem:catalog}
There are increasing index sets $I_n\subset I$ such that $j(i),m_i,N_i\leqslant n$ for every $i\in I_n$, and
\begin{equation}\label{ns:eq:tau-support}
 K_i\subset\{\tau>\frac{1}{n}\},\qquad
 \#\{i\in I_n:S_i=S\}\leqslant n.
\end{equation}
Each family $(K_i)_{i\in I_n}$ is locally finite in $H$, and
\begin{equation}\label{ns:eq:catalogue-layer}
 W_n=\sum\limits_{i\in I_n}w_i\in\mathscr R_H^\infty(H,\mathbb R).
\end{equation}
There are constants $a_n>0$, $R_n\geqslant1$, and $\widehat B_n,E_n<\infty$ such that
\begin{align}
 |s_l|&\geqslant a_n,
 &K_i&\subset\overline B(0,R_n),\label{ns:eq:anRn}\\
 0&\leqslant W_n(x)\leqslant \widehat B_n,
 &\sum\limits_{i\in I_n}\|Dw_i(x)\|&\leqslant E_n.\label{ns:eq:W-bounds}
\end{align}
The coordinate estimate in~\eqref{ns:eq:anRn} holds for every $i\in I_n$, $(s,t)\in\operatorname{supp}\beta_{j(i)}$,
and $1\leqslant l\leqslant m_i$, and the inclusion holds for every $i\in I_n$;
the estimates in~\eqref{ns:eq:W-bounds} hold for every $x\in H$.
For every $a\in U$, there are an open neighborhood $V\subset U$ of $a$ and $n\geqslant1$ such that $W_n|_V\geqslant1$.
If $W_n(x)>0$, then
\begin{equation}\label{ns:eq:Ln}
 j_1(u;x)\leqslant L_n:=n(2+2R_n).
\end{equation}
\end{lemma}
\begin{proof}
\noindent\firstproofstep{}\label{step:catalogue:1}\quad
Given $a\in U$, use Definition~\ref{ns:def:class} to choose
$u|_{V^{(0)}}=b\circ P_{S_0}|_{V^{(0)}}$, where $a\in V^{(0)}\subset U$,
$S_0\subset\operatorname{supp}a$, $b\in C^1(O,Z)$, and $P_{S_0}(V^{(0)})\subset O$.
Choose $V_a\in\mathcal V$ containing $a$, and apply Lemma~\ref{ns:lem:base} to the open neighborhood
$V^{(0)}\cap V_a\cap\{\tau>\frac{\tau(a)}{2}\}$ with the prescribed sets $S_0,O$.
This gives $S\supset S_0$ and rational rectangles $P_S(a)\in B'\Subset B$ such that
\begin{equation*}
 \begin{gathered}
 \pi(s,t)=(s_{S_0},t),\qquad \overline B\subset\pi^{-1}(O),\\
 P_S^{-1}(\overline B)\subset V^{(0)}\cap V_a\cap\{\tau>\frac{\tau(a)}{2}\}.
 \end{gathered}
\end{equation*}
Let $j$ be the index of the triple $(|S|,B',B)$, and let $i=(S,j)$. Set
\begin{equation*}
 V_i=V_a,\qquad V_i^{(0)}=V^{(0)},\qquad O_i=\pi^{-1}(O),\qquad b_i=b\circ\pi.
\end{equation*}
Since $P_{S_0}(V^{(0)})\subset O$ and $\operatorname{supp}\beta_j\Subset B$,
these choices satisfy~\eqref{ns:eq:catalogue-representation}. Thus $i\in I$,
$P_{S_i}(a)\in B'_{j(i)}$, and $K_i\subset\{\tau>\frac{\tau(a)}{2}\}$. Therefore
\begin{equation}\label{ns:eq:catalog-cover}
 U=\bigcup\limits_{i\in I}P_{S_i}^{-1}(B'_{j(i)}).
\end{equation}
Set
\begin{equation*}
 I_n=\{i\in I:j(i),m_i,N_i\leqslant n,\quad K_i\subset\{\tau>\frac{1}{n}\}\}.
\end{equation*}
For $n>\max\limits\{j(i),m_i,N_i,\frac{2}{\tau(a)}\}$,
one has $i\in I_n$ and $W_n\geqslant w_i=1$ on $P_{S_i}^{-1}(B'_{j(i)})$.

\noindent\proofstep{}\label{step:catalogue:2}\quad
Let $\mathcal J_n=\{j\leqslant n:m_j\leqslant n\}$.
This set is finite, and each $\operatorname{supp}\beta_j$ is compact and disjoint from the corresponding coordinate hyperplanes. We may therefore choose
$a_n>0$ and $R_n\geqslant1$ such that, for every $j\in\mathcal J_n$ and $(s,t)\in\operatorname{supp}\beta_j$,
\begin{equation*}
 |s_l|\geqslant a_n,\qquad t\leqslant R_n^2.
\end{equation*}
The first inequality holds for every $1\leqslant l\leqslant m_j$.
If $x\in K_i$, then $t=q(x)=\|x\|^2$, and~\eqref{ns:eq:anRn} follows.
Fix $a\in H$ and set
\begin{equation*}
 F=\{\gamma:|a_\gamma|\geqslant\frac{a_n}{2}\},\qquad
 |F|\leqslant\frac{4\|a\|^2}{a_n^2}<\infty.
\end{equation*}
If $x\in K_i\cap B(a,\frac{a_n}{4})$ and $\gamma\in S_i$, then
\begin{equation*}
 |a_\gamma|\geqslant|x_\gamma|-|x_\gamma-a_\gamma|
 >\frac{3a_n}{4}>\frac{a_n}{2},
\end{equation*}
Hence $S_i\subset F$, and therefore
\begin{equation*}
 \#\left\{i\in I_n:K_i\cap B\left(a,\frac{a_n}{4}\right)\ne\varnothing\right\}
 \leqslant n2^{|F|}<\infty.
\end{equation*}
This proves that $(K_i)_{i\in I_n}$ is locally finite.
Since $w_i\in\mathscr R_H^\infty(H,\mathbb R)$, Lemma~\ref{ns:lem:class-closure} gives
$W_n\in\mathscr R_H^\infty(H,\mathbb R)$.
If $w_i(x)\ne0$ or $Dw_i(x)\ne0$, then $x\in K_i$,
$\|x\|\leqslant R_n$, and
$S_i\subset\{\gamma:|x_\gamma|\geqslant a_n\}$.
This coordinate set has at most $\lfloor\frac{R_n^2}{a_n^2}\rfloor$ elements. Consequently, the number of such indices $i$ is at most
\begin{equation*}
 N_n=n2^{\left\lfloor\frac{R_n^2}{a_n^2}\right\rfloor}
\end{equation*}
Using $0\leqslant w_i\leqslant1$ and~\eqref{ns:eq:DP}, we may take
\begin{equation*}
 \widehat B_n=N_n,\qquad
 E_n=N_n(1+2R_n)\max\limits\bigl(\{1\}\cup\{\|D\beta_j\|_\infty:j\in\mathcal J_n\}\bigr),
\end{equation*}
which proves~\eqref{ns:eq:W-bounds}. If $W_n(x)>0$, choose $i\in I_n$ such that $w_i(x)>0$.
The chain rule then gives
\begin{equation*}
j_1(u;x)\leqslant N_i(1+\|DP_{S_i}(x)\|)
 \leqslant n(2+2R_n)=L_n.\qedhere
\end{equation*}
\end{proof}

Inspired by \cite[Section~4, Lemma~12, p.~574]{HajekJohanis2010}, set, for every $n\geqslant1$,
\begin{align}
 \vartheta_n&=\chi(W_n)\prod\limits_{j<n}(1-\chi(W_j)),\label{ns:eq:vartheta}\\
 \phi_{n,i}&=w_i c(W_n)\prod\limits_{j<n}(1-\chi(W_j)),
       \quad\text{for every }i\in I_n.\label{ns:eq:partitionphi}
\end{align}
Since $w_i,W_j\in\mathscr R_H^\infty(H,\mathbb R)$ and $c,\chi\in C^\infty(\mathbb R)$,
Lemma~\ref{ns:lem:class-closure} yields
$\vartheta_n,\phi_{n,i}\in\mathscr R_H^\infty(H,\mathbb R)$.
Both families are nonnegative because $w_i\geqslant0$, $c\geqslant0$, and $0\leqslant\chi\leqslant1$.
The identity $tc(t)=\chi(t)$ gives
\begin{equation*}
 \sum\limits_{i\in I_n}\phi_{n,i}
 =W_nc(W_n)\prod\limits_{j<n}(1-\chi(W_j))
 =\chi(W_n)\prod\limits_{j<n}(1-\chi(W_j))=\vartheta_n.
\end{equation*}
For each $a\in U$, Lemma~\ref{ns:lem:catalog} provides an open neighborhood $V$ and an index $h$ such that $W_h|_V\geqslant1$,
so that $\chi(W_h)|_V=1$. Applying Lemma~\ref{lem:successive-partition} with $b_n=\chi(W_n)$,
we find that $(\vartheta_n)$ is locally finite in $U$ and $\sum\limits_n\vartheta_n=1$.
For every $n>h$, equation~\eqref{ns:eq:partitionphi} gives $\phi_{n,i}|_V=0$.
For each $1\leqslant n\leqslant h$, the supports of
$(\phi_{n,i})_{i\in I_n}$ are contained in the locally finite
family $(K_i)_{i\in I_n}$.
Since all terms with $n>h$ vanish on $V$, the family
$(\phi_{n,i})_{n\geqslant1,\,i\in I_n}$ is locally finite in $U$.
It follows that
\begin{equation}\label{ns:eq:partition-one}
 \sum\limits_{n\geqslant1}\sum\limits_{i\in I_n}\phi_{n,i}=1,
 \qquad \operatorname{supp}_U\phi_{n,i}\subset K_i.
\end{equation}

\begin{proposition}\label{ns:prop:partition}
Every open cover of $U$ admits a subordinate locally finite partition of unity with all terms in
$\mathscr R_H^\infty(U,\mathbb R)$.
\end{proposition}
\begin{proof}
For the given open cover $\mathcal V$, take $u=0$ and $\tau=1$, and construct $(\phi_{n,i})$ using Lemma~\ref{ns:lem:catalog} and~\eqref{ns:eq:partitionphi}.
The local finiteness established above and~\eqref{ns:eq:partition-one} give the required partition of unity. Subordination follows from
$\operatorname{supp}_U\phi_{n,i}\subset K_i\subset V_i\in\mathcal V$.
\end{proof}

\begin{lemma}\label{ns:lem:small-positive}
For every $\tau\in C(U,(0,\infty))$, there is
$\rho\in\mathscr R_H^\infty(U,\mathbb R)$ such that, for every $x\in U$,
\begin{equation}\label{ns:eq:small-positive}
 0<\rho(x),\qquad j_1(\rho;x)<\tau(x).
\end{equation}
\end{lemma}
\begin{proof}
Take $u=0$ and construct the functions $\vartheta_n$ corresponding to $\tau$ using Lemma~\ref{ns:lem:catalog} and~\eqref{ns:eq:vartheta}.
In view of~\eqref{eq:successive-partition-derivative} and
$\|D(\chi\circ W_j)\|_\infty\leqslant K_\chi E_j$, write
\begin{equation}\label{ns:eq:Mn}
 M_n=K_\chi\sum\limits_{j\leqslant n}E_j,
 \qquad \|D\vartheta_n\|_\infty\leqslant M_n,
\end{equation}
and set
\begin{equation*}
 \rho=\sum\limits_{n\geqslant1}\frac{2^{-n}}{n(1+M_n)}\vartheta_n.
\end{equation*}
The local finiteness of $(\vartheta_n)_{n\geqslant1}$ and Lemma~\ref{ns:lem:class-closure}
give $\rho\in\mathscr R_H^\infty(U,\mathbb R)$.
Since $\vartheta_n\geqslant0$, $\sum\limits_n\vartheta_n=1$, and all coefficients are positive, $\rho(x)>0$ for every $x\in U$. If
$\vartheta_n(x)\ne0$ or $D\vartheta_n(x)\ne0$, then
$W_n(x)>\frac{1}{4}$, and~\eqref{ns:eq:tau-support} gives $\tau(x)>\frac{1}{n}$.
Thus
\begin{equation*}
j_1(\rho;x)
 \leqslant\sum\limits_{\vartheta_n(x)\ne0\,\text{or}\,D\vartheta_n(x)\ne0}
       \frac{2^{-n}}{n}
 <\tau(x)\sum\limits_{n\geqslant1}2^{-n}=\tau(x).\qedhere
\end{equation*}
\end{proof}

\subsection[Holomorphic families and Gaussian approximation]{Absolutely summable holomorphic families and weighted Gaussian approximation}\label{ns:sec:tools}

\begin{lemma}\label{ns:aux:absolute-holomorphic}
Let $E,F$ be complex Banach spaces, let $V\subset E$ be open, and let $f_i\in\mathcal H(V,F)$ for $i\in I$. Suppose that every $a\in V$ has a ball $B(a,R)\subset V$ on which
\begin{equation}\label{ns:eq:local-absolute-bound}
 \sup\limits_{z\in B(a,R)}\sum\limits_{i\in I}\|f_i(z)\|\leqslant M<\infty,
\end{equation}
Then $\mathcal F=(f_i)_{i\in I}\in\mathcal H(V,\ell^1(I,F))$,
$f=\sum\limits_i f_i\in\mathcal H(V,F)$, and
\begin{equation}\label{ns:eq:directional-summation}
 Df(a)[h]=\sum\limits_iDf_i(a)[h],\qquad
 \sum\limits_i\|Df_i(a)[h]\|<\infty.
\end{equation}
\end{lemma}
\begin{proof}
Set $Y=\ell^1(I,F)$, and let $W\subset Y^*$ be the subspace spanned by finitely supported families of elements of $F^*$.
By the Hahn--Banach theorem, for every $u\in Y$,
\begin{equation*}
 \|u\|_Y=\sup\limits_{\substack{\varphi\in W\\\|\varphi\|\leqslant1}}|\varphi(u)|.
\end{equation*}
Condition~\eqref{ns:eq:local-absolute-bound} implies that $\mathcal F:V\to Y$ is locally bounded,
and each $\varphi\circ\mathcal F$ is a finite sum of holomorphic functions.
Apply \cite[Theorem~1.3, p.~781]{ArendtNikolski2000} on each complex line,
and then use the criterion for locally bounded G\^ateaux holomorphic mappings
in \cite[Proposition~8.6 and Theorem~8.7, p.~61]{Mujica1986}, to obtain $\mathcal F\in\mathcal H(V,Y)$.
If $\pi_i:Y\to F$ denotes the $i$th coordinate projection, then
\begin{equation*}
 \pi_iD\mathcal F(a)[h]=Df_i(a)[h],\qquad
 \sum\limits_i\|Df_i(a)[h]\|=\|D\mathcal F(a)[h]\|_Y<\infty.
\end{equation*}
Composition with the bounded linear summation operator
\begin{equation*}
 \Sigma:Y\to F,\qquad \Sigma((u_i)_{i\in I})=\sum\limits_i u_i
\end{equation*}
proves that $f=\Sigma\circ\mathcal F$ is holomorphic and yields~\eqref{ns:eq:directional-summation}.
\end{proof}

For a finite set $S\subset\Gamma$, write
\begin{equation}\label{ns:eq:coordinate-monomials}
 p_S(z)=\prod\limits_{\gamma\in S}z_\gamma^2,
 \qquad p_\varnothing=1.
\end{equation}

\begin{lemma}\label{ns:lem:pairs}
For $z\in H_{\mathbb C}$ and $x,h\in H$, one has
\begin{align}
 \sum\limits_{S\in\operatorname{Fin}(\Gamma)}|p_S(z)|
 &\leqslant e^{\|z\|^2},\label{ns:eq:pairs-value}\\
 \sum\limits_{S\in\operatorname{Fin}(\Gamma)}|Dp_S(x)[h]|
 &\leqslant2\|x\|e^{\|x\|^2}\|h\|.\label{ns:eq:pairs-derivative}
\end{align}
\end{lemma}
\begin{proof}
For every finite $S\subset\Gamma$,
\begin{equation*}
 Dp_S(x)[h]
 =2\sum\limits_{\gamma\in S}x_\gamma h_\gamma
       \prod\limits_{\delta\in S\setminus\{\gamma\}}x_\delta^2.
\end{equation*}
Thus, for every finite set $F\subset\Gamma$,
\begin{equation*}
 \sum\limits_{S\subset F}|p_S(z)|
 =\prod\limits_{\gamma\in F}(1+|z_\gamma|^2)\leqslant e^{\|z\|^2},
\end{equation*}
\begin{equation*}
 \begin{aligned}
 \sum\limits_{S\subset F}|Dp_S(x)[h]|
 &\leqslant2\sum\limits_{\gamma\in F}|x_\gamma h_\gamma|
        \sum\limits_{T\subset F\setminus\{\gamma\}}
                      \prod\limits_{\delta\in T}|x_\delta|^2\\
 &=2\sum\limits_{\gamma\in F}|x_\gamma h_\gamma|
        \prod\limits_{\delta\in F\setminus\{\gamma\}}(1+|x_\delta|^2)\\
 &\leqslant2e^{\|x\|^2}\sum\limits_{\gamma\in F}|x_\gamma h_\gamma|
 \leqslant2\|x\|e^{\|x\|^2}\|h\|.
 \end{aligned}
\end{equation*}
For every finite family $\mathcal A\subset\operatorname{Fin}(\Gamma)$,
the set $F=\bigcup\limits_{S\in\mathcal A}S$ is finite and contains
every member of $\mathcal A$.
Taking the supremum over such $F$ gives~\eqref{ns:eq:pairs-value} and~\eqref{ns:eq:pairs-derivative}, respectively.
\end{proof}

\begin{lemma}\label{ns:lem:weighted-sum-estimate}
Let $n\in\mathbb N$, $d\in\mathbb N_0$, and $S_i\in\operatorname{Fin}(\Gamma)$ satisfy
$\#\{i:S_i=S\}\leqslant n$. Let $Q_j\in C^1(U,\mathbb R)$ for $1\leqslant j\leqslant d$, and set
\begin{equation*}
 L=\sum\limits_{j=1}^d\|DQ_j\|_\infty<\infty,\qquad
 P_i(x)=(x_{S_i},q(x),Q_1(x),\ldots,Q_d(x)).
\end{equation*}
Let $\delta>0$. If $E_i\in C^1(\mathbb R^{|S_i|+d+1},Z)$ satisfies
\begin{equation}\label{ns:eq:weighted-coefficient-bound}
 \|E_i(w)\|+\|DE_i(w)\|\leqslant\delta e^{-2\|w\|},
\end{equation}
then $e_i=p_{S_i}(E_i\circ P_i)$ satisfies
\begin{equation}\label{ns:eq:weighted-sum-bound}
 \sup\limits_{x\in U}\sum\limits_i\|e_i(x)\|
 +\sup\limits_{x\in U}\sup\limits_{\|h\|\leqslant1}
       \sum\limits_i\|De_i(x)[h]\|\leqslant n(4+L)\delta.
\end{equation}
\end{lemma}
\begin{proof}
Put $R=\|x\|$. The product and chain rules give
\begin{equation*}
 \begin{aligned}
 De_i(x)[h]={}&Dp_{S_i}(x)[h]E_i(P_i(x))\\
 &+p_{S_i}(x)DE_i(P_i(x))[DP_i(x)[h]].
 \end{aligned}
\end{equation*}
Since $\|P_i(x)\|\geqslant R^2$ and
$\|DP_i(x)\|\leqslant1+2R+L$, the multiplicity bound and
Lemma~\ref{ns:lem:pairs} yield
\begin{equation*}
 \begin{aligned}
 \sum\limits_i\|De_i(x)[h]\|
 &\leqslant\delta e^{-2R^2}
   \left(\sum\limits_i|Dp_{S_i}(x)[h]|
       +(1+2R+L)\|h\|\sum\limits_i|p_{S_i}(x)|\right)\\
 &\leqslant n\delta(1+4R+L)e^{-R^2}\|h\|.
 \end{aligned}
\end{equation*}
Together with the value estimate, this gives
\begin{equation*}
 \sum\limits_i\|e_i(x)\|\leqslant n\delta e^{-R^2},\qquad
 \sum\limits_i\|De_i(x)[h]\|
 \leqslant n\delta(1+4R+L)e^{-R^2}\|h\|.
\end{equation*}
Moreover,
\begin{equation*}
 2R\leqslant1+R^2\leqslant e^{R^2},\qquad
 (1+4R+L)e^{-R^2}
 \leqslant(1+L)e^{-R^2}+2\leqslant3+L,
\end{equation*}
and hence~\eqref{ns:eq:weighted-sum-bound} follows.
\end{proof}

\begin{lemma}\label{ns:lem:Gaussian}
Let $d\in\mathbb N$, $b\in C_c^1(\mathbb R^d,Z)$, and $\kappa_1,\ldots,\kappa_d>0$. For $\lambda>0$, set
\begin{equation}\label{ns:eq:Gaussian}
 \begin{gathered}
 C_\kappa=\pi^{-\frac{d}{2}}\prod\limits_{j=1}^d\sqrt{\kappa_j},\\
 a_\lambda(w)=C_\kappa\lambda^{\frac{d}{2}}
   \int\limits_{\mathbb R^d}b(y)
   \exp\left(-\lambda\sum\limits_{j=1}^d\kappa_j(w_j-y_j)^2\right)\,\mathrm{d}y.
 \end{gathered}
\end{equation}
Then $a_\lambda:\mathbb C^d\to Z_{\mathbb C}$ is entire and takes values in $Z$ at real points.
At real points, we continue to write $a_\lambda,Da_\lambda$ for $a_\lambda|_{\mathbb R^d}$ and its Fr\'echet derivative.
For every $a>0$,
\begin{equation}\label{ns:eq:Gaussian-weighted}
 \lim\limits_{\lambda\to\infty}\sup\limits_{w\in\mathbb R^d}
 e^{a\|w\|}\bigl(\|a_\lambda(w)-b(w)\|
                  +\|Da_\lambda(w)-Db(w)\|\bigr)=0.
\end{equation}
Furthermore, if
$\operatorname{Re}\sum\limits_j\kappa_j(w_j-y_j)^2\geqslant\beta>0$
for every $y\in\operatorname{supp}b$, then
\begin{equation}\label{ns:eq:Gaussian-offreal}
 \|a_\lambda(w)\|\leqslant C_\kappa\|b\|_{L^1}\lambda^{\frac{d}{2}}e^{-\beta\lambda}.
\end{equation}
\end{lemma}
\begin{proof}
Inspired by \cite[Lemma~3 and its proof, pp.~3--4]{Johanis2015}, write
\begin{equation*}
 K_\lambda(t)=C_\kappa\lambda^{\frac{d}{2}}
       e^{-\lambda\sum\limits_j\kappa_jt_j^2},\qquad
 (v*K_\lambda)(x)=\int\limits_{\mathbb R^d}v(x-t)K_\lambda(t)\,\mathrm{d}t.
\end{equation*}
Take $v=b$ and $v=Db$ in turn. Since $\int\limits_{\mathbb R^d}K_\lambda=1$, for every $r>0$ one has
\begin{equation*}
 \|v*K_\lambda-v\|_\infty
 \leqslant\sup\limits_{\|t\|\leqslant r}
       \|v(\,\cdot-t)-v\|_\infty
       +2\|v\|_\infty\int\limits_{\|t\|>r}K_\lambda(t)\,\mathrm{d}t.
\end{equation*}
The mapping $v$ is uniformly continuous, and the substitution
$s=\sqrt\lambda\,t$ gives
\begin{equation*}
 \int\limits_{\|t\|>r}K_\lambda(t)\,\mathrm{d}t
 =C_\kappa\int\limits_{\|s\|>\sqrt\lambda r}
       e^{-\sum\limits_j\kappa_js_j^2}\,\mathrm{d}s\longrightarrow0.
\end{equation*}
For $h\in\mathbb R^d$, differentiation under the integral gives
\begin{equation*}
 Da_\lambda(x)[h]
 =\int\limits_{\mathbb R^d}Db(x-t)[h]K_\lambda(t)\,\mathrm{d}t.
\end{equation*}
Indeed, the difference quotients are bounded by
$\|Db\|_\infty\|h\|K_\lambda(t)$, an integrable function.
Letting first $\lambda\to\infty$ and then $r\downarrow0$ in the preceding
uniform estimate, we obtain
\begin{equation}\label{ns:eq:Gaussian-real-derivative}
 a_\lambda=b*K_\lambda,\qquad Da_\lambda=(Db)*K_\lambda,
 \qquad \lim\limits_{\lambda\to\infty}\|v*K_\lambda-v\|_\infty=0.
\end{equation}
To prove~\eqref{ns:eq:Gaussian-weighted}, choose $R\geqslant1$ with $\operatorname{supp}b\subset B(0,R)$,
and set $c=\frac{1}{4}\min\limits_j\kappa_j$.
For $v=b$ or $v=Db$, every $x\in\mathbb R^d$ with $\|x\|\geqslant2R$ satisfies $v(x)=0$ and
\begin{equation*}
 \|(v*K_\lambda)(x)\|
 \leqslant C_\kappa\|v\|_{L^1}\lambda^{\frac{d}{2}}e^{-c\lambda\|x\|^2},
\end{equation*}
where $\|v\|_{L^1}=\int\limits_{\mathbb R^d}\|v(y)\|\,\mathrm{d}y$.
For $\lambda\geqslant1$,
\begin{equation*}
 \sup\limits_{t\geqslant2R}
 e^{at}\lambda^{\frac{d}{2}}e^{-c\lambda t^2}
 \leqslant\left(\sup\limits_{t\geqslant0}e^{at-\frac{ct^2}{2}}\right)
          \lambda^{\frac{d}{2}}e^{-2c\lambda R^2}
 \longrightarrow0.
\end{equation*}
Combining this with uniform $C^1$ convergence on the ball $\{\|x\|\leqslant2R\}$ proves
\eqref{ns:eq:Gaussian-weighted}.

Fix a compact set $L\subset\mathbb C^d$.
For every multi-index $\alpha\in\mathbb N_0^d$, continuity and compactness give
\begin{equation*}
	M_{L,\alpha}
	:=C_\kappa\lambda^{\frac{d}{2}}
	\sup\limits_{\substack{w\in L\\
			y\in\operatorname{supp}b}}
	\left|
	\partial_w^\alpha
	\exp\left(
	-\lambda\sum\limits_{j=1}^d
	\kappa_j(w_j-y_j)^2
	\right)
	\right|
	<\infty.
\end{equation*}
Consequently, for every $w\in L$ and $y\in\mathbb R^d$,
\begin{equation*}
	\begin{aligned}
		&\left\|
		\partial_w^\alpha\left(
		C_\kappa\lambda^{\frac{d}{2}}b(y)
		\exp\left(
		-\lambda\sum\limits_{j=1}^d
		\kappa_j(w_j-y_j)^2
		\right)
		\right)
		\right\|\leqslant M_{L,\alpha}\|b(y)\|.
	\end{aligned}
\end{equation*}
Since $b$ has compact support, the function
$M_{L,\alpha}\|b(\cdot)\|$ belongs to
$L^1(\mathbb R^d)$.
Differentiation under the integral is therefore valid
for every multi-index $\alpha$.
In particular, for every $1\leqslant j\leqslant d$,

\begin{equation*}
 \frac{\partial a_\lambda}{\partial w_j}(w)
 =-2\lambda\kappa_j C_\kappa\lambda^{\frac{d}{2}}
   \int\limits_{\mathbb R^d}(w_j-y_j)b(y)
       e^{-\lambda\sum\limits_\ell\kappa_\ell(w_\ell-y_\ell)^2}\,\mathrm{d}y.
\end{equation*}
The same argument applies to every complex derivative, so $a_\lambda$ is entire.
If $\operatorname{Re}\sum\limits_j\kappa_j(w_j-y_j)^2
\geqslant\beta$ for every $y\in\operatorname{supp}b$, then
\begin{equation*}
 \begin{aligned}
 \|a_\lambda(w)\|
 &\leqslant C_\kappa\lambda^{\frac{d}{2}}
   \int\limits_{\mathbb R^d}\|b(y)\|
   e^{-\lambda\operatorname{Re}\sum\limits_j\kappa_j(w_j-y_j)^2}\,\mathrm{d}y\leqslant C_\kappa\|b\|_{L^1}\lambda^{\frac{d}{2}}e^{-\beta\lambda}.
 \end{aligned}
\end{equation*}
\end{proof}

\subsection[Uniform analytic approximation]{Uniform analytic approximation}\label{ns:sec:constant}

\begin{lemma}\label{ns:lem:entire-layer}
For each $W_n$ in Lemma~\ref{ns:lem:catalog} and every $\alpha>0$,
there is an entire function $H_n:H_{\mathbb C}\to\mathbb C$ such that
\begin{equation*}
 H_n(H)\subset\mathbb R,\qquad
 \sup\limits_{\|z\|\leqslant R}|H_n(z)|<\infty,
\end{equation*}
where the second assertion holds for every $R>0$, and whose restriction to $H$ satisfies
\begin{equation}\label{ns:eq:Hn-approx}
 \|H_n|_H-W_n\|_{\infty,H}
 +\|D(H_n|_H)-DW_n\|_{\infty,H}<\alpha.
\end{equation}
\end{lemma}
\begin{proof}
At real points, we continue to write $H_n,DH_n$ for $H_n|_H$ and its Fr\'echet derivative.
Recall that $\mathcal J_n=\{j\leqslant n:m_j\leqslant n\}$.
For $j\in\mathcal J_n$, set
\begin{equation}\label{ns:eq:divided-cutoff}
 r_j(s,t)=\frac{\beta_j(s,t)}{\prod\limits_{l=1}^{m_j}s_l^2}.
\end{equation}
Set $r_j=0$ wherever the denominator vanishes.
Since $\operatorname{supp}\beta_j$ is disjoint from each hyperplane $s_l=0$, this extension satisfies
$r_j\in C_c^\infty(\mathbb R^{m_j+1})$. Choose
$0<e<\frac{\alpha}{4n}$. Lemma~\ref{ns:lem:Gaussian} gives an entire function $a_j$ such that
\begin{equation*}
 |a_j(w)-r_j(w)|+\|Da_j(w)-Dr_j(w)\|<e e^{-2\|w\|}.
\end{equation*}
Set
\begin{equation}\label{ns:eq:entire-layer-definition}
 H_n(z)=\sum\limits_{i\in I_n}p_{S_i}(z)a_{j(i)}(P_{S_i}(z)).
\end{equation}
Fix $R>0$ and set
\begin{equation*}
 M_R=\max\limits_{j\in\mathcal J_n}
       \sup\limits_{\substack{s\in\mathbb C^{m_j},\,\|s\|\leqslant R\\
                             t\in\mathbb C,\,|t|\leqslant R^2}}
       |a_j(s,t)|<\infty,
\end{equation*}
with the convention $M_R=0$ if $\mathcal J_n=\varnothing$. Since
$\|z_{S_i}\|\leqslant\|z\|$ and $|q(z)|\leqslant\|z\|^2$, we obtain
\begin{equation*}
 \sup\limits_{\|z\|\leqslant R}\sum\limits_{i\in I_n}
 |p_{S_i}(z)a_{j(i)}(P_{S_i}(z))|
 \leqslant nM_Re^{R^2}.
\end{equation*}
Lemma~\ref{ns:aux:absolute-holomorphic} shows that $H_n$ is entire and bounded on bounded balls.
At real points, $p_{S_i}(r_{j(i)}\circ P_{S_i})=w_i$, so
$W_n=\sum\limits_{i\in I_n}p_{S_i}(r_{j(i)}\circ P_{S_i})$.
The supports of the summands are contained in the locally finite family $(K_i)_{i\in I_n}$; hence the sum may be differentiated as a finite sum near each point.
Differentiate $H_n$ in each direction using~\eqref{ns:eq:directional-summation},
and apply Lemma~\ref{ns:lem:weighted-sum-estimate} with $d=0$, to obtain
\begin{equation*}
\|H_n-W_n\|_{\infty,H}+\|DH_n-DW_n\|_{\infty,H}\leqslant4ne<\alpha.\qedhere
\end{equation*}
\end{proof}

\begin{theorem}\label{ns:thm:constant}
Let $u\in\mathscr R_H^1(U,Z)$ and $\varepsilon_*>0$. There is
$v\in C^\omega(U,Z)$ such that
\begin{equation}\label{ns:eq:constant-claim}
 \|v-u\|_{\infty,U}<\varepsilon_*,\qquad
 \|Dv-Du\|_{\infty,U}<\varepsilon_*.
\end{equation}
\end{theorem}
\begin{proof}
Inspired by \cite[proofs of Lemmas~5--6, pp.~74--78]{Whitney1934} and
\cite[proof of Theorem~1, pp.~4--6]{Johanis2015},
apply Lemma~\ref{ns:lem:catalog} with $\mathcal V=\{U\}$ and $\tau=1$, and retain its notation.

\noindent\firstproofstep{}\label{step:hilbert-approximation:1}\quad
Set
\begin{equation}\label{ns:eq:approximation-constants}
 \begin{gathered}
 C_c=\max\limits\{1,\|c\|_\infty,\|c'\|_\infty\},\qquad
 \mathcal T_n=K_\chi\sum\limits_{j\leqslant n}(E_j+1),\\
 \mathcal C_n=2L_nC_c(E_n+2)(1+\mathcal T_n),\qquad
 \alpha_n=\min\limits\left\{\frac{1}{32},
        \frac{\varepsilon_*2^{-n-4}}{1+\mathcal C_n}\right\}.
 \end{gathered}
\end{equation}
Use Lemma~\ref{ns:lem:entire-layer} to choose $H_n$ satisfying~\eqref{ns:eq:Hn-approx}
with $\alpha=\alpha_n$. Then
\begin{equation}\label{ns:eq:Hn-real-bounds}
 \|H_n\|_{\infty,H}\leqslant \widehat B_n+1,\qquad
 \|D(H_n|_H)\|_{\infty,H}\leqslant E_n+1.
\end{equation}
On the real space $H$, write
\begin{equation*}
 P_{n-1}=\prod\limits_{j<n}(1-\chi(H_j)),\qquad
 \theta_n=\chi(H_n)P_{n-1}.
\end{equation*}
For every $a\in U$, Lemma~\ref{ns:lem:catalog} provides an open neighborhood $V$ and an index $h$ such that $W_h|_V\geqslant1$.
Thus $H_h>\frac{31}{32}$ on $V$, so $\chi(H_h)|_V=1$ and $P_{n-1}|_V=0$ for $n>h$.
Applying Lemma~\ref{lem:successive-partition} with $b_n=\chi(H_n)$ shows that $(\theta_n)$ is locally finite and
\begin{equation}\label{ns:eq:theta-sum}
 \sum\limits_{n\geqslant1}\theta_n=1,\qquad
 \|DP_{n-1}\|_\infty
 \leqslant K_\chi\sum\limits_{j<n}(E_j+1)\leqslant\mathcal T_n.
\end{equation}
The product rule gives, for $h\in H$,
\begin{equation*}
 DP_{n-1}(x)[h]
 =-\sum\limits_{j<n}\chi'(H_j(x))DH_j(x)[h]
       \prod\limits_{\substack{\ell<n\\\ell\ne j}}(1-\chi(H_\ell(x))).
\end{equation*}
Since $0\leqslant\chi\leqslant1$, taking norms and using
\eqref{ns:eq:Hn-real-bounds} gives the derivative estimate in
\eqref{ns:eq:theta-sum}.
On the neighborhood $V$ above, $c(H_n)P_{n-1}W_nu$ vanishes identically for every $n>h$.
Consequently, the following expression is a finite sum of $C^1$ mappings near each point:
\begin{equation}\label{ns:eq:u0-definition}
 u_0=\sum\limits_{n\geqslant1}c(H_n)P_{n-1}W_nu.
\end{equation}
Since $tc(t)=\chi(t)$,
\begin{equation*}
 u_0-u=\sum\limits_{n\geqslant1}R_n^{(0)},\qquad
 R_n^{(0)}=c(H_n)P_{n-1}(W_n-H_n)u.
\end{equation*}
If $H_n(x)\leqslant\frac{1}{4}$, then $j_1(R_n^{(0)};x)=0$; otherwise,
$W_n(x)>\frac{1}{4}-\frac{1}{32}>0$, and hence $j_1(u;x)\leqslant L_n$.
On the set $H_n>\frac{1}{4}$, we have
\begin{equation*}
 \begin{gathered}
 j_1(c\circ H_n;x)\leqslant C_c(E_n+2),\qquad
 j_1(P_{n-1};x)\leqslant1+\mathcal T_n,\\
 j_1(W_n-H_n;x)<\alpha_n,\qquad j_1(u;x)\leqslant L_n.
 \end{gathered}
\end{equation*}
Thus the first-order product estimate~\eqref{eq:jetproduct} gives
\begin{equation*}
 j_1(R_n^{(0)};x)
 \leqslant L_nC_c(E_n+2)(1+\mathcal T_n)\alpha_n.
\end{equation*}
Therefore
\begin{equation}\label{ns:eq:u0-u-bound}
 \|u_0-u\|_{\infty,U}+\|Du_0-Du\|_{\infty,U}
 \leqslant\sum\limits_{n\geqslant1}\mathcal C_n\alpha_n
 <\frac{\varepsilon_*}{16}.
\end{equation}

\noindent\proofstep{}\label{step:hilbert-approximation:2}\quad
Set
\begin{equation*}
 k_j=\widehat B_j+2,\qquad
 \zeta(t)=1-\sigma\left(\frac{t^2-1}{3}\right),\qquad
 \zeta_n(t_1,\ldots,t_n)=\prod\limits_{j=1}^n\zeta\left(\frac{t_j}{k_j}\right).
\end{equation*}
Then $\zeta=1$ on $[-1,1]$, $\zeta=0$ on $\mathbb R\setminus(-2,2)$, and
$\zeta_n(H_1(x),\ldots,H_n(x))=1$ at every real point $x$.
For $i\in I_n$, set
\begin{equation}\label{ns:eq:Bni}
 \begin{split}
 B_{n,i}(s,t_0,t_1,\ldots,t_n)
 ={}&r_{j(i)}(s,t_0)b_i(s,t_0)c(t_n)\cdot\prod\limits_{j<n}(1-\chi(t_j))\,
          \zeta_n(t_1,\ldots,t_n).
 \end{split}
\end{equation}
Here $r_{j(i)}b_i$ is extended by zero outside $O_i$. Since
$\operatorname{supp}r_{j(i)}\Subset O_i$,
\begin{equation*}
 B_{n,i}\in C_c^1(\mathbb R^{m_i+n+1},Z),\qquad
 \operatorname{supp}B_{n,i}\subset
 \operatorname{supp}\beta_{j(i)}\times\prod\limits_{j=1}^n[-2k_j,2k_j].
\end{equation*}
In particular, on this support,
\begin{equation}\label{ns:eq:literal-support}
 |s_l|\geqslant a_n,\qquad
 t_j\leqslant\frac{1}{2}.
\end{equation}
The first estimate holds for every $1\leqslant l\leqslant m_i$, and the second for every $j<n$.
Define the entire mapping $P_{n,S}$ and the locally finite sum $U_n$ by
\begin{align}
 P_{n,S}(z)&=(z_S,q(z),H_1(z),\ldots,H_n(z)),\label{ns:eq:PnS}\\
 U_n(x)&=\sum\limits_{i\in I_n}p_{S_i}(x)B_{n,i}(P_{n,S_i}(x)).\label{ns:eq:Un}
\end{align}
The support of the $i$th summand in~\eqref{ns:eq:Un} is contained in $K_i$. Since $(K_i)_{i\in I_n}$ is locally finite, $U_n$ is a $C^1$ mapping.
For $x\in K_i$, the identities
$p_{S_i}(x)r_{j(i)}(P_{S_i}(x))=w_i(x)$,
$b_i(P_{S_i}(x))=u(x)$, and $\zeta_n(H_1(x),\ldots,H_n(x))=1$ give
\begin{equation*}
 p_{S_i}(x)B_{n,i}(P_{n,S_i}(x))
 =w_i(x)u(x)c(H_n(x))P_{n-1}(x).
\end{equation*}
Both sides vanish on $U\setminus K_i$. Hence, for $x\in U$,
\begin{equation*}
 U_n(x)=c(H_n(x))P_{n-1}(x)u(x)\sum\limits_{i\in I_n}w_i(x)
       =c(H_n(x))P_{n-1}(x)W_n(x)u(x),
\end{equation*}
and therefore
\begin{equation}\label{ns:eq:Un-identity}
 U_n\in C^1(H,Z),\qquad U_n|_U=c(H_n)P_{n-1}W_nu,
 \qquad u_0=\sum\limits_{n\geqslant1}U_n|_U.
\end{equation}

\noindent\proofstep{}\label{step:hilbert-approximation:3}\quad
Inspired by \cite[equation~(9) and the proof of Claim~1.5, pp.~87, 89--92]{CepedelloHajek2001} and
\cite[proof of Lemma~3(5), pp.~134--139]{AFK2012}, set
\begin{equation}\label{ns:eq:cj}
 c_j=\frac{2^{-j}}{1+\sup\limits_{\|z\|\leqslant j}|H_j(z)|^2}>0.
\end{equation}
For every $R>0$ and every integer $J\geqslant R$,
\begin{equation}\label{ns:eq:weighted-entire-tail}
 \sup\limits_{\|z\|\leqslant R}\sum\limits_{j>J}c_j|H_j(z)|^2\leqslant2^{-J}.
\end{equation}
In dimension $m+n+1$, choose the weights
$\kappa^{(n,m)}=(1,\ldots,1,c_1,\ldots,c_n)$, where the first $m+1$ entries equal $1$.
Write
\begin{align}
 \Lambda_{n,i}&=1+\pi^{-\frac{m_i+n+1}{2}}
       \left(\prod\limits_{j=1}^n\sqrt{c_j}\right)\|B_{n,i}\|_{L^1},\label{ns:eq:Lambda}\\
 \mathcal D_n&=\sum\limits_{j=1}^n\|D(H_j|_H)\|_{\infty,H},\qquad
 \widehat C_n=n(4+\mathcal D_n),\qquad
 e_n=\frac{\varepsilon_*2^{-n-4}}{1+\widehat C_n}.\label{ns:eq:en}
\end{align}
For each $(n,i)$, approximate $B_{n,i}$ using~\eqref{ns:eq:Gaussian}.
By~\eqref{ns:eq:Gaussian-weighted}, we may choose an integer $\lambda_{n,i}\geqslant1$
such that the corresponding entire mapping $A_{n,i}$ satisfies
\begin{align}
 &\sup\limits_{w\in\mathbb R^{m_i+n+1}}e^{2\|w\|}
   \bigl(\|A_{n,i}(w)-B_{n,i}(w)\|
      +\|DA_{n,i}(w)-DB_{n,i}(w)\|\bigr)<e_n,\label{ns:eq:scale-real}\\
 &n\Lambda_{n,i}(1+\lambda_{n,i})^{\frac{2n+1}{2}}
       \leqslant2^{-n}e^{\sqrt{\lambda_{n,i}}}.\label{ns:eq:scale-growth}
\end{align}
Condition~\eqref{ns:eq:scale-growth} can be imposed simultaneously because
\begin{equation*}
 \lim\limits_{\lambda\to\infty}
 \frac{\log(n\Lambda_{n,i})+n\log2
       +\frac{2n+1}{2}\log(1+\lambda)}{\sqrt\lambda}=0.
\end{equation*}

\noindent\proofstep{}\label{step:hilbert-approximation:4}\quad
Set
\begin{equation}\label{ns:eq:analytic-summands-and-S}
 \begin{gathered}
 F_{n,i}(z)=p_{S_i}(z)A_{n,i}(P_{n,S_i}(z)),\\
 \mathcal S(z)=\|\operatorname{Im}z\|^2+|\operatorname{Im}q(z)|^2
       +\sum\limits_{j\geqslant1}c_j|\operatorname{Im}H_j(z)|^2.
 \end{gathered}
\end{equation}
By~\eqref{ns:eq:weighted-entire-tail} and
$|\operatorname{Im}H_j(z)|\leqslant|H_j(z)|$, the series defining $\mathcal S$ converges uniformly on every bounded ball.
Its terms are continuous, so $\mathcal S$ is continuous. Since $q$ and each $H_j$ are real-valued at real points, $\mathcal S|_H=0$.
Fix $x^0\in U$, and choose $h$ such that $W_h\geqslant1$ on a neighborhood of $x^0$. Then
$H_h(x^0)>\frac{31}{32}$. Set
\begin{equation}\label{ns:eq:phase-constant}
 \beta=\frac{1}{32}\min\limits\{c_h,a_1^2,\ldots,a_h^2\}>0,
\end{equation}
and, for $1\leqslant n\leqslant h$, let
\begin{equation*}
 \begin{gathered}
 J_n=\{\gamma:|x^0_\gamma|\geqslant\frac{a_n}{2}\},\qquad
 \mathcal E_n=\{i\in I_n:S_i\subset J_n\},\\
 \mathcal E=\{(n,i):n\leqslant h,\ i\in\mathcal E_n\}.
 \end{gathered}
\end{equation*}
Since $|J_n|\leqslant\frac{4\|x^0\|^2}{a_n^2}$ and
$|\mathcal E_n|\leqslant n2^{|J_n|}$, the set $\mathcal E$ is finite.
Fix $R>\|x^0\|+1$. The continuity of $H_h$ and $\mathcal S$ at $x^0$,
together with $H_h(x^0)>\frac{31}{32}$ and $\mathcal S(x^0)=0$, allows us to choose a sufficiently small complex ball $V$ centered at $x^0$
such that $V\subset B(0,R)$, $V\cap H\subset U$, and, for every $z\in V$,
\begin{equation}\label{ns:eq:common-complex-neighborhood}
 \operatorname{Re}H_h(z)>\frac{3}{4},\qquad
 \mathcal S(z)<\beta,\qquad
 \|z-x^0\|<\frac{1}{4}\min\limits_{n\leqslant h}a_n.
\end{equation}
Since $\mathcal E$ is finite and each $F_{n,i}$ is continuous at $x^0$, we may shrink the same complex ball $V$ further
so that, for every $(n,i)\in\mathcal E$,
\begin{equation*}
 \sup\limits_{z\in V}\|F_{n,i}(z)\|
 \leqslant1+\|F_{n,i}(x^0)\|.
\end{equation*}
For $y=(s,t_0,\ldots,t_n)\in\operatorname{supp}B_{n,i}$, write
\begin{equation}\label{ns:eq:gaussian-phase}
 \mathcal Q_{n,S_i}(z;y)=\sum\limits_{\gamma\in S_i}(z_\gamma-s_\gamma)^2
 +(q(z)-t_0)^2+\sum\limits_{j=1}^nc_j(H_j(z)-t_j)^2.
\end{equation}
The identity
$\operatorname{Re}\bigl((w-t)^2\bigr)=(\operatorname{Re}w-t)^2-(\operatorname{Im}w)^2$
valid for every $t\in\mathbb R$, and the definition of $\mathcal S$ give
\begin{equation*}
 \begin{aligned}
 \operatorname{Re}\mathcal Q_{n,S_i}(z;y)
 \geqslant{}&\sum\limits_{\gamma\in S_i}(\operatorname{Re}z_\gamma-s_\gamma)^2
       +(\operatorname{Re}q(z)-t_0)^2\\
 &+\sum\limits_{j=1}^nc_j(\operatorname{Re}H_j(z)-t_j)^2
       -\mathcal S(z).
 \end{aligned}
\end{equation*}
If $n>h$, then $t_h\leqslant\frac{1}{2}$, and hence
\begin{equation*}
 \operatorname{Re}\mathcal Q_{n,S_i}(z;y)
 \geqslant c_h(\operatorname{Re}H_h(z)-t_h)^2-\mathcal S(z)
 >\frac{c_h}{16}-\beta\geqslant\beta.
\end{equation*}
If $n\leqslant h$ and $i\notin\mathcal E_n$, choose $\gamma\in S_i\setminus J_n$. Then
\begin{equation*}
 |\operatorname{Re}z_\gamma-s_\gamma|
 \geqslant|s_\gamma|-|x^0_\gamma|-|\operatorname{Re}z_\gamma-x^0_\gamma|
 >\frac{a_n}{4},
\end{equation*}
Keeping the term $(\operatorname{Re}z_\gamma-s_\gamma)^2$
in the preceding lower bound and using $\mathcal S(z)<\beta$,
we obtain
$\operatorname{Re}\mathcal Q_{n,S_i}(z;y)>\frac{a_n^2}{16}-\beta\geqslant\beta$.
Thus, for every $(n,i)\notin\mathcal E$, $z\in V$, and $y\in\operatorname{supp}B_{n,i}$,
\begin{equation}\label{ns:eq:gaussian-phase-lower-bound}
 \operatorname{Re}\mathcal Q_{n,S_i}(z;y)\geqslant\beta.
\end{equation}

\noindent\proofstep{}\label{step:hilbert-approximation:5}\quad
By~\eqref{ns:eq:Gaussian-offreal}, \eqref{ns:eq:scale-growth}, \eqref{ns:eq:gaussian-phase-lower-bound}, and $m_i\leqslant n$,
for $(n,i)\notin\mathcal E$ we have
\begin{equation*}
 \begin{aligned}
 \|F_{n,i}(z)\|
 &\leqslant |p_{S_i}(z)|\Lambda_{n,i}
          \lambda_{n,i}^{\frac{2n+1}{2}}e^{-\beta\lambda_{n,i}}\\
 &\leqslant\frac{2^{-n}}{n} e^{\frac{1}{4\beta}}|p_{S_i}(z)|,
 \end{aligned}
\end{equation*}
where we have used
\begin{equation*}
 \sqrt\lambda-\beta\lambda
 =\frac{1}{4\beta}-\beta\left(\sqrt\lambda-\frac{1}{2\beta}\right)^2
 \leqslant\frac{1}{4\beta}.
\end{equation*}
The multiplicity bound in~\eqref{ns:eq:tau-support} and~\eqref{ns:eq:pairs-value} give
\begin{equation}\label{ns:eq:local-holomorphic-family-bound}
 \sup\limits_{z\in V}\sum\limits_{n\geqslant1}\sum\limits_{i\in I_n}\|F_{n,i}(z)\|
 \leqslant\sum\limits_{(n,i)\in\mathcal E}\sup\limits_{z\in V}\|F_{n,i}(z)\|
       +e^{R^2+\frac{1}{4\beta}}<\infty.
\end{equation}
By Lemma~\ref{ns:aux:absolute-holomorphic}, the mappings
\begin{equation}\label{ns:eq:G-formal}
 G=\sum\limits_{n\geqslant1}\sum\limits_{i\in I_n}F_{n,i},\qquad
 G_n=\sum\limits_{i\in I_n}F_{n,i}
\end{equation}
are holomorphic on $V$. Denote the neighborhood constructed at $a\in U$
by $V_a$ and the local sum by $G^{(a)}$. For $z\in V_a\cap V_b$,
absolute summability gives
\begin{equation*}
 G^{(a)}(z)=\sum\limits_{n\geqslant1}\sum\limits_{i\in I_n}F_{n,i}(z)
          =G^{(b)}(z).
\end{equation*}
Thus $\mathcal V=\bigcup\limits_{a\in U}V_a$ is a complex neighborhood of $U$,
and $G(z)=G^{(a)}(z)$ for $z\in V_a$ defines $G\in\mathcal H(\mathcal V,Z_{\mathbb C})$.
Since $F_{n,i}(x)\in Z$ for every $x\in U$ and $Z$ is closed in $Z_{\mathbb C}$,
we have $G(U)\subset Z$, and hence $v=G|_U\in C^\omega(U,Z)$.

\noindent\proofstep{}\label{step:hilbert-approximation:6}\quad
We continue to denote the restrictions of $G_n$ and $U_n$ to $U$ by $G_n$ and $U_n$.
On $\mathbb R^{m_i+n+1}$, set $E_{n,i}=A_{n,i}|_{\mathbb R^{m_i+n+1}}-B_{n,i}$.
By Lemma~\ref{ns:aux:absolute-holomorphic}, $G$ and each $G_n$ may be differentiated term by term in each direction.
The supports of the summands in~\eqref{ns:eq:Un} are contained in the locally finite family $(K_i)_{i\in I_n}$, so $U_n$ may also be differentiated term by term.
Thus
\begin{equation*}
 G_n-U_n=\sum\limits_{i\in I_n}p_{S_i}(E_{n,i}\circ P_{n,S_i}).
\end{equation*}
Apply Lemma~\ref{ns:lem:weighted-sum-estimate} to these summands with
$Q_j=H_j|_U$ for $j=1,\ldots,n$, using
$\sum\limits_{j=1}^n\|DQ_j\|_{\infty,U}\leqslant\mathcal D_n$.
Equation~\eqref{ns:eq:scale-real} supplies the exponentially weighted coefficient estimate required in that lemma. Consequently,
\begin{equation}\label{ns:eq:Gn-Un-bound}
 \|G_n-U_n\|_{\infty,U}+\|D(G_n-U_n)\|_{\infty,U}
 \leqslant \widehat C_ne_n.
\end{equation}
By Step~\ref{step:hilbert-approximation:1} and~\eqref{ns:eq:Un-identity},
$u_0=\sum\limits_nU_n$ is a finite sum near each real point.
Fix $x\in U$ and $h\in H$. Grouping the absolutely summable family defining $G$ and its derivatives in the direction $h$ according to $n$ yields
\begin{equation*}
 \begin{aligned}
 v(x)-u_0(x)&=\sum\limits_{n\geqslant1}(G_n(x)-U_n(x)),\\
 (Dv(x)-Du_0(x))[h]&=\sum\limits_{n\geqslant1}D(G_n-U_n)(x)[h].
 \end{aligned}
\end{equation*}
Therefore
\begin{equation}\label{ns:eq:v-u0}
 \|v-u_0\|_{\infty,U}+\|Dv-Du_0\|_{\infty,U}
 \leqslant\sum\limits_{n\geqslant1}\widehat C_ne_n
 <\frac{\varepsilon_*}{16}.
\end{equation}
Combining~\eqref{ns:eq:v-u0} and~\eqref{ns:eq:u0-u-bound}, we obtain
\begin{equation*}
\|v-u\|_{\infty,U}+\|Dv-Du\|_{\infty,U}
 <\frac{\varepsilon_*}{8}<\varepsilon_*.\qedhere
\end{equation*}
\end{proof}

\subsection[Fine analytic approximation]{Fine analytic approximation}\label{ns:sec:fine}

\begin{theorem}\label{ns:thm:fine}
Let $u\in\mathscr R_H^1(U,Z)$ and $\epsilon\in C(U,(0,\infty))$.
There is $v\in C^\omega(U,Z)$ such that, for every $x\in U$,
\begin{equation*}
 \|v(x)-u(x)\|<\epsilon(x),\qquad
 \|Dv(x)-Du(x)\|<\epsilon(x).
\end{equation*}
\end{theorem}

\begin{proof}[Proof of Theorem~\ref{ns:thm:fine}]
Inspired by \cite[Chapter~7, Lemma~34, p.~405]{HajekJohanis2014} and
\cite[proof of Theorem~2, p.~222]{Kurzweil1954}, set $M(x)=1+j_1(u;x)$ and
$\tau(x)=\min\limits\{\frac{1}{8},\frac{\epsilon(x)}{4M(x)}\}$.
By Lemma~\ref{ns:lem:small-positive}, choose
$\rho\in\mathscr R_H^\infty(U,\mathbb R)$ such that
\begin{equation}\label{ns:eq:rho-fine}
 0<\rho(x),\qquad j_1(\rho;x)<\tau(x).
\end{equation}
Since $\rho>0$, apply Lemma~\ref{ns:lem:class-closure} to $(\rho,u)$ and
the smooth mapping $(t,z)\mapsto(t^{-1},t^{-1}z)$ on $(0,\infty)\times Z$ to obtain
$(\rho^{-1},\rho^{-1}u)\in\mathscr R_H^1(U,\mathbb R\oplus Z)$.
Apply Theorem~\ref{ns:thm:constant} to this mapping with $\varepsilon_*=\frac{1}{4}$ to obtain an analytic mapping $(A,V)$.
Set $d=A-\rho^{-1}$ and $e=V-\rho^{-1}u$. Since
$\max\limits_{j=0,1}\max\limits\{\|D^jd\|_{\infty,U},\|D^je\|_{\infty,U}\}<\frac{1}{4}$,
we have
\begin{equation}\label{ns:eq:de-jet-bounds}
 j_1(d;x)<\frac{1}{2},\qquad j_1(e;x)<\frac{1}{2}.
\end{equation}
Since $\rho<\frac{1}{8}$, we have $A>\rho^{-1}-\frac{1}{4}>0$; set $B=A^{-1}$ and $v=BV$.
The identity $B=\rho-\rho dB$ gives
\begin{equation*}
 j_1(B;x)\leqslant j_1(\rho;x)
       +j_1(\rho;x)j_1(d;x)j_1(B;x),
\end{equation*}
Since $j_1(\rho;x)j_1(d;x)<\frac{1}{16}$,
\begin{equation*}
 j_1(B;x)
 \leqslant\frac{j_1(\rho;x)}{1-j_1(\rho;x)j_1(d;x)}
 \leqslant2j_1(\rho;x).
\end{equation*}
Since $A=\rho^{-1}+d$ and $V=\rho^{-1}u+e$,
\begin{equation*}
 \begin{aligned}
 B&=\frac{\rho}{1+\rho d},\qquad DB=-A^{-2}DA,\\
 v-u&=B(e-du),\\
 D(v-u)(x)[h]
 &=DB(x)[h](e(x)-d(x)u(x))\\
 &\quad+B(x)\bigl(De(x)[h]-Dd(x)[h]u(x)-d(x)Du(x)[h]\bigr).
 \end{aligned}
\end{equation*}
Applying \eqref{eq:jetproduct}, we obtain
\begin{align*}
j_1(v-u;x)
 &\leqslant2j_1(\rho;x)\bigl(j_1(e;x)+j_1(d;x)j_1(u;x)\bigr)\\
 &<j_1(\rho;x)(1+j_1(u;x))
 <M(x)\tau(x)\leqslant\frac{\epsilon(x)}{4}.\qedhere
\end{align*}
\end{proof}

\section[\texorpdfstring{$C^0$}{C0} Approximation]{Fine \texorpdfstring{$C^0$}{C0} Real-Analytic Approximation on Spaces Admitting a Separating Polynomial}\label{sec:nonseparable-c0}

\subsection[Polynomial embeddings into Hilbert spaces]{Polynomial topological embeddings into Hilbert spaces and their characterization}\label{subsec:general-polynomial-embedding}\label{subsec:hilbert-embedding-characterization}

\begin{theorem}\label{thm:general-polynomial-hilbert-embedding}\label{prop:hilbert-embedding-separator}\label{cor:polynomial-hilbert-characterization}
For a real Banach space $X$, the following conditions are equivalent:
\begin{enumerate}[label=\textnormal{(\roman*)},ref=(\roman*)]
\item\label{cond:hilbert-characterization-separator} $X$ admits a separating polynomial;
\item\label{cond:hilbert-characterization-polynomial} $X$ admits a polynomial topological embedding into a real Hilbert space;
\item\label{cond:hilbert-characterization-analytic} $X$ admits a real-analytic topological embedding into a real Hilbert space.
\end{enumerate}
Here a topological embedding $J:X\to H$ means that $J$ is continuous and injective and that $J^{-1}:J(X)\to X$ is continuous.
The mapping $J$ in condition~\ref{cond:hilbert-characterization-polynomial} admits a complexification to an entire mapping
$X_{\mathbb C}\to H_{\mathbb C}$.
\end{theorem}

\begin{lemma}\label{lem:homogeneous-separator}
If a real Banach space $X$ admits a separating polynomial, then there are an even integer $d\geqslant2$, a continuous symmetric
$d$-linear form $A$, and a constant $B\geqslant1$ such that the functions
\begin{equation*}
 Q(x)=A(x,\ldots,x),\qquad \mu(x)=Q(x)^{\frac{1}{d}}
\end{equation*}
satisfy the following estimates, the first two for every $x\in X$ and the third for every $x\in X\setminus\{0\}$:
\begin{equation}\label{eq:separator-gauge}
 \|x\|^d\leqslant Q(x),\qquad
 \|x\|\leqslant\mu(x)\leqslant B\|x\|,\qquad
 \|D\mu(x)\|\leqslant B.
\end{equation}
Moreover, $\mu$ is $B$-Lipschitz and
$\mu|_{X\setminus\{0\}}\in C^\infty(X\setminus\{0\})$.
\end{lemma}
\begin{proof}
If $X=\{0\}$, take $d=2$, $A=0$, $Q=0$, and $B=1$.
If $X\ne\{0\}$, \cite[proof of $(iii)\Rightarrow(iv)$ in Corollary~2, p.~420]{FabianEtAl1989}
provides a positive homogeneous separating polynomial of even degree. After multiplication by a positive constant, write it as
$Q(x)=A(x^d)\geqslant\|x\|^d$.
Then $\mu=Q^{\frac{1}{d}}$ is smooth on $X\setminus\{0\}$, and, for every $x\in X\setminus\{0\}$,
\begin{equation*}
 \begin{gathered}
 D\mu(x)[h]=\frac{1}{d}Q(x)^{\frac{1}{d}-1}DQ(x)[h]
          =Q(x)^{\frac{1}{d}-1}A(x^{d-1},h),\\
 \|D\mu(x)\|\leqslant\|A\|.
 \end{gathered}
\end{equation*}
Take
\begin{equation*}
 B=\max\limits\left\{1,\|A\|,
                  \sup\limits_{\|x\|=1}Q(x)^{\frac{1}{d}}\right\}.
\end{equation*}
If $0\notin[x,y]$, integration along the segment gives
\begin{equation*}
 |\mu(y)-\mu(x)|
 =\left|\int\limits_0^1D\mu(x+t(y-x))[y-x] \,\mathrm{d}t\right|
 \leqslant B\|y-x\|.
\end{equation*}
If $0\in[x,y]$, then
\begin{equation*}
 |\mu(x)-\mu(y)|\leqslant\mu(x)+\mu(y)
 \leqslant B(\|x\|+\|y\|)=B\|x-y\|.
\end{equation*}
Thus~\eqref{eq:separator-gauge} holds.
\end{proof}

A real Banach space admitting a separating polynomial is superreflexive
\cite[discussion preceding Theorem~2.2 in Section~2, p.~151]{GonzaloJaramillo1997}.
For a superreflexive real Banach space $X$, fix an equivalent uniformly convex norm of power type
\cite[Theorems~3.1(a) and~3.2, pp.~337, 340]{Pisier1975},
\cite[Theorem~10.1, p.~382; Theorem~11.37, p.~454]{Pisier2016}.
By \cite[proof of Lemma~2, p.~133]{JohnTorunczykZizler1981},
there are an exponent $p_0\in(1,\infty)$ and a constant $C_{\mathrm{PRI}}\geqslant1$
such that every nonseparable closed subspace of $X$ admits a projectional resolution
satisfying~\eqref{eq:general-pri-norms}--\eqref{eq:general-pri-continuity}
and the $\ell^{p_0}$ estimate in the proof of Lemma~\ref{lem:general-pri-increments},
with the same constant $C_{\mathrm{PRI}}$.
Set $p=\max\limits\{2,p_0\}$.

Let $E$ be a nonseparable reflexive Banach space and let $\kappa=\operatorname{dens}E$.
Since $E$ is weakly compactly generated, choose a projectional resolution
of the identity with projections of norm at most one
\cite[Lemma~6, pp.~44--45]{AmirLindenstrauss1968}.
Reindex it as in \cite[Section~1, p.~514]{Kalenda2020}
to obtain a well-ordered set $A$
satisfying~\eqref{eq:well-order-small-parts}, with $|A|=\kappa$,
and pairwise commuting bounded linear projections $(P_\alpha^-,P_\alpha^+)_{\alpha\in A}$ such that, for every $\alpha\in A$,
\begin{equation}\label{eq:general-pri-norms}
 \|P_\alpha^-\|\leqslant1,\qquad \|P_\alpha^+\|\leqslant1,
\end{equation}
\begin{equation}\label{eq:general-pri-commutation}
 \begin{gathered}
 P_\alpha^-P_\alpha^+=P_\alpha^-,\\
 P_\alpha^+P_\beta^-=P_\alpha^+\quad\text{for every }\beta\in A\text{ with }\alpha<\beta,
 \end{gathered}
\end{equation}
\begin{equation}\label{eq:general-pri-density}
 \operatorname{dens}(P_\alpha^+E)
 \leqslant\max\limits\{\aleph_0,|A_{\leqslant\alpha}|\}<\kappa,
\end{equation}
and
\begin{equation}\label{eq:general-pri-continuity}
 P_\alpha^-E
 =\overline{\bigcup\limits_{\beta<\alpha}P_\beta^+E},\qquad
 E=\overline{\bigcup\limits_{\alpha\in A}P_\alpha^+E}.
\end{equation}
When $\alpha$ is the least element of $A$, the right-hand side of the first identity is interpreted as $\{0\}$, so $P_\alpha^-=0$.
All closures are norm closures, and each closed subspace carries the restriction of the fixed norm above.

\begin{lemma}\label{lem:general-pri-increments}
Let $E$ be a nonseparable closed subspace of $X$, equipped with the projectional resolution of the identity above. Set, for every $\alpha\in A$,
\begin{equation}\label{eq:general-pri-increments}
 R_\alpha=P_\alpha^+-P_\alpha^-,\qquad K=C_{\mathrm{PRI}}^p.
\end{equation}
Then, for every $x\in E$,
\begin{align}
 &\sum\limits_{\alpha\in A}\|R_\alpha x\|^p
       \leqslant K\|x\|^p,\label{eq:general-pri-sum}\\
 &x\in\overline{\operatorname{span}}
       \{R_\alpha x:\alpha\in A\}.
       \label{eq:general-pri-span}
\end{align}
\end{lemma}
\begin{proof}
Since $p\geqslant p_0$, the choice of the uniform constant $C_{\mathrm{PRI}}$ gives, for every $x\in E$,
\begin{equation*}
 \left(\sum\limits_{\alpha\in A}\|R_\alpha x\|^p\right)^{\frac{1}{p}}
 \leqslant
 \left(\sum\limits_{\alpha\in A}\|R_\alpha x\|^{p_0}\right)^{\frac{1}{p_0}}
 \leqslant C_{\mathrm{PRI}}\|x\|.
\end{equation*}
Raising this inequality to the $p$th power yields~\eqref{eq:general-pri-sum}.
Fix $x\in E$ and set
\begin{equation*}
 M_x=\overline{\operatorname{span}}\{R_\alpha x:\alpha\in A\}.
\end{equation*}
If $\alpha$ is not the least element, then~\eqref{eq:general-pri-continuity} implies that, for every $\eta>0$,
there are $\beta<\alpha$ and $y\in P_\beta^+E$ such that
$\|P_\alpha^-x-y\|<\eta$.
Using~\eqref{eq:general-pri-commutation}, $P_\beta^+y=y$, and $\|P_\beta^+\|\leqslant1$, we obtain
\begin{equation*}
 \begin{aligned}
 \|P_\beta^+x-P_\alpha^-x\|
 &=\|P_\beta^+(P_\alpha^-x-y)-(P_\alpha^-x-y)\|\\
 &\leqslant2\|P_\alpha^-x-y\|<2\eta.
 \end{aligned}
\end{equation*}
Therefore
\begin{equation*}
 P_\alpha^-x\in\overline{\{P_\beta^+x:\beta<\alpha\}}.
\end{equation*}
If there is an $\alpha\in A$ with $P_\alpha^+x\notin M_x$, choose the least such $\alpha$.
Every preceding $P_\beta^+x$ belongs to the closed subspace $M_x$, so the preceding inclusion gives $P_\alpha^-x\in M_x$;
when $\alpha$ is the least element, the same conclusion follows from $P_\alpha^-=0$.
Hence
\begin{equation*}
 P_\alpha^+x=P_\alpha^-x+R_\alpha x\in M_x,
\end{equation*}
a contradiction. Thus every $P_\alpha^+x$ belongs to $M_x$.
Finally, the second identity in~\eqref{eq:general-pri-continuity} implies that, for every $\eta>0$, there are
$\alpha\in A$ and $y\in P_\alpha^+E$ with $\|x-y\|<\eta$.
Consequently,
\begin{equation*}
 \|P_\alpha^+x-x\|
 =\|P_\alpha^+(x-y)-(x-y)\|\leqslant2\|x-y\|<2\eta.
\end{equation*}
Since $M_x$ is closed, $x\in M_x$, which is~\eqref{eq:general-pri-span}.
\end{proof}

\begin{lemma}\label{lem:general-support-map}
Let $X$ be a superreflexive real Banach space with the fixed norm and exponent $p$ above.
There are an index set $I$, vectors $v_i\in X$ with $\|v_i\|\leqslant1$, and a bounded linear mapping
$T:X\to\ell^p(I;\mathbb R)$ such that, for every $x\in X$,
\begin{equation}\label{eq:general-support-map-properties}
 \|T\|\leqslant1,\qquad
 x\in\overline{\operatorname{span}}
          \{v_i:(Tx)_i\ne0\}.
\end{equation}
In particular, $T$ is injective.
\end{lemma}
\begin{proof}
Inspired by \cite[proof of Theorem~10, p.~703]{Johanis2017},
fix the norm and exponent $p$ above and argue by induction on the cardinal $\operatorname{dens}E$ over all closed subspaces $E$ of $X$.
For $E=\{0\}$, take $I_E=\varnothing$ and $T_E=0$.

\firstproofstep{}\label{step:density-coordinates:1}
Suppose that $E\ne\{0\}$ is separable. Choose a dense sequence $(v_n)_{n\geqslant1}$ in its unit sphere.
By the Hahn--Banach theorem, choose $f_j\in E^*$ such that
\begin{equation*}
 \|f_j\|=1,\qquad f_j(v_j)=1.
\end{equation*}
For $x\ne0$, choose $j$ with $\|v_j-\frac{x}{\|x\|}\|<\frac{1}{2}$. Then
\begin{equation*}
 f_j\left(\frac{x}{\|x\|}\right)
 \geqslant1-\left\|v_j-\frac{x}{\|x\|}\right\|>\frac{1}{2}.
\end{equation*}
Set $I_E=\mathbb N^2$, associate the vector $v_{j,n}=v_n$ with the coordinate $(j,n)$, and define
\begin{equation}\label{eq:general-support-map-separable}
 (T_Ex)_{j,n}=2^{-j-n}f_j(x).
\end{equation}
Then $T_E$ is linear and
\begin{equation*}
 \|T_Ex\|_p^p
 \leqslant\left(\sum\limits_{j\geqslant1}2^{-pj}\right)^2\|x\|^p
 \leqslant\|x\|^p.
\end{equation*}
If $x\ne0$, choose $j$ with $f_j(x)\ne0$. Then, for every $n\geqslant1$,
\begin{equation*}
 (T_Ex)_{j,n}\ne0,
\end{equation*}
and consequently
\begin{equation*}
 x\in E=\overline{\operatorname{span}}\{v_{j,n}:n\geqslant1\}
 \subset\overline{\operatorname{span}}
       \{v_{k,n}:(T_Ex)_{k,n}\ne0\}.
\end{equation*}

\proofstep{}\label{step:density-coordinates:2}
Suppose that $E$ is nonseparable, put $\kappa=\operatorname{dens}E$, and assume that the assertion holds for every closed subspace of density less than
$\kappa$. Choose the projectional resolution in~\eqref{eq:general-pri-norms}--\eqref{eq:general-pri-continuity},
retaining the notation $R_\alpha,K$ from~\eqref{eq:general-pri-increments}.
For every $\alpha\in A$, $R_\alpha$ is a bounded linear projection, so its range $E_\alpha=R_\alpha E$ is closed.
By~\eqref{eq:general-pri-commutation}, $E_\alpha\subset P_\alpha^+E$.
Together with~\eqref{eq:general-pri-density}, this gives
\begin{equation}\label{eq:general-increment-density}
 \operatorname{dens}E_\alpha
 \leqslant\operatorname{dens}(P_\alpha^+E)
 \leqslant\max\limits\{\aleph_0,|A_{\leqslant\alpha}|\}<\kappa.
\end{equation}
The inductive hypothesis provides
\begin{equation*}
 T_\alpha:E_\alpha\to\ell^p(I_\alpha;\mathbb R),\qquad
 \|T_\alpha\|\leqslant1,\qquad v_{\alpha,i}\in E_\alpha,\quad\|v_{\alpha,i}\|\leqslant1,
\end{equation*}
such that, for every $y\in E_\alpha$,
\begin{equation}\label{eq:general-child-support}
 y\in\overline{\operatorname{span}}
       \{v_{\alpha,i}:(T_\alpha y)_i\ne0\}.
\end{equation}
Inspired by \cite[Section~3.2, equation~(3.3), pp.~28--29]{KaniaKochanek2016} and
\cite[proof of the stated theorem, p.~139]{BenyaminiStarbird1976}, take disjoint copies of the sets $I_\alpha$ and define
\begin{equation}\label{eq:general-support-map-induction}
 I_E=\{(\alpha,i):\alpha\in A,\ i\in I_\alpha\},\qquad
 (T_Ex)_{\alpha,i}=K^{-\frac{1}{p}}(T_\alpha R_\alpha x)_i.
\end{equation}
By Lemma~\ref{lem:general-pri-increments},
\begin{equation*}
 \|T_Ex\|_p^p
 =K^{-1}\sum\limits_{\alpha\in A}\|T_\alpha R_\alpha x\|_p^p
 \leqslant K^{-1}\sum\limits_{\alpha\in A}\|R_\alpha x\|^p
 \leqslant\|x\|^p.
\end{equation*}
Thus $T_E:E\to\ell^p(I_E)$ is bounded and linear. Fix $x\in E$ and set
\begin{equation*}
 N_x=\overline{\operatorname{span}}
       \{v_{\alpha,i}:(T_Ex)_{\alpha,i}\ne0\}.
\end{equation*}
Equations~\eqref{eq:general-child-support} and~\eqref{eq:general-support-map-induction}
give $R_\alpha x\in N_x$ for every $\alpha$; equation~\eqref{eq:general-pri-span} then yields $x\in N_x$.
This completes the induction on density. Finally, if $Tx=0$, then
$x\in\overline{\operatorname{span}}\varnothing=\{0\}$.
\end{proof}

\begin{proof}[Proof of Theorem~\ref{thm:general-polynomial-hilbert-embedding}]
If $X=\{0\}$, take $H=\mathbb R$ and $J(0)=0$; all three conditions hold. Henceforth suppose that $X\ne\{0\}$. We first prove that
condition~\ref{cond:hilbert-characterization-separator} implies condition~\ref{cond:hilbert-characterization-polynomial}.
Passing to an equivalent norm does not affect the property of being a topological embedding. By superreflexivity, fix the equivalent norm $\|\cdot\|$ used in Lemma~\ref{lem:general-pri-increments}.
Denote the original norm by $\|\cdot\|_0$, and choose $b_*>0$ such that $\|x\|\leqslant b_*\|x\|_0$.
Lemma~\ref{lem:homogeneous-separator} provides an even integer $d\geqslant2$ and
$Q_0(x)=\mathcal Q_0(x^d)\geqslant\|x\|_0^d$.
Set $\mathcal Q=b_*^d\mathcal Q_0$, $Q=b_*^dQ_0$, and $C=\max\limits\{1,\|\mathcal Q\|\}$,
where the multilinear norm is taken with respect to $\|\cdot\|$. Then
\begin{equation}\label{eq:general-separating-form}
 Q(x)=\mathcal Q(x^d),\qquad
 \|x\|^d\leqslant Q(x)\leqslant C\|x\|^d.
\end{equation}
Choose $I,T,(v_i)_{i\in I}$ as in Lemma~\ref{lem:general-support-map}.

\firstproofstep{}\label{step:polynomial-embedding:1}
Choose $r\in\mathbb N$ with $2r\geqslant p$, and set
\begin{equation}\label{eq:general-hilbert-power}
 A(x)=((Tx)_i^r)_{i\in I}\in\ell^2(I;\mathbb R).
\end{equation}
For $u^1,\ldots,u^r\in\ell^p(I;\mathbb R)$, H\"older's inequality gives
\begin{equation}\label{eq:general-coordinate-product}
 \left\|\left(\prod\limits_{j=1}^r u_i^j\right)_{i\in I}\right\|_2
 \leqslant\prod\limits_{j=1}^r\|u^j\|_{2r}
 \leqslant\prod\limits_{j=1}^r\|u^j\|_p.
\end{equation}
Thus $M(u^1,\ldots,u^r)=(\prod\limits_{j=1}^ru_i^j)_{i\in I}$ is a bounded $r$-linear mapping,
and $A(x)=M(Tx,\ldots,Tx)$, so $A$ is a continuous $r$-homogeneous polynomial.
Moreover, $(Tx)_i^r\ne0$ if and only if $(Tx)_i\ne0$. Equation~\eqref{eq:general-support-map-properties} therefore gives
\begin{equation}\label{eq:general-hilbert-power-properties}
 \|A(x)\|_2\leqslant\|Tx\|_p^r\leqslant\|x\|^r,\qquad
 x\in\overline{\operatorname{span}}\{v_i:A_i(x)\ne0\}.
\end{equation}
For $1\leqslant k\leqslant d-1$, define the multilinear mapping
\begin{equation*}
 \begin{aligned}
 \mathcal B_k(u^1,\ldots,u^k;z_1,\ldots,z_{d-k})
 =\biggl(&\left(\prod\limits_{\nu=1}^ku^\nu_{i_\nu}\right)\\
 &\cdot\mathcal Q(z_1,\ldots,z_{d-k},v_{i_1},\ldots,v_{i_k})
                                       \biggr)_{(i_1,\ldots,i_k)\in I^k},
 \end{aligned}
\end{equation*}
where $u^\nu\in\ell^2(I)$ and $z_j\in X$. Since $\|v_i\|\leqslant1$, repeated summation of nonnegative terms gives
\begin{equation*}
 \begin{aligned}
 \|\mathcal B_k(u^1,\ldots,u^k;z_1,\ldots,z_{d-k})\|_2^2
 &\leqslant\|\mathcal Q\|^2\prod\limits_{j=1}^{d-k}\|z_j\|^2
   \sum\limits_{(i_1,\ldots,i_k)\in I^k}
       \prod\limits_{\nu=1}^k|u^\nu_{i_\nu}|^2\\
 &=\|\mathcal Q\|^2\prod\limits_{j=1}^{d-k}\|z_j\|^2
       \prod\limits_{\nu=1}^k\|u^\nu\|_2^2.
 \end{aligned}
\end{equation*}
Taking square roots, we obtain
\begin{equation}\label{eq:general-multilinear-bound}
 \|\mathcal B_k(u^1,\ldots,u^k;z_1,\ldots,z_{d-k})\|_2
 \leqslant\|\mathcal Q\|\prod\limits_{\nu=1}^k\|u^\nu\|_2
                       \prod\limits_{j=1}^{d-k}\|z_j\|.
\end{equation}
Consequently, $B_k(x)=\mathcal B_k(A(x),\ldots,A(x);x,\ldots,x)$
is a continuous $(kr+d-k)$-homogeneous polynomial with coordinates
\begin{equation}\label{eq:general-mixed-polynomial}
 B_k(x)_{i_1,\ldots,i_k}
 =\left(\prod\limits_{\nu=1}^kA_{i_\nu}(x)\right)
   \mathcal Q(x^{d-k},v_{i_1},\ldots,v_{i_k}).
\end{equation}
Equation~\eqref{eq:general-multilinear-bound} yields
\begin{equation}\label{eq:general-mixed-polynomial-bound}
 \|B_k(x)\|_2\leqslant\|\mathcal Q\|\|A(x)\|_2^k\|x\|^{d-k}.
\end{equation}
Set
\begin{equation}\label{eq:general-polynomial-embedding}
 \begin{gathered}
 H=\mathbb R\oplus_2\ell^2(I;\mathbb R)
       \oplus_2\left(\bigoplus\limits_{k=1}^{d-1}
                         \ell^2(I^k;\mathbb R)\right)_2,\\
 J(x)=(Q(x),A(x),B_1(x),\ldots,B_{d-1}(x)).
 \end{gathered}
\end{equation}
Then $H$ is a real Hilbert space and $J$ is a continuous polynomial of degree at most $(d-1)r+1$. Its componentwise complexification is entire.

\proofstep{}\label{step:polynomial-embedding:2}
Given a finite set $F\subset I$ and real numbers $(c_i)_{i\in F}$, set
\begin{equation}\label{eq:general-finite-center}
 b=\sum\limits_{i\in F}c_iv_i.
\end{equation}
Write the ambient coordinates as $y=(s,z,w_1,\ldots,w_{d-1})$, and define the open set
\begin{equation}\label{eq:general-recovery-domain}
 \Omega_F=\{(s,z,w_1,\ldots,w_{d-1})\in H:z_i\ne0\text{ for every }i\in F\}.
\end{equation}
For the fixed representation in~\eqref{eq:general-finite-center}, define
\begin{equation}\label{eq:general-recovery-function}
 \begin{split}
 \mathcal R_{F,c}(s,z,w_1,\ldots,w_{d-1})
 ={}&s+Q(b)\\
 &+\sum\limits_{k=1}^{d-1}(-1)^k\binom{d}{k}
    \sum\limits_{(i_1,\ldots,i_k)\in F^k}
      c_{i_1}\cdots c_{i_k}
      \frac{(w_k)_{i_1,\ldots,i_k}}{z_{i_1}\cdots z_{i_k}}.
 \end{split}
\end{equation}
Since $F$ is finite and $z_i\ne0$ on $\Omega_F$ for every $i\in F$,
we have $\mathcal R_{F,c}\in C^\omega(\Omega_F,\mathbb R)$.
If $J(x)\in\Omega_F$, then~\eqref{eq:general-mixed-polynomial} gives
\begin{equation*}
 \frac{B_k(x)_{i_1,\ldots,i_k}}
      {A_{i_1}(x)\cdots A_{i_k}(x)}
 =\mathcal Q(x^{d-k},v_{i_1},\ldots,v_{i_k}).
\end{equation*}
Since $d$ is even, multilinear expansion yields
\begin{equation}\label{eq:general-recovery-identity}
 \begin{aligned}
 \mathcal R_{F,c}(J(x))
 &=Q(x)+Q(b)+\sum\limits_{k=1}^{d-1}(-1)^k\binom{d}{k}
                                    \mathcal Q(x^{d-k},b^k)\\
 &=Q(x-b).
 \end{aligned}
\end{equation}

\proofstep{}\label{step:polynomial-embedding:3}
At the origin, $J(0)=0$ and
\begin{equation}\label{eq:general-embedding-origin}
 \|x\|^d\leqslant Q(x)\leqslant\|J(x)-J(0)\|_H.
\end{equation}
Fix $a\ne0$ and $\eta>0$.
Inspired by \cite[proof of Theorem~2, pp.~696--698]{Johanis2017}, use~\eqref{eq:general-hilbert-power-properties}
to choose a nonempty finite set $F\subset\{i:A_i(a)\ne0\}$ and real coefficients $(c_i)$ such that
the vector $b$ in~\eqref{eq:general-finite-center} satisfies
\begin{equation}\label{eq:general-center-approximation}
 \|a-b\|<\frac{\eta}{4C^{\frac{1}{d}}}\leqslant\frac{\eta}{4},\qquad
 Q(a-b)<\left(\frac{\eta}{4}\right)^d.
\end{equation}
Since $J(a)\in\Omega_F$, equation~\eqref{eq:general-recovery-identity} and
the continuity of $\mathcal R_{F,c}$ give an open neighborhood $J(a)\in W\subset\Omega_F$ such that, for every $y\in W$,
\begin{equation*}
 \mathcal R_{F,c}(y)<\left(\frac{\eta}{2}\right)^d.
\end{equation*}
If $J(x)\in W$, then
\begin{equation*}
 \|x-a\|\leqslant\|x-b\|+\|b-a\|
 \leqslant Q(x-b)^{\frac{1}{d}}+\|b-a\|
 <\frac{\eta}{2}+\frac{\eta}{4}<\eta.
\end{equation*}
For $a=0$, equation~\eqref{eq:general-embedding-origin} allows us to take $W=B_H(0,\eta^d)$.
Thus, for every $a\in X$ and $\eta>0$,
there is an open set $W\subset H$ such that
\begin{equation}\label{eq:general-inverse-neighborhood}
 J(a)\in W,\qquad J^{-1}(W)\subset B_X(a,\eta).
\end{equation}
If $J(x)=J(a)$, then~\eqref{eq:general-inverse-neighborhood} gives
\begin{equation*}
 x\in\bigcap\limits_{\eta>0}B_X(a,\eta)=\{a\}.
\end{equation*}
Thus $J$ is injective; equation~\eqref{eq:general-inverse-neighborhood} also shows that $J^{-1}:J(X)\to X$ is continuous.

\proofstep{}\label{step:polynomial-embedding:4}
Every continuous polynomial is real-analytic, so condition~\ref{cond:hilbert-characterization-polynomial} implies condition~\ref{cond:hilbert-characterization-analytic}.
It remains to deduce condition~\ref{cond:hilbert-characterization-separator} from condition~\ref{cond:hilbert-characterization-analytic}.
In this step, $\|\cdot\|$ denotes the original norm on $X$. Let $J:X\to H$ be a real-analytic topological embedding,
and set $p(x)=\|J(x)-J(0)\|^2$.
On a neighborhood of $0$, there is an expansion into homogeneous polynomials
\begin{equation*}
 p(x)=\sum\limits_{n=1}^{\infty}p_n(x),\qquad
 \sum\limits_{n=1}^{\infty}\sup\limits_{\|x\|\leqslant r}|p_n(x)|<\infty
\end{equation*}
for sufficiently small $r>0$. Since $J^{-1}$ is continuous at $J(0)$, there is $\eta>0$ such that
\begin{equation*}
 \|J(x)-J(0)\|<\eta\quad\Longrightarrow\quad\|x\|<r.
\end{equation*}
Choose $N\in\mathbb N$ such that
\begin{equation*}
 q_N=\sum\limits_{n=1}^Np_n,\qquad
 \sup\limits_{\|x\|\leqslant r}|q_N(x)-p(x)|<\frac{\eta^2}{2}.
\end{equation*}
Since $\|x\|=r$ implies $p(x)\geqslant\eta^2$, the continuous polynomial $Q(x)=q_N(rx)$ satisfies
\begin{equation*}
Q(0)=0,\qquad
 \inf\limits_{\|x\|=1}Q(x)
 \geqslant\eta^2-\frac{\eta^2}{2}=\frac{\eta^2}{2}>0.\qedhere
\end{equation*}
\end{proof}

\subsection[Fine approximation of continuous maps]{Fine analytic approximation of continuous mappings}\label{subsec:general-c0-transfer}

\begin{proposition}\label{prop:general-c0-transfer}
Suppose that a real Banach space $X$ admits a real-analytic topological embedding $J:X\to H$ into a real Hilbert space $H$.
For every subset $A\subset X$, every real Banach space $Y$,
every $f\in C(A,Y)$, and every $\varepsilon\in C(A,(0,\infty))$,
there are an open set $U\subset X$ containing $A$ and $g\in C^\omega(U,Y)$
such that, for every $x\in A$,
\begin{equation}\label{eq:general-c0-transfer-conclusion}
 \|g(x)-f(x)\|<\varepsilon(x).
\end{equation}
\end{proposition}
\begin{proof}
If $A=\varnothing$, take $U=X$ and $g=0$; if $X=\{0\}$ and $A\ne\varnothing$, take $U=X$ and $g=f$.
Henceforth suppose that $A\ne\varnothing$ and $X\ne\{0\}$. Set
\begin{equation*}
 B=J(A),\qquad \widehat f=f\circ(J|_A)^{-1},\qquad
 \widehat\varepsilon=\varepsilon\circ(J|_A)^{-1}.
\end{equation*}
Since $J|_A$ is a homeomorphism, both $\widehat f$ and $\widehat\varepsilon$ are continuous.
Apply the open-cover and convex-combination construction of Lemma~\ref{lem:arbitrary-subset-neighborhood} to $B,\widehat f,\widehat\varepsilon$,
using the subordinate partition of unity from Proposition~\ref{ns:prop:partition}.
More explicitly, on the resulting open set $V\supset B$, choose a locally finite partition of unity
$(\theta_i)_{i\in I_0}\subset\mathscr R_H^\infty(V,\mathbb R)$,
let $b_i\in B$ be the point corresponding to the associated member of the cover, and define
\begin{equation}\label{eq:general-c0-auxiliary-maps}
 u(y)=\sum\limits_{i\in I_0}\theta_i(y)\widehat f(b_i),\qquad
 \delta(y)=\frac{1}{8}\sum\limits_{i\in I_0}\theta_i(y)\widehat\varepsilon(b_i).
\end{equation}
Lemma~\ref{ns:lem:class-closure}, together with $\theta_i\geqslant0$ and $\sum\limits_i\theta_i=1$, gives
\begin{equation*}
 u\in\mathscr R_H^\infty(V,Y),\qquad
 \delta\in\mathscr R_H^\infty(V,\mathbb R),\qquad\delta>0.
\end{equation*}
For $x\in A$, set $I_x=\{i:\theta_i(J(x))>0\}$.
Subordination and~\eqref{eq:arbitrary-subset-local-control} give, for $i\in I_x$,
\begin{equation*}
 \|\widehat f(b_i)-f(x)\|<\frac{\widehat\varepsilon(b_i)}{8},\qquad
 \widehat\varepsilon(b_i)<2\varepsilon(x).
\end{equation*}
Since $I_x$ is finite and $\sum\limits_{i\in I_x}\theta_i(J(x))=1$,
\begin{equation*}
 \begin{aligned}
 \|u(J(x))-f(x)\|
 &\leqslant\sum\limits_{i\in I_x}\theta_i(J(x))
             \|\widehat f(b_i)-f(x)\|\\
 &<\frac{1}{8}\sum\limits_{i\in I_x}\theta_i(J(x))\widehat\varepsilon(b_i)
 =\delta(J(x))<\frac{\varepsilon(x)}{4}.
 \end{aligned}
\end{equation*}
Thus, for every $x\in A$,
\begin{equation}\label{eq:general-c0-approximation}
 \|u(J(x))-f(x)\|<\delta(J(x))<\frac{\varepsilon(x)}{4}.
\end{equation}
By Theorem~\ref{ns:thm:fine}, choose $v\in C^\omega(V,Y)$ such that $\|v(y)-u(y)\|<\delta(y)$.
Set $U=J^{-1}(V)$ and $g=v\circ J|_U$. Then $U$ is open and contains $A$, $g\in C^\omega(U,Y)$, and, for every $x\in A$,
\begin{align*}
\|g(x)-f(x)\|
 &\leqslant\|v(J(x))-u(J(x))\|+\|u(J(x))-f(x)\|\\
 &<2\delta(J(x))<\varepsilon(x).\qedhere
\end{align*}
\end{proof}

\begin{proof}[Proof of Theorem~\ref{thm:main-nonseparable}]
By Theorem~\ref{thm:general-polynomial-hilbert-embedding}, $X$ admits a polynomial
topological embedding into a Hilbert space. Proposition~\ref{prop:general-c0-transfer} then gives
\eqref{eq:main-nonseparable-approximation}.
\end{proof}
\section[\texorpdfstring{$C^1$}{C1} Approximation with Decompositions]{Fine \texorpdfstring{$C^1$}{C1} Approximation on Spaces with Unconditional Finite-Dimensional Decompositions}\label{sec:ufdd}

\subsection[A sufficient condition for \texorpdfstring{$C^1$}{C1} approximation]{A sufficient condition for fine \texorpdfstring{$C^1$}{C1} approximation}\label{sec:transfer}

Let $X$ be a real Banach space, let $H=\ell^2(\Gamma_H;\mathbb R)$, and let
$J:X\to H$ be a real-analytic topological embedding.
The class $\mathscr R_H^k$ and the quadratic polynomial $q$
are given by Definition~\ref{ns:def:class} and~\eqref{ns:eq:qPS}, respectively; for $y\in H$, we have $q(y)=\|y\|_H^2$.
When $H=\{0\}$, the corresponding assertions follow immediately by taking constant mappings.

In what follows, the domains are assumed to be nonempty.

\begin{definition}\label{def:liftable-class}
Let $U\subset X$ be open, let $Z$ be a real Banach space, let $k\in\{1,\infty\}$, and let $u:U\to Z$.
If, for every $a\in U$, there exist an open neighborhood $W\subset H$ of $J(a)$ and
a mapping $\widetilde u\in\mathscr R_H^k(W,Z)$ such that
\begin{equation}\label{eq:local-lift-representation}
 J^{-1}(W)\subset U,\qquad
 u|_{J^{-1}(W)}=\widetilde u\circ J|_{J^{-1}(W)},
\end{equation}
then we write $u\in\mathscr F_J^k(U,Z)$ and say that $u$ admits
\emph{local $\mathscr R_H^k$-extensions along $J$}.
Equivalently, the mapping $u\circ(J|_U)^{-1}:J(U)\to Z$ satisfies
\begin{equation*}
 \widetilde u|_{W\cap J(U)}
 =\bigl(u\circ(J|_U)^{-1}\bigr)|_{W\cap J(U)}.
\end{equation*}
The term local extension below always includes the requirement
$\widetilde u\in\mathscr R_H^k(W,Z)$.
\end{definition}

\begin{lemma}\label{lem:liftable-closure}
We have $\mathscr F_J^k(U,Z)\subset C^k(U,Z)$. The classes $\mathscr F_J^k$ are stable under restriction to open subsets,
finite products, continuous multilinear operations, postcomposition with $C^k$ mappings, and locally finite sums.
If $V\subset U$ is open, $\psi\in\mathscr F_J^k(U,\mathbb R)$,
$v\in\mathscr F_J^k(V,Z)$, and $\operatorname{supp}_U\psi\subset V$,
then the extension of $\psi v$ by zero to $U\setminus V$ belongs to $\mathscr F_J^k(U,Z)$.
\end{lemma}
\begin{proof}
Fix $a\in U$ and choose $W,\widetilde u$ as in
\eqref{eq:local-lift-representation}. On $J^{-1}(W)$, the chain rule gives
\begin{equation*}
 u=\widetilde u\circ J,\qquad
 Du(x)[h]=D\widetilde u(J(x))[DJ(x)[h]].
\end{equation*}
Since $\widetilde u\in C^k(W,Z)$ and $J\in C^\omega(X,H)$,
we have $u\in C^k(U,Z)$.
If $U'\subset U$ is open and $a\in U'$, continuity of
$J^{-1}:J(X)\to X$ gives an open set $O\subset H$ such that
$J(a)\in O$ and $J^{-1}(O)\subset U'$.
Then
\begin{equation*}
 W'=W\cap O,\qquad J^{-1}(W')\subset U',\qquad
 u|_{J^{-1}(W')}=(\widetilde u|_{W'})\circ J|_{J^{-1}(W')}
\end{equation*}
is the required local extension representation for $u|_{U'}$.

For finitely many $u_j\in\mathscr F_J^k(U,Z_j)$, choose
$W_j,\widetilde u_j$ at $a$ and set $W=\bigcap\limits_{j=1}^NW_j$.
On $W$, define
\begin{equation*}
 \widetilde U(y)=(\widetilde u_1(y),\ldots,\widetilde u_N(y)),\qquad
 \widetilde v(y)=M(\widetilde u_1(y),\ldots,\widetilde u_N(y)),
\end{equation*}
where $M:Z_1\times\cdots\times Z_N\to Z$ is continuous and multilinear.
Lemma~\ref{ns:lem:class-closure} places these mappings in their respective
classes $\mathscr R_H^k$, and, on $J^{-1}(W)$,
\begin{equation*}
 \widetilde U\circ J=(u_1,\ldots,u_N),\qquad
 \widetilde v\circ J=M(u_1,\ldots,u_N).
\end{equation*}
If $\Phi\in C^k(O_Z,Z')$ and $u(U)\subset O_Z$, where $O_Z\subset Z$ is open,
choose a local extension $\widetilde u\in\mathscr R_H^k(W,Z)$ of $u$ at $a$.
Set $W'=W\cap\widetilde u^{-1}(O_Z)$ and define
\begin{equation*}
 \widetilde v=\Phi\circ\widetilde u|_{W'}\in\mathscr R_H^k(W',Z'),\qquad
 \widetilde v\circ J=\Phi\circ u\quad\hbox{on }J^{-1}(W').
\end{equation*}

For a locally finite family $(u_i)_{i\in I}$, choose an open neighborhood
$V_0$ of $a$ and a finite set $I_a$ such that $u_i|_{V_0}=0$ for $i\notin I_a$.
Choose an open neighborhood $W_0$ of $J(a)$ with $J^{-1}(W_0)\subset V_0$.
For the finitely many local extensions $\widetilde u_i\in\mathscr R_H^k(W_i,Z)$,
set
\begin{equation*}
 W=W_0\cap\bigcap\limits_{i\in I_a}W_i,\qquad
 \widetilde s=\sum\limits_{i\in I_a}\widetilde u_i|_W\in\mathscr R_H^k(W,Z).
\end{equation*}
For $x\in J^{-1}(W)$,
\begin{equation*}
 \widetilde s(J(x))=\sum\limits_{i\in I_a}u_i(x)
                  =\sum\limits_{i\in I}u_i(x).
\end{equation*}

Finally, let
\begin{equation*}
 w(x)=\begin{cases}\psi(x)v(x),&x\in V,\\0,&x\in U\setminus V.\end{cases}
\end{equation*}
The open sets $V$ and $U\setminus\operatorname{supp}_U\psi$ cover $U$.
At $a\in V$, the product construction above gives a common neighborhood
$W$ and
\begin{equation*}
 \widetilde w(y)=\widetilde\psi(y)\widetilde v(y),\qquad
 w|_{J^{-1}(W)}=\widetilde w\circ J|_{J^{-1}(W)}.
\end{equation*}
At $a\in U\setminus\operatorname{supp}_U\psi$, choose $W$ with
$J^{-1}(W)\subset U\setminus\operatorname{supp}_U\psi$ and take
$\widetilde w=0$. These representations prove $w\in\mathscr F_J^k(U,Z)$.
\end{proof}

\begin{lemma}\label{lem:global-lift}
If $u\in\mathscr F_J^k(U,Z)$, then there exist an open set $V\subset H$ and
a mapping $\widetilde u\in\mathscr R_H^k(V,Z)$ such that $J(U)\subset V$ and
$u=\widetilde u\circ J$ on $U$.
\end{lemma}
\begin{proof}
For each $a\in U$, choose $W_a,u_a$ as in \eqref{eq:local-lift-representation}, and set
$V=\bigcup\limits_{a\in U}W_a$. By Proposition~\ref{ns:prop:partition}, there is a locally finite partition of unity
\begin{equation*}
 0\leqslant\theta_i\in\mathscr R_H^\infty(V,\mathbb R),\qquad
 \sum\limits_i\theta_i=1,\qquad \operatorname{supp}_V\theta_i\subset W_{a_i}.
\end{equation*}
Define
\begin{equation*}
 v_i(y)=\begin{cases}\theta_i(y)u_{a_i}(y),&y\in W_{a_i},\\
                     0,&y\in V\setminus W_{a_i}.
        \end{cases}
\end{equation*}
Lemma~\ref{ns:lem:class-closure} gives $v_i\in\mathscr R_H^k(V,Z)$, and
$\operatorname{supp}_Vv_i\subset\operatorname{supp}_V\theta_i$ proves local
finiteness. Also, $J^{-1}(V)=U$, since $V=\bigcup\limits_{a\in U}W_a$
and $J^{-1}(W_a)\subset U$. Hence
$\widetilde u=\sum\limits_i v_i\in\mathscr R_H^k(V,Z)$.
If $\theta_i(J(x))\ne0$, then $J(x)\in W_{a_i}$, so
$u_{a_i}(J(x))=u(x)$. Thus
\begin{equation*}
\widetilde u(J(x))=\sum\limits_i\theta_i(J(x))u(x)=u(x).\qedhere
\end{equation*}
\end{proof}

\begin{proposition}\label{prop:analytic-liftable}
Every $u\in\mathscr F_J^1(U,Z)$ admits real-analytic approximation in the fine $C^1$ topology.
\end{proposition}
\begin{proof}
Given $\eta\in C(U,(0,\infty))$, use Lemma~\ref{lem:global-lift} to write
$u=\widetilde u\circ J$, $\widetilde u\in\mathscr R_H^1(V,Z)$.
The function
\begin{equation*}
 a(J(x))=\frac{\eta(x)}{1+\|DJ(x)\|}
\end{equation*}
is the composition of the positive continuous function $\frac{\eta}{1+\|DJ\|}$ with $J^{-1}|_{J(U)}$,
and is therefore continuous and positive on $J(U)$.
Apply Lemma~\ref{lem:arbitrary-subset-neighborhood} with
$X=H$, $A=J(U)$, $f=0$, and $\varepsilon=2a$.
Intersecting the resulting neighborhood with $V$, we obtain
$J(U)\subset V_0\subset V$ and $\delta\in C(V_0,(0,\infty))$ such that
\begin{equation*}
 \delta(J(x))<\frac{\eta(x)}{2(1+\|DJ(x)\|)}.
\end{equation*}
By Theorem~\ref{ns:thm:fine}, choose $\widetilde v\in C^\omega(V_0,Z)$ such that
\begin{equation*}
 \max\limits\{\|\widetilde v-\widetilde u\|,
                  \|D\widetilde v-D\widetilde u\|\}(y)<\delta(y).
\end{equation*}
Set $v=\widetilde v\circ J|_U$. Since $J(U)\subset V_0$, this defines
$v\in C^\omega(U,Z)$, and
\begin{equation*}
 Dv(x)-Du(x)
 =\bigl(D\widetilde v(J(x))-D\widetilde u(J(x))\bigr)\circ DJ(x).
\end{equation*}
Consequently,
\begin{align*}
\|v(x)-u(x)\|&<\delta(J(x))<\eta(x),\\
 \|Dv(x)-Du(x)\|
 &\leqslant\delta(J(x))\|DJ(x)\|<\eta(x).\qedhere
\end{align*}
\end{proof}

\begin{lemma}\label{lem:liftable-localization}
Every open cover of $U$ admits a subordinate locally finite partition of unity
$(\varphi_{n,j})_{n\geqslant1,\,j\in I_n}\subset\mathscr F_J^\infty(U,\mathbb R)$
and constants $D_n<\infty$ such that, for every $x\in U$,
\begin{equation}\label{eq:liftable-layer-bounds}
 0\leqslant\varphi_{n,j}\leqslant1,\qquad
 \sum\limits_{n,j}\varphi_{n,j}=1,\qquad
 \sum\limits_{j\in I_n}\|D\varphi_{n,j}(x)\|\leqslant D_n.
\end{equation}
\end{lemma}
\begin{proof}
For each $a\in U$, choose a member $O_a$ of the given cover containing $a$.
By continuity of $DJ$ at $a$, there are an open neighborhood $a\in V_a\subset O_a$ and an integer $m_a\geqslant1$
such that $\|DJ(x)\|\leqslant m_a$ for every $x\in V_a$.
Since $J^{-1}:J(X)\to X$ is continuous, choose an ambient open neighborhood $W_a\ni J(a)$ such that
$J^{-1}(W_a)\subset V_a$. Then $\|DJ(x)\|\leqslant m_a$ on $J^{-1}(W_a)$.
On $V=\bigcup\limits_aW_a$, take $u=0$ and $\tau=1$,
and use Lemma~\ref{ns:lem:catalog} and~\eqref{ns:eq:partitionphi}
to construct a partition of unity $(\theta_{n,j})$ subordinate to $(W_a)$.
Write
$P_{n-1}=\prod\limits_{\ell<n}(1-\chi(W_\ell))$, so that
$\theta_{n,j}=w_{n,j}c(W_n)P_{n-1}$. Then
\begin{equation*}
 \begin{aligned}
 D\theta_{n,j}
 &=c(W_n)P_{n-1}Dw_{n,j}
   +w_{n,j}c'(W_n)P_{n-1}DW_n\\
 &\quad+w_{n,j}c(W_n)DP_{n-1},\\
 \|DP_{n-1}\|&\leqslant\|\chi'\|_\infty
                         \sum\limits_{\ell<n}\|DW_\ell\|.
 \end{aligned}
\end{equation*}
Summing over $j$ and using~\eqref{ns:eq:W-bounds}, we obtain the finite constant
\begin{equation*}
 H_n=\|c\|_\infty E_n+\widehat B_n\|c'\|_\infty E_n
       +\widehat B_n\|c\|_\infty\|\chi'\|_\infty\sum\limits_{l<n}E_l.
\end{equation*}
Thus $\operatorname{supp}_V\theta_{n,j}\subset W_{a(n,j)}$ and, for every $y\in V$,
\begin{equation*}
 \sum\limits_j\|D\theta_{n,j}(y)\|\leqslant H_n.
\end{equation*}
For each fixed $n$, group the functions $\theta_{n,j}\circ J$
according to the value of the integer $m_{a(n,j)}$.
If $D(\theta_{n,j}\circ J)(x)\ne0$, then
$J(x)\in\operatorname{supp}_V\theta_{n,j}\subset W_{a(n,j)}$, and hence
$\|DJ(x)\|\leqslant m_{a(n,j)}$. Therefore, for each fixed pair
$(n,m)\in\mathbb N^2$, the chain rule yields
\begin{equation*}
 \sum\limits_{j:\,m_{a(n,j)}=m}\|D(\theta_{n,j}\circ J)(x)\|
 \leqslant m\sum\limits_j\|D\theta_{n,j}(J(x))\|
 \leqslant mH_n.
\end{equation*}
Fix a bijection $\nu:\mathbb N\to\mathbb N^2$, writing
$\nu(k)=(n(k),m(k))$. Define
\begin{equation*}
 \begin{gathered}
 I_k=\{j:m_{a(n(k),j)}=m(k)\},\qquad
 \varphi_{k,j}=\theta_{n(k),j}\circ J,\\
 D_k=m(k)H_{n(k)}.
 \end{gathered}
\end{equation*}
Then $\sum\limits_{j\in I_k}\|D\varphi_{k,j}(x)\|\leqslant D_k$,
which gives the derivative bound in~\eqref{eq:liftable-layer-bounds}.
For every $x\in U$, some neighborhood of $J(x)$ meets only finitely many sets $\operatorname{supp}_V\theta_{n,j}$.
Its inverse image under $J$ shows that the pulled-back family remains locally finite near $x$.
The reindexing does not change these functions, and
\begin{equation*}
 \operatorname{supp}_U(\theta_{n,j}\circ J)
 \subset J^{-1}(\operatorname{supp}_V\theta_{n,j})
 \subset J^{-1}(W_{a(n,j)})\subset O_{a(n,j)}.
\end{equation*}
For $a\in U$, choose an open neighborhood $W\subset V$ of $J(a)$ with
$J^{-1}(W)\subset U$. Then
\begin{equation*}
 \widetilde\varphi=\theta_{n,j}|_W\in\mathscr R_H^\infty(W,\mathbb R),\qquad
 (\theta_{n,j}\circ J)|_{J^{-1}(W)}
 =\widetilde\varphi\circ J|_{J^{-1}(W)}.
\end{equation*}
Hence $\theta_{n,j}\circ J\in\mathscr F_J^\infty(U,\mathbb R)$.
\end{proof}

\begin{theorem}\label{thm:criterion}
Let $J:X\to H$ be a real-analytic topological embedding. Suppose that there exist
$(\Phi_n)_{n\geqslant1}\subset\mathscr F_J^1(X,X)$ and $C<\infty$ such that
\begin{equation}\label{eq:static-criterion}
 \lim\limits_{n\to\infty}\sup\limits_{x\in X}\|\Phi_n(x)-x\|=0,
 \qquad
 \sup\limits_{n\geqslant1}\sup\limits_{x\in X}\|D\Phi_n(x)\|\leqslant C.
\end{equation}
Then, for every open set $U\subset X$, every real Banach space $Y$, every $f\in C^1(U,Y)$,
and every $\varepsilon\in C(U,(0,\infty))$, there exists $g\in C^\omega(U,Y)$ such that, for every $x\in U$,
\begin{equation*}
 \max\limits_{j=0,1}\|D^jg(x)-D^jf(x)\|<\varepsilon(x).
\end{equation*}
\end{theorem}
\begin{proof}
\smallskip\noindent\firstproofstep{}\label{step:transfer-c1:1}
For each $a\in U$, choose $r_a>0$ such that $B(a,3r_a)\subset U$ and, for every $x\in B(a,3r_a)$,
\begin{equation}\label{eq:transfer-local-oscillation}
 \varepsilon(x)>\frac{\varepsilon(a)}{2},\qquad
 \|Df(x)-Df(a)\|<\lambda_a:=\frac{\varepsilon(a)}{16(C+1)}.
\end{equation}
Take the partition of unity in Lemma~\ref{lem:liftable-localization} subordinate to $(B(a,r_a))_{a\in U}$.
Write its indices as $i=(n,j)$ and denote the corresponding centers and radii by $a_i,r_i$.
For each $i=(n,j)$, choose
\begin{equation}\label{eq:transfer-static-scale}
 0<\delta_i<\min\limits\left\{r_i,
 \frac{2^{-n}\varepsilon(a_i)}{64(1+D_n)(1+\lambda_{a_i})}\right\},
 \qquad
 \sup\limits_X\|\Phi_{\nu_i}-I_X\|<\delta_i.
\end{equation}
Define the local affine mapping $T_i$ centered at $a_i$ by~\eqref{eq:affine-difference-gluing-definition},
and write $R_i=f-T_i$. Set $u_i=R_i\circ\Phi_{\nu_i}$ on $B(a_i,2r_i)$.
By~\eqref{eq:transfer-static-scale}, $\Phi_{\nu_i}(B(a_i,2r_i))\subset B(a_i,3r_i)$.
Since $R_i\in C^1(B(a_i,3r_i),Y)$, the postcomposition property in Lemma~\ref{lem:liftable-closure}
gives $u_i\in\mathscr F_J^1(B(a_i,2r_i),Y)$.
Set $e_i=u_i-R_i$. For $x\in B(a_i,2r_i)$ and $0\leqslant t\leqslant1$,
\begin{equation*}
 \|x+t(\Phi_{\nu_i}(x)-x)-a_i\|
 \leqslant\|x-a_i\|+\|\Phi_{\nu_i}(x)-x\|<3r_i.
\end{equation*}
Since $DR_i(z)=Df(z)-Df(a_i)$, integration along this segment and the
chain rule give
\begin{equation*}
 \begin{aligned}
 e_i(x)&=\int\limits_0^1
  DR_i(x+t(\Phi_{\nu_i}(x)-x))[\Phi_{\nu_i}(x)-x] \,\mathrm{d}t,\\
 De_i(x)&=DR_i(\Phi_{\nu_i}(x))\circ D\Phi_{\nu_i}(x)-DR_i(x).
 \end{aligned}
\end{equation*}
Consequently,
\begin{equation}\label{eq:transfer-local-difference-bounds}
 \|e_i\|<\lambda_{a_i}\delta_i,\qquad
 \|De_i\|\leqslant(C+1)\lambda_{a_i}=\frac{\varepsilon(a_i)}{16}
 \quad\text{on }B(a_i,2r_i).
\end{equation}

\smallskip\noindent\proofstep{}\label{step:transfer-c1:2}
In Lemma~\ref{lem:affine-difference-gluing}, take $V_i=B(a_i,2r_i)$ and $v_i=u_i$,
and define $A,b,h$ by~\eqref{eq:affine-difference-gluing-definition}.
Since $y\mapsto z_0$ is an $\mathscr R_H^1$-extension of every constant mapping $x\mapsto z_0$ and $\operatorname{supp}_U\varphi_i\subset B(a_i,r_i)$,
Lemma~\ref{lem:liftable-closure} gives
$b\in\mathscr F_J^1(U,Y)$, $A\in\mathscr F_J^1(U,\mathcal L(X,Y))$.
We use~\eqref{eq:partition-difference-derivative} and write
$(e_iD\varphi_i)[v]=D\varphi_i[v]e_i$; products involving local mappings are extended by zero outside their domains.
Whenever $\varphi_i(x)\ne0$ or $D\varphi_i(x)\ne0$,
we have $x\in B(a_i,r_i)$ and hence $\varepsilon(a_i)<2\varepsilon(x)$.
For such $i=(n,j)$, \eqref{eq:transfer-static-scale} gives
\begin{equation*}
 \|e_i(x)\|<\lambda_{a_i}\delta_i
 <\frac{2^{-n}\varepsilon(x)}{32(1+D_n)}
       \frac{\lambda_{a_i}}{1+\lambda_{a_i}}
 <\frac{\varepsilon(x)}{32}.
\end{equation*}
Therefore
\begin{equation*}
 \begin{aligned}
 \|h(x)-f(x)\|
 &\leqslant\sum\limits_i\varphi_i(x)\|e_i(x)\|
 <\frac{\varepsilon(x)}{32}\sum\limits_i\varphi_i(x)
 =\frac{\varepsilon(x)}{32},\\
 \left\|\sum\limits_i\varphi_i(x)De_i(x)\right\|
 &\leqslant\sum\limits_i\varphi_i(x)\frac{\varepsilon(a_i)}{16}
 <\frac{\varepsilon(x)}{8}.
 \end{aligned}
\end{equation*}
On the other hand, by~\eqref{eq:liftable-layer-bounds} and~\eqref{eq:transfer-static-scale},
\begin{equation*}
 \begin{aligned}
 \left\|\sum\limits_i e_i(x)D\varphi_i(x)\right\|
 &\leqslant\sum\limits_{n\geqslant1}\sum\limits_{j\in I_n}
             \lambda_{a_{n,j}}\delta_{n,j}\|D\varphi_{n,j}(x)\|\\
 &\leqslant\frac{\varepsilon(x)}{32}
       \sum\limits_{n\geqslant1}2^{-n}\frac{D_n}{1+D_n}
 \leqslant\frac{\varepsilon(x)}{32}.
 \end{aligned}
\end{equation*}
Consequently,
\begin{equation}\label{eq:transfer-affine-approximation}
 \max\limits_{j=0,1}\|D^jh(x)-D^jf(x)\|<\frac{5\varepsilon(x)}{32}.
\end{equation}

\smallskip\noindent\proofstep{}\label{step:transfer-c1:3}
By Proposition~\ref{prop:analytic-liftable}, $(b,A)$ admits fine $C^1$ real-analytic approximation.
Lemma~\ref{lem:affine-analytic-approximation} therefore gives $g\in C^\omega(U,Y)$ such that
$\max\limits_{j=0,1}\|D^jg-D^jh\|<\frac{\varepsilon}{8}$.
Combining this with~\eqref{eq:transfer-affine-approximation}, we obtain
$\max\limits_{j=0,1}\|D^jg-D^jf\|<\frac{9\varepsilon}{32}<\varepsilon$.
\end{proof}

\subsection[Decompositions and polynomial estimates]{Unconditional decompositions and polynomial estimates}\label{subsec:ufdd-polynomial-estimates}

\begin{lemma}\label{lem:unconditional-renorming}
Let $(E_i)_{i\in I}$ be a family of closed subspaces of a real Banach space $(X,\|\cdot\|_0)$.
Suppose that every $x\in X$ has a unique representation as an unconditional sum, convergent as a net of finite partial sums,
$x=\sum\limits_{i\in I}u_i$, where $u_i\in E_i$. Then $R_ix=u_i$ defines a bounded linear projection, and
\begin{equation*}
 \|x\|_u=\sup\limits_{\substack{F\in\operatorname{Fin}(I)\\ |a_i|\leqslant1}}
              \left\|\sum\limits_{i\in F}a_iR_ix\right\|_0
\end{equation*}
is an equivalent norm. For every $(a_i)\in\ell^\infty(I)$, the multiplier sum converges unconditionally and satisfies
\begin{equation}\label{eq:unconditional-decomposition-multiplier}
 \left\|\sum\limits_{i\in I}a_iR_ix\right\|_u
 \leqslant\left(\sup\limits_{i\in I}|a_i|\right)\|x\|_u.
\end{equation}
\end{lemma}
\begin{proof}
Inspired by~\cite[vol.~I, Proposition~1.c.7 and the preceding discussion of renorming, pp.~18--19; Definition~1.g.1, pp.~47--48]{LindenstraussTzafriri1977},
let $\mathcal U$ be the space of all unconditionally summable families $u=(u_i)\in\prod\limits_iE_i$, equipped with the norm
\begin{equation*}
 \|u\|_{\mathcal U}
 =\sup\limits_{F\in\operatorname{Fin}(I)}
       \left\|\sum\limits_{i\in F}u_i\right\|_0.
\end{equation*}
Write $s_F=\sum\limits_{i\in F}u_i$. Since the net of partial sums converges, there exist a finite set $F_0$ and $C<\infty$
such that $\|s_G\|_0\leqslant C$ whenever $G\supset F_0$. Hence
\begin{equation*}
 s_F=s_{F\cup F_0}-s_{F_0\setminus F},\qquad
 \sup\limits_F\|s_F\|_0
 \leqslant C+\max\limits_{E\subset F_0}\|s_E\|_0<\infty.
\end{equation*}
If $(u^{(n)})$ is a Cauchy sequence, closedness of each $E_i$ gives
$u_i=\lim\limits_{n\to\infty}u_i^{(n)}\in E_i$ and
\begin{equation*}
 \sup\limits_F\left\|\sum\limits_{i\in F}(u_i-u_i^{(n)})\right\|_0
 \leqslant\liminf\limits_{m\to\infty}
       \|u^{(n)}-u^{(m)}\|_{\mathcal U}\longrightarrow0.
\end{equation*}
Given $\epsilon>0$, choose $n$ and $F_0\in\operatorname{Fin}(I)$ such that
\begin{gather*}
 \sup\limits_{F\in\operatorname{Fin}(I)}
 \left\|\sum\limits_{i\in F}(u_i-u_i^{(n)})\right\|_0
 <\frac{\epsilon}{2},\\
 \left\|\sum\limits_{i\in F}u_i^{(n)}\right\|_0
 <\frac{\epsilon}{2}\quad\text{for every }F\in\operatorname{Fin}(I\setminus F_0).
\end{gather*}
For $F\in\operatorname{Fin}(I\setminus F_0)$, we have
\begin{equation*}
 \left\|\sum\limits_{i\in F}u_i\right\|_0
 \leqslant\left\|\sum\limits_{i\in F}(u_i-u_i^{(n)})\right\|_0
 +\left\|\sum\limits_{i\in F}u_i^{(n)}\right\|_0<\epsilon.
\end{equation*}
Completeness of $X$ and this unconditional Cauchy condition on the tails imply that $u\in\mathcal U$.
Moreover, $\|u-u^{(n)}\|_{\mathcal U}\to0$, so $\mathcal U$ is complete.

The summation operator $S:\mathcal U\to X$ satisfies $\|Su\|_0\leqslant\|u\|_{\mathcal U}$.
Existence and uniqueness of the decomposition show that it is a linear bijection. By the open mapping theorem,
its inverse is bounded. Thus there is $K\geqslant1$ such that
\begin{equation*}
 \|u\|_{\mathcal U}\leqslant K\|Su\|_0,
 \qquad\left\|\sum\limits_{i\in F}R_ix\right\|_0\leqslant K\|x\|_0.
\end{equation*}
Uniqueness of the decomposition implies that $R_i$ is linear and $R_iR_j=\delta_{ij}R_i$.
Taking $F=\{i\}$ in the preceding inequality also gives $\|R_i\|\leqslant K$.
For a family of signs $\epsilon\in\{-1,1\}^F$, write $F_\pm=\{i\in F:\epsilon_i=\pm1\}$. Then
\begin{equation*}
 \left\|\sum\limits_{i\in F}\epsilon_iR_ix\right\|_0
 \leqslant\left\|\sum\limits_{i\in F_+}R_ix\right\|_0
 +\left\|\sum\limits_{i\in F_-}R_ix\right\|_0
 \leqslant2K\|x\|_0.
\end{equation*}
For $|a_i|\leqslant1$, define
\begin{equation*}
 \lambda_\epsilon=\prod\limits_{i\in F}\frac{1+\epsilon_i a_i}{2},
 \qquad\text{for every }\epsilon\in\{-1,1\}^F.
\end{equation*}
Then $\lambda_\epsilon\geqslant0$,
$\sum\limits_\epsilon\lambda_\epsilon=1$, and
$\sum\limits_\epsilon\lambda_\epsilon\epsilon_i=a_i$. Hence
\begin{equation*}
 \left\|\sum\limits_{i\in F}a_iR_ix\right\|_0
 =\left\|\sum\limits_\epsilon\lambda_\epsilon
                \sum\limits_{i\in F}\epsilon_iR_ix\right\|_0
 \leqslant\sum\limits_\epsilon\lambda_\epsilon\,2K\|x\|_0=2K\|x\|_0.
\end{equation*}
Consequently,
\begin{equation*}
 \|x\|_0\leqslant\|x\|_u\leqslant2K\|x\|_0.
\end{equation*}
Let $a\in\ell^\infty(I)$. For any finite sets $F_0$ and $F\subset I\setminus F_0$, we have
\begin{equation*}
 \left\|\sum\limits_{i\in F}a_iR_ix\right\|_0
 \leqslant2K\|a\|_\infty
 \left\|x-\sum\limits_{i\in F_0}R_ix\right\|_0.
\end{equation*}
As $F_0$ increases to $I$, the right-hand side tends to zero.
Thus $\bigl(\sum\limits_{i\in F}a_iR_ix\bigr)_F$ is a Cauchy net.
Denote its limit by $T_ax$. Then $R_iT_ax=a_iR_ix$, and hence
\begin{equation*}
\|T_ax\|_u
 =\sup\limits_{\substack{F\in\operatorname{Fin}(I)\\|b_i|\leqslant1}}
       \left\|\sum\limits_{i\in F}b_ia_iR_ix\right\|_0
 \leqslant\left(\sup\limits_i|a_i|\right)\|x\|_u.\qedhere
\end{equation*}
\end{proof}

\begin{remark}\label{rem:equivalent-norm-bounds}
Suppose that $\alpha\|x\|_0\leqslant\|x\|_u\leqslant\beta\|x\|_0$, where
$0<\alpha\leqslant\beta<\infty$. If, with respect to $\|\cdot\|_u$, the mappings $\Phi_n$ satisfy
$\sup\limits_X\|\Phi_n-I_X\|_u\leqslant\delta_n$,
$\sup\limits_X\|D\Phi_n\|_u\leqslant C$, then
\begin{equation}\label{eq:equivalent-norm-geometric-bounds}
 \sup\limits_X\|\Phi_n-I_X\|_0\leqslant\frac{\delta_n}{\alpha},\qquad
 \sup\limits_X\|D\Phi_n\|_0\leqslant\frac{\beta C}{\alpha}.
\end{equation}
Real analyticity, the topological embedding property, and the local $\mathscr R_H^k$-extension property along $J$ are preserved under equivalent renormings.
If a homogeneous polynomial satisfies $Q_0(x)\geqslant\|x\|_0^d$, then
$\beta^dQ_0(x)\geqslant\|x\|_u^d$.
\end{remark}

\begin{theorem}\label{thm:ufdd}
Suppose that $X$ admits a separating polynomial and an unconditional finite-dimensional decomposition indexed by an arbitrary set.
Then there exist an index set $\Gamma_H$, a real Hilbert space
$H=\ell^2(\Gamma_H;\mathbb R)$ with its canonical coordinates,
an entire real-analytic topological embedding $J:X\to H$,
a norm $\|\cdot\|$ equivalent to the given norm on $X$, and
$(\Phi_n)_{n\geqslant1}\subset\mathscr F_J^\infty(X,X)$
such that, for every $n\geqslant1$,
\begin{equation}\label{eq:ufdd-final-static-bounds}
 \sup\limits_X\|\Phi_n-I_X\|<2^{-n},\qquad
 \sup\limits_X\|D\Phi_n\|\leqslant2.
\end{equation}
Near each point, every $\Phi_n$ takes values in the sum of finitely many blocks of the original decomposition.
For every $n\geqslant1$ and $a\in X$, there exist an open neighborhood
$V\subset X$ of $a$, a finite set $S\subset\operatorname{supp}J(a)$,
an open neighborhood $O\subset\mathbb R^S$ of $(J(a))_S$, and
$B\in C^\infty(O,X)$ such that, for every $x\in V$,
\begin{equation}\label{eq:ufdd-final-active-lift}
 (J(x))_S\in O,\qquad
 \Phi_n(x)=B((J(x))_S).
\end{equation}
\end{theorem}

Fix a real Banach space $X$ admitting a separating polynomial and an unconditional finite-dimensional decomposition $(E_i)_{i\in I}$, with $E_i\ne\{0\}$ for every $i\in I$.
The set $\{i\in I:R_ix\ne0\}$ is called the block support of $x$ with respect to this decomposition.
By Lemma~\ref{lem:unconditional-renorming}, fix an equivalent norm satisfying
\eqref{eq:unconditional-decomposition-multiplier}, still denoted by $\|\cdot\|$.
For $A\subset I$, write
\begin{equation*}
 \Pi_Ax=\sum\limits_{i\in A}R_ix,\qquad
 \Pi_A X=\overline{\operatorname{span}}\bigcup\limits_{i\in A}E_i.
\end{equation*}
Then $\|\Pi_A\|\leqslant1$, and changes of signs on the blocks are isometries.
For every $x\in X$, choose finite sets $F_n\subset I$ such that
$\|x-\Pi_{F_n}x\|<\frac{1}{n}$.
If $i\notin F_n$, then
\begin{equation*}
 \|R_ix\|=\|R_i(x-\Pi_{F_n}x)\|
 \leqslant\|x-\Pi_{F_n}x\|<\frac{1}{n}.
\end{equation*}
Thus $\{i:R_ix\ne0\}\subset\bigcup\limits_{n\geqslant1}F_n$; in particular,
every vector has at most countably many nonzero components with respect
to the decomposition.
By Lemma~\ref{lem:homogeneous-separator}, fix an even integer $d\geqslant2$, a continuous symmetric
$d$-linear form $A$, and $B\geqslant1$ such that, for every $x\in X$,
\begin{equation}\label{eq:ufdd-separator-bounds}
 \begin{gathered}
 Q(x)=A(x^d),\qquad \mu(x)=Q(x)^{\frac{1}{d}},\qquad
 \|x\|^d\leqslant Q(x),\\
 \|x\|\leqslant\mu(x)\leqslant B\|x\|.
 \end{gathered}
\end{equation}
Moreover, for every $x\in X\setminus\{0\}$,
\begin{equation*}
 \|D\mu(x)\|\leqslant B.
\end{equation*}
For every $A_0\subset I$, the subspace $\Pi_{A_0}X$ carries the restrictions of this norm and of $Q$.

\begin{lemma}\label{lem:ufdd-polynomial-arrays}
Let $(p_\gamma)_{\gamma\in\Gamma}$ be a family of continuous scalar-valued polynomials on a real Banach space $Z$,
each of degree at most $N$, and suppose that, for every $R>0$,
\begin{equation}\label{eq:ufdd-polynomial-array-bound}
 \sup\limits_{\|z\|\leqslant R}
 \sum\limits_{\gamma\in\Gamma}|p_\gamma(z)|^2<\infty.
\end{equation}
Then $z\mapsto(p_\gamma(z))_\gamma$ is a continuous
$\ell^2(\Gamma)$-valued polynomial of degree at most $N$.
More generally, let $(H_\gamma)_{\gamma\in\Gamma}$ be a family
of real Hilbert spaces, and let $p_\gamma:Z\to H_\gamma$ be
continuous polynomials of degree at most $N$.
If
\begin{equation*}
 \sup\limits_{\|z\|\leqslant R}
 \sum\limits_{\gamma\in\Gamma}
 \|p_\gamma(z)\|_{H_\gamma}^2<\infty
 \qquad\text{for every }R>0,
\end{equation*}
then $z\mapsto(p_\gamma(z))_{\gamma\in\Gamma}$ is a continuous
polynomial of degree at most $N$ with values in
$\bigoplus\limits_{\gamma\in\Gamma}^{\ell^2}H_\gamma$.
\end{lemma}
\begin{proof}
Write $F(z)=(p_\gamma(z))_\gamma$. By assumption, $F$ is bounded on every ball.
Fix distinct real numbers $t_0,\ldots,t_N$, and set
\begin{equation*}
 V=(t_\ell^j)_{0\leqslant\ell,j\leqslant N},\qquad
 V^{-1}=(c_{j\ell})_{0\leqslant j,\ell\leqslant N},\qquad
 F_j(z)=\sum\limits_{\ell=0}^Nc_{j\ell}F(t_\ell z).
\end{equation*}
Vandermonde inversion gives $(F_j)_\gamma=p_\gamma^{(j)}$, where
$p_\gamma^{(j)}$ is the homogeneous component of degree $j$ of $p_\gamma$.
Hence $F=\sum\limits_{j=0}^NF_j$, and each $F_j$ is bounded on bounded balls.
For $j\geqslant1$, set
\begin{equation*}
 A_j(z_1,\ldots,z_j)
 =\frac{1}{2^jj!}
 \sum\limits_{\epsilon\in\{-1,1\}^j}
 \left(\prod\limits_{\nu=1}^j\epsilon_\nu\right)
 F_j\left(\sum\limits_{\nu=1}^j\epsilon_\nu z_\nu\right).
\end{equation*}
By~\cite[Theorems~1.10 and~2.2 and Proposition~2.9, pp.~6, 12, and~14--15]{Mujica1986},
coordinatewise polarization gives
$(A_j)_\gamma=A_{j,\gamma}$, where $A_{j,\gamma}$ is the symmetric $j$-linear
form associated with $p_\gamma^{(j)}$.
Since each scalar coordinate of $A_j$ is symmetric and $j$-linear,
and the coordinate functionals separate points of $\ell^2(\Gamma)$,
the mapping $A_j$ is symmetric and $j$-linear.
Moreover, $F_j(z)=A_j(z^j)$. If $\|z_\nu\|\leqslant1$,
the polarization formula gives
\begin{equation*}
 \|A_j(z_1,\ldots,z_j)\|
 \leqslant\frac{1}{j!}\sup\limits_{\|z\|\leqslant j}\|F_j(z)\|<\infty.
\end{equation*}
Therefore $A_j$ is bounded.
Since $F_0$ is constant, $F$ is a continuous polynomial.
For the Hilbert-space-valued case, choose an orthonormal basis $(e_{\gamma,\nu})_\nu$ of each $H_\gamma$.
The identity
\begin{equation*}
 \sum\limits_{\gamma,\nu}|\langle p_\gamma(z),e_{\gamma,\nu}\rangle|^2
 =\sum\limits_\gamma\|p_\gamma(z)\|^2
\end{equation*}
allows us to apply the scalar assertion to all coordinates $\langle p_\gamma,e_{\gamma,\nu}\rangle$, proving the last assertion.
\end{proof}

\begin{lemma}\label{lem:ufdd-block-summability}
Set $q_i=Q\circ R_i$. Then
\begin{equation}\label{eq:ufdd-block-summability}
 \sum\limits_{i\in I}q_i(x)\leqslant\|A\|\|x\|^d,
 \qquad
 \left(\sum\limits_{i\in I}q_i(x)^2\right)^{\frac{1}{2}}
 \leqslant\|A\|\|x\|^d.
\end{equation}
In particular, $x\mapsto(q_i(x))_{i\in I}$ is a continuous $\ell^2(I)$-valued homogeneous polynomial of degree $d$.
\end{lemma}
\begin{proof}
Inspired by~\cite[proof of Theorem~13, equation~(5.3), p.~158]{DimantDineen1998} and
\cite[proof of Proposition~1.3, pp.~735--736]{DimantGonzalo2001},
fix a finite set $F\subset I$, take independent random variables
$\epsilon_i^s$, for $i\in F$ and $1\leqslant s<d$, each taking the values $-1,1$ with equal probability, and set
\begin{equation*}
 \begin{gathered}
 u_s=\sum\limits_{i\in F}\epsilon_i^sR_ix
 \quad\text{for every }1\leqslant s<d,\\
 u_d=\sum\limits_{i\in F}
 \left(\prod\limits_{s=1}^{d-1}\epsilon_i^s\right)R_ix.
 \end{gathered}
\end{equation*}
By~\eqref{eq:unconditional-decomposition-multiplier}, for every $1\leqslant s\leqslant d$,
$\|u_s\|\leqslant\|x\|$.
The families of random variables corresponding to distinct $s$ are independent, whereas for fixed $s$,
\begin{equation*}
 \mathbb E(\epsilon_i^s\epsilon_j^s)
 =\begin{cases}1,&i=j,\\0,&i\ne j.\end{cases}
\end{equation*}
Multilinearity of $A$ and finite summation give
\begin{equation*}
 \begin{aligned}
 \mathbb E A(u_1,\ldots,u_d)
 &=\sum\limits_{i_1,\ldots,i_d\in F}
   A(R_{i_1}x,\ldots,R_{i_d}x)
   \prod\limits_{s=1}^{d-1}
   \mathbb E(\epsilon_{i_s}^s\epsilon_{i_d}^s)\\
 &=\sum\limits_{i\in F} A(R_ix,\ldots,R_ix)
 =\sum\limits_{i\in F}Q(R_ix).
 \end{aligned}
\end{equation*}
Therefore,
\begin{equation*}
 \sum\limits_{i\in F}Q(R_ix)
 \leqslant\mathbb E|A(u_1,\ldots,u_d)|
 \leqslant\|A\|\|x\|^d.
\end{equation*}
Since $q_i(x)\geqslant0$, taking the supremum over finite sets $F$ gives the first inequality in~\eqref{eq:ufdd-block-summability}.
Moreover,
\begin{equation*}
 \sum\limits_iq_i(x)^2
 \leqslant\left(\sum\limits_iq_i(x)\right)^2
 \leqslant\|A\|^2\|x\|^{2d}.
\end{equation*}
The Hilbert-valued polynomial assertion follows from Lemma~\ref{lem:ufdd-polynomial-arrays}.
\end{proof}

By Theorem~\ref{thm:general-polynomial-hilbert-embedding}, choose a polynomial topological embedding
$J_\star:X\to H_\star$. Fix
\begin{equation}\label{eq:ufdd-baseline-embedding}
 \begin{gathered}
 H_0=H_\star\oplus_2\ell^2(I),\\
 J_0(x)=\bigl(J_\star(x),(q_i(x))_{i\in I}\bigr),\qquad
 \widehat J_0(x)=\bigl(J_0(x),1+\|J_0(x)\|^2\bigr).
 \end{gathered}
\end{equation}
Here $J_0:X\to H_0$ and $\widehat J_0:X\to H_0\oplus_2\mathbb R$.
Both are continuous polynomial topological embeddings of finite degree, and each $q_i$ is a scalar coordinate of $J_0$.

\subsection[Approximations with locally finite-dimensional range]{Approximations of the identity with locally finite-dimensional range}\label{subsec:ufdd-tail-cutoffs}

\begin{lemma}\label{lem:ufdd-type-two}
There exists $T\geqslant1$ such that, for every finite family of vectors $(x_j)$ and independent symmetric random signs $(\epsilon_j)$,
\begin{equation}\label{eq:ufdd-type-two}
 \left(\mathbb E\left\|\sum\limits_j\epsilon_jx_j\right\|^2\right)^{\frac{1}{2}}
 \leqslant T\left(\sum\limits_j\|x_j\|^2\right)^{\frac{1}{2}}.
\end{equation}
If a finite family $(u_j)$ has pairwise disjoint block supports, then
\begin{equation}\label{eq:ufdd-disjoint-type-two}
 \left\|\sum\limits_ju_j\right\|
 \leqslant T\left(\sum\limits_j\|u_j\|^2\right)^{\frac{1}{2}}.
\end{equation}
The same constant $T$ applies to every subspace $\Pi_{A_0}X$, $A_0\subset I$.
\end{lemma}
\begin{proof}
By~\cite[Chapter~5, Theorem~64, pp.~269--270; proof of Theorem~67, p.~273]{HajekJohanis2014},
we may choose $T$ such that~\eqref{eq:ufdd-type-two} holds.
If the block supports are pairwise disjoint, then~\eqref{eq:unconditional-decomposition-multiplier}
shows that blockwise changes of signs and their inverses are contractions, hence isometries. Thus
\begin{equation*}
 \left\|\sum\limits_j u_j\right\|
 =\left(\mathbb E\left\|\sum\limits_j\epsilon_ju_j\right\|^2\right)^{\frac{1}{2}}
 \leqslant T\left(\sum\limits_j\|u_j\|^2\right)^{\frac{1}{2}}.
\end{equation*}
For every $A_0\subset I$, restricting~\eqref{eq:ufdd-type-two} to vectors in $\Pi_{A_0}X$ preserves the constant $T$.
\end{proof}

Inspired by~\cite[Section~1.A, proof of Basic Lemma~1, pp.~297--300]{Moulis1971},
use the step function $\sigma$ in~\eqref{ns:eq:smooth-step-function}, and set
$c_\sigma=\frac{\|\sigma'\|_\infty}{\log2}$. For an integer $L\geqslant1$, define
\begin{equation}\label{eq:ufdd-logarithmic-cutoff}
 \chi_L(t)=
 \begin{cases}
 \sigma\left(1+\frac{\log_2t}{L}\right),&t>0,\\
 0,&t=0.
 \end{cases}
\end{equation}
Then $\chi_L(t)=0$ for every $0\leqslant t\leqslant2^{-L}$, and $\chi_L(t)=1$ for every $t\geqslant1$.
For $t>0$, differentiation gives
\begin{equation*}
 \chi_L'(t)=\frac{1}{Lt\log2}
             \sigma'\left(1+\frac{\log_2t}{L}\right).
\end{equation*}
Since $\chi_L$ is zero near $0$, this implies, for every $t\geqslant0$,
\begin{equation}\label{eq:ufdd-logarithmic-cutoff-bounds}
 |t\chi_L'(t)|\leqslant\frac{c_\sigma}{L}.
\end{equation}
Fix a well-ordering of $I$; if $I$ is infinite, choose one satisfying~\eqref{eq:well-order-small-parts}.
Throughout this section, $<$ and $\leqslant$ refer to this well-ordering. Write
\begin{equation*}
 T_i=\Pi_{\{j:j\geqslant i\}},\qquad
 \Pi_i=\Pi_{\{j:j\leqslant i\}},\qquad
 T_{i+}=\Pi_{\{j:j>i\}}.
\end{equation*}

\begin{lemma}\label{lem:ufdd-tail-cutoff}
For $r>0$, set
\begin{equation}\label{eq:ufdd-ordered-tail-map}
 \Theta_{r,L}(x)
 =\sum\limits_{i\in I}
 \chi_L\left(\frac{\mu(T_ix)}{r}\right)R_ix.
\end{equation}
For every fixed $r>0$ and $L\geqslant1$, we have $\Theta_{r,L}\in C^\infty(X,X)$ and
\begin{equation}\label{eq:ufdd-ordered-tail-bounds}
 \|\Theta_{r,L}(x)-x\|\leqslant r,\qquad
 \|D\Theta_{r,L}(x)\|
 \leqslant1+\frac{2c_\sigma BT}{\sqrt L}.
\end{equation}
For fixed $r,L$, define $V_0$ and, for every $i\in I$ and $m\in\mathbb N$, $V_{i,m}$ by
\begin{equation}\label{eq:ufdd-tail-chart-cover}
 \begin{aligned}
 V_0&=\{x:B\|x\|<r2^{-L}\},\\
 V_{i,m}&=\{x:q_i(x)>0,\ B\|T_{i+}x\|<r2^{-L},\ \|x\|<m\}.
 \end{aligned}
\end{equation}
These open sets cover $X$.
We have $\Theta_{r,L}=0$ on $V_0$, and on $V_{i,m}$,
\begin{equation}\label{eq:ufdd-prefix-localization}
 \Theta_{r,L}=\Pi_i\Theta_{r,L},\qquad
 D\Theta_{r,L}=\Pi_iD\Theta_{r,L}.
\end{equation}
\end{lemma}
\begin{proof}
\smallskip\noindent\firstproofstep{}\label{step:ufdd-truncation:1}
Write the coefficient as $\beta_r(Q(T_ix))$, where
\begin{equation}\label{eq:ufdd-tail-coefficient}
 \begin{gathered}
 \beta_r(s)=\chi_L\left(\frac{s^{\frac{1}{d}}}{r}\right)
 \quad\text{for every }s\geqslant0,\\
 \beta_r(s)=0\quad\text{for every }s<0.
 \end{gathered}
\end{equation}
Fix $r,L$. The function $\beta_r$ vanishes on a neighborhood of the origin, and hence belongs to $C^\infty(\mathbb R)$.
The mapping $P(x)=(Q(T_ix))_{i\in I}$ is a continuous $\ell^\infty(I)$-valued homogeneous polynomial of degree $d$,
with associated multilinear mapping
$(x_1,\ldots,x_d)\mapsto(A(T_ix_1,\ldots,T_ix_d))_i$, whose norm is at most $\|A\|$.
Set $N_{\beta_r}(s)=(\beta_r(s_i))_i$. For $k\geqslant0$, define coordinatewise
\begin{equation*}
 B_k(s)[h_1,\ldots,h_k]
 =\left(\beta_r^{(k)}(s_i)\prod\limits_{\nu=1}^k h_{\nu,i}\right)_i,
 \qquad B_0=N_{\beta_r}.
\end{equation*}
If $\|s\|_\infty,\|s+h\|_\infty\leqslant R$, the one-dimensional mean value theorem and Taylor's formula give
\begin{equation*}
 \begin{aligned}
 \|B_k(s+h)-B_k(s)\|
 &\leqslant\sup\limits_{|t|\leqslant R}|\beta_r^{(k+1)}(t)|\,\|h\|_\infty,\\
 \|B_k(s+h)-B_k(s)-B_{k+1}(s)[h]\|
 &\leqslant\frac{1}{2}\sup\limits_{|t|\leqslant R}|\beta_r^{(k+2)}(t)|
           \,\|h\|_\infty^2.
 \end{aligned}
\end{equation*}
For $k\geqslant1$, each $B_k(s)$ is a bounded multilinear mapping.
In the preceding estimates, the vector norm is used for $k=0$, and the multilinear operator norm otherwise.
Hence $DB_k=B_{k+1}$, so $N_{\beta_r}$ is smooth. It follows that
$b=N_{\beta_r}\circ P\in C^\infty(X,\ell^\infty(I))$.
The block multiplier mapping
\begin{equation}\label{eq:ufdd-bilinear-multiplier}
 M(a,x)=\sum\limits_{i\in I}a_iR_ix,\qquad
 \|M(a,x)\|\leqslant\|a\|_\infty\|x\|
\end{equation}
is bounded and bilinear. Therefore $\Theta_{r,L}(x)=M(b(x),x)$ belongs to
$C^\infty(X,X)$, with
\begin{equation*}
 \begin{aligned}
 D\Theta_{r,L}(x)[h]&=M(Db(x)[h],x)+M(b(x),h),\\
 (Db(x)[h])_i&=\beta_r'(Q(T_ix))DQ(T_ix)[T_ih].
 \end{aligned}
\end{equation*}
These identities take place in $X$ and $\ell^\infty(I)$, respectively,
so the derivative is obtained before expanding the multiplier sums.

\smallskip\noindent\proofstep{}\label{step:ufdd-truncation:2}
Write $a_i=\chi_L(\frac{\mu(T_ix)}{r})$. If $\{i:a_i\ne1\}\ne\varnothing$, set
$i_0=\min\limits\{i:a_i\ne1\}$. Then
\begin{equation*}
 \mu(T_{i_0}x)<r,\qquad
 x-\Theta_{r,L}(x)=\sum\limits_{i\geqslant i_0}(1-a_i)R_iT_{i_0}x.
\end{equation*}
Since $0\leqslant a_i\leqslant1$, \eqref{eq:unconditional-decomposition-multiplier} gives
\begin{equation*}
 \|x-\Theta_{r,L}(x)\|
 \leqslant\|T_{i_0}x\|\leqslant\mu(T_{i_0}x)<r.
\end{equation*}
If $a_i=1$ for every $i$, then $\Theta_{r,L}(x)=x$.
\smallskip\noindent\proofstep{}\label{step:ufdd-truncation:3}
For $h\in X$, we have
\begin{equation*}
 \begin{aligned}
 D\Theta_{r,L}(x)[h]
 &=\sum\limits_i\chi_L\left(\frac{\mu(T_ix)}{r}\right)R_ih+v,\\
 v&=\sum\limits_i
 \frac{1}{r}\chi_L'\left(\frac{\mu(T_ix)}{r}\right)
 D\mu(T_ix)[T_ih]R_ix.
 \end{aligned}
\end{equation*}
By~\eqref{eq:unconditional-decomposition-multiplier} and $0\leqslant\chi_L\leqslant1$,
\begin{equation*}
 \left\|\sum\limits_i\chi_L\left(\frac{\mu(T_ix)}{r}\right)R_ih\right\|
 \leqslant\|h\|.
\end{equation*}
When $T_ix=0$, the corresponding term is defined to be zero.
For every integer $k$ with $0\leqslant k<L$, set
\begin{equation*}
 A_k=\{i:r2^{-k-1}<\mu(T_ix)\leqslant r2^{-k}\}.
\end{equation*}
The sets $A_k$ are pairwise disjoint, and the index of every nonzero term in $v$ belongs to one of them.
If $A_k$ is nonempty, denote its least index by $i_k$. Then
\begin{equation*}
 \|\Pi_{A_k}x\|\leqslant\|T_{i_k}x\|\leqslant r2^{-k}.
\end{equation*}
For $i\in A_k$, \eqref{eq:ufdd-logarithmic-cutoff-bounds} gives
\begin{equation*}
 \frac{1}{r}\left|\chi_L'\left(\frac{\mu(T_ix)}{r}\right)
 D\mu(T_ix)[T_ih]\right|
 \leqslant\frac{c_\sigma B}{L\mu(T_ix)}\|h\|
 \leqslant\frac{c_\sigma B2^{k+1}}{Lr}\|h\|.
\end{equation*}
Applying~\eqref{eq:unconditional-decomposition-multiplier} again, we obtain
\begin{equation*}
 \|\Pi_{A_k}v\|
 \leqslant\frac{c_\sigma B2^{k+1}}{Lr}\|h\|\,\|\Pi_{A_k}x\|
 \leqslant\frac{2c_\sigma B}{L}\|h\|.
\end{equation*}
By~\eqref{eq:ufdd-disjoint-type-two},
\begin{equation*}
 \|v\|\leqslant
 T\left(\sum\limits_{k=0}^{L-1}\|\Pi_{A_k}v\|^2\right)^{\frac{1}{2}}
 \leqslant\frac{2c_\sigma BT}{\sqrt L}\|h\|.
\end{equation*}
Therefore,
\begin{equation*}
 \|D\Theta_{r,L}(x)\|
 \leqslant1+\frac{2c_\sigma BT}{\sqrt L}.
\end{equation*}

\smallskip\noindent\proofstep{}\label{step:ufdd-truncation:4}
For $x\ne0$, choose a nonempty finite subset $F$ of its block support such that
$\|x-\Pi_Fx\|<\frac{r2^{-L}}{B}$. Set $i=\max\limits F$. Then
\begin{equation*}
 q_i(x)>0,\qquad
 \|T_{i+}x\|\leqslant\|x-\Pi_Fx\|<\frac{r2^{-L}}{B}.
\end{equation*}
Increasing $m$ gives $x\in V_{i,m}$, while $0\in V_0$.
If $x\in V_{i,m}$ and $j>i$, then
\begin{equation*}
 \mu(T_jx)\leqslant B\|T_jx\|
 \leqslant B\|T_{i+}x\|<r2^{-L}.
\end{equation*}
Thus $\chi_L(\frac{\mu(T_jx)}{r})=0$, and hence
$\Theta_{r,L}|_{V_{i,m}}=\Pi_i\Theta_{r,L}|_{V_{i,m}}$.
Differentiating on the open set $V_{i,m}$ gives \eqref{eq:ufdd-prefix-localization}.
If $x\in V_0$, then, for every $j\in I$,
$\mu(T_jx)\leqslant B\|x\|<r2^{-L}$, so $\Theta_{r,L}(x)=0$.
\end{proof}

\begin{lemma}\label{lem:ufdd-fixed-chart-polynomials}
For $p\in\mathbb R[t]$, set
\begin{equation}\label{eq:ufdd-polynomial-tail-map}
 S_p(x)=\sum\limits_{j\in I}p(Q(T_jx))R_jx,
 \qquad K_{i,p}=\Pi_iS_p.
\end{equation}
Both mappings are continuous vector-valued polynomials. Fix $r,L,m$, and, for $\beta_r$ defined in~\eqref{eq:ufdd-tail-coefficient}, write
\begin{equation*}
 \zeta_0=\sup\limits_{[0,B^dm^d]}|p-\beta_r|,
 \qquad
 \zeta_1=\sup\limits_{[0,B^dm^d]}|p'-\beta_r'|.
\end{equation*}
On $\|x\|<m$, we have
\begin{equation}\label{eq:ufdd-fixed-chart-estimates}
 \begin{aligned}
 \|S_p(x)-\Theta_{r,L}(x)\|&\leqslant m\zeta_0,\\
 \|DS_p(x)-D\Theta_{r,L}(x)\|
 &\leqslant\zeta_0+d\|A\|m^d\zeta_1.
 \end{aligned}
\end{equation}
On $V_{i,m}$, the estimates in~\eqref{eq:ufdd-fixed-chart-estimates} also hold for $K_{i,p}-\Theta_{r,L}$.
Consequently, for every $\tau>0$, one can choose $p\in\mathbb Q[t]$, independently of $i$, such that, for every $i\in I$ and $x\in V_{i,m}$,
\begin{equation*}
 \max\limits_{j=0,1}\|D^jK_{i,p}(x)-D^j\Theta_{r,L}(x)\|<\tau.
\end{equation*}
\end{lemma}
\begin{proof}
Use the multiplication $(ab)_j=a_jb_j$ on $\ell^\infty(I)$.
Compose the $\ell^\infty(I)$-valued polynomial
$P(x)=(Q(T_jx))_j$ from Step~\ref{step:ufdd-truncation:1} with
$a\mapsto p(a)=(p(a_j))_j$, and then apply the bounded bilinear mapping
$M$ in~\eqref{eq:ufdd-bilinear-multiplier}. More explicitly, if
$p(t)=\sum\limits_{\nu=0}^Na_\nu t^\nu$, then
\begin{equation*}
 \begin{gathered}
 S_p(x)=\sum\limits_{\nu=0}^Na_\nu M(P(x)^\nu,x),\qquad
 K_{i,p}(x)=\Pi_iS_p(x),\\
 \deg S_p,\deg K_{i,p}\leqslant dN+1,
 \end{gathered}
\end{equation*}
where powers in $\ell^\infty(I)$ are coordinatewise and $P(x)^0=(1)_{i\in I}$.
Thus both mappings are continuous polynomials.
Set $e=p-\beta_r$. For every $\|x\|<m$ and $h\in X$,
\begin{equation*}
 \begin{aligned}
 S_p(x)-\Theta_{r,L}(x)
 &=\sum\limits_j e(Q(T_jx))R_jx,\\
 D(S_p-\Theta_{r,L})(x)[h]
 &=\sum\limits_j e(Q(T_jx))R_jh
   +\sum\limits_j e'(Q(T_jx))DQ(T_jx)[T_jh]R_jx.
 \end{aligned}
\end{equation*}
Since $0\leqslant Q(T_jx)\leqslant B^dm^d$, we have
\begin{equation*}
 |e(Q(T_jx))|\leqslant\zeta_0,\qquad
 |e'(Q(T_jx))|\leqslant\zeta_1,\qquad
 \|DQ(T_jx)T_j\|\leqslant d\|A\|m^{d-1}.
\end{equation*}
Applying~\eqref{eq:unconditional-decomposition-multiplier} to each block multiplier sum above gives
\begin{equation*}
 \|S_p(x)-\Theta_{r,L}(x)\|\leqslant m\zeta_0,\qquad
 \|D(S_p-\Theta_{r,L})(x)[h]\|
 \leqslant(\zeta_0+d\|A\|m^d\zeta_1)\|h\|,
\end{equation*}
which is~\eqref{eq:ufdd-fixed-chart-estimates}.
By \eqref{eq:ufdd-prefix-localization}, applying the projection $\Pi_i$, with $\|\Pi_i\|\leqslant1$,
to the difference and its derivative on $V_{i,m}$ gives the estimates for $K_{i,p}$.
Finally, given $\tau>0$, set $M=B^dm^d$ and
$\rho=\frac{\tau}{2(1+m+d\|A\|m^d)}$.
By the Weierstrass theorem, choose a polynomial $v$ such that
$\|v-\beta_r'\|_{C([0,M])}<\frac{\rho}{2(1+M)}$.
Set $p_0(t)=\beta_r(0)+\int\limits_0^t v(s)\,\mathrm{d}s$. Then
\begin{equation*}
 \max\limits_{j=0,1}\|p_0^{(j)}-\beta_r^{(j)}\|_{C([0,M])}<\frac{\rho}{2}.
\end{equation*}
Approximate the finitely many coefficients of $p_0$ by rational numbers so that the resulting $p\in\mathbb Q[t]$ satisfies
$\max\limits_{j=0,1}\|p^{(j)}-p_0^{(j)}\|_{C([0,M])}<\frac{\rho}{2}$.
Then $\max\limits\{\zeta_0,\zeta_1\}<\rho$, and~\eqref{eq:ufdd-fixed-chart-estimates} gives
\begin{equation*}
 \max\limits_{j=0,1}\|D^jK_{i,p}(x)-D^j\Theta_{r,L}(x)\|
 <(1+m+d\|A\|m^d)\rho=\frac{\tau}{2}<\tau.
\end{equation*}
This choice depends only on $r,L,m,\tau$, so the same polynomial $p$ applies to every $i$.
\end{proof}

\begin{corollary}\label{cor:ufdd-quantitative-partition}
Let $J_0,\widehat J_0$ be as in \eqref{eq:ufdd-baseline-embedding}.
For an arbitrary open cover of $X$, apply Lemma~\ref{lem:liftable-localization} with $U=X$ and $J=J_0$.
Near every point, each member of the resulting partition of unity is a smooth function of finitely many scalar coordinates of $\widehat J_0$,
and these coordinates can be chosen to be nonzero at that point.
\end{corollary}
\begin{proof}
Denote the resulting partition of unity by $(\varphi_{n,j})$. Fix $a\in X$ and $(n,j)$.
By Definitions~\ref{def:liftable-class} and~\ref{ns:def:class},
on some neighborhood of $a$ we may write
\begin{equation*}
 \varphi_{n,j}(x)
 =b\bigl((J_0(x))_S,\|J_0(x)\|^2\bigr),
 \qquad S\subset\operatorname{supp}J_0(a),\quad \#S<\infty,
\end{equation*}
where $b$ is smooth. By~\eqref{eq:ufdd-baseline-embedding},
the additional coordinate of $\widehat J_0(x)$ is $1+\|J_0(x)\|^2$, so we may also write
\begin{equation*}
 \varphi_{n,j}(x)
 =\widetilde b\bigl((J_0(x))_S,1+\|J_0(x)\|^2\bigr),
 \qquad \widetilde b(s,t)=b(s,t-1).
\end{equation*}
The selected coordinates $(J_0)_\gamma$, $\gamma\in S$,
are nonzero at $a$, and the additional coordinate
$1+\|J_0(\,\cdot\,)\|^2$ is everywhere positive.
This yields the required representation in nonzero polynomial coordinates.
\end{proof}

\begin{definition}\label{def:ufdd-polynomial-features}
Let $X,Z$ be real Banach spaces, and let
$\mathcal E=(\mathcal P_k:X\to H_k)_{k\in\mathbb N}$ be a countable family of continuous Hilbert-valued polynomials,
where each real Hilbert space $H_k$ is equipped with fixed orthonormal coordinates.
A mapping $u:X\to Z$ is said to have a \emph{representation in nonzero polynomial coordinates} relative to $\mathcal E$ if, for every
$a\in X$, there exist finitely many scalar coordinates $p_1,\ldots,p_N$ of these polynomials,
an open neighborhood $V$ of $a$, an open set $O\subset\mathbb R^N$, and $b\in C^\infty(O,Z)$ such that
$(p_1,\ldots,p_N)(V)\subset O$, $p_\nu(a)\ne0$ for every $1\leqslant\nu\leqslant N$, and, for every $x\in V$,
\begin{equation}\label{eq:ufdd-active-feature-representation}
 u(x)=b(p_1(x),\ldots,p_N(x)).
\end{equation}
The case $N=0$ is allowed and means that $u|_V$ is constant.
A finite family may be extended by zero polynomials.
\end{definition}

\begin{proposition}
\label{prop:ufdd-strengthened-approximation}
There exists a countable polynomial family $\mathcal E_X$ containing $\widehat J_0$ such that, for every
$\delta,\varepsilon>0$, there is $\Phi_{\delta,\varepsilon}\in C^\infty(X,X)$ satisfying
\begin{equation}\label{eq:ufdd-strengthened-approximation}
 \sup\limits_{x\in X}\|\Phi_{\delta,\varepsilon}(x)-x\|<\delta,
 \qquad
 \sup\limits_{x\in X}\|D\Phi_{\delta,\varepsilon}(x)\|
 \leqslant1+\varepsilon.
\end{equation}
On a neighborhood of every point, this mapping takes values in the sum of finitely many original blocks and has a representation
in nonzero polynomial coordinates relative to $\mathcal E_X$.
The polynomial family $\mathcal E_X$ is independent of $\delta,\varepsilon$.
\end{proposition}
\begin{proof}
We argue by induction on the cardinality $|I|$ of the block index set.

\smallskip\noindent\firstproofstep{}\label{step:uniform-identity:1}
If $X=\{0\}$, take $\mathcal E_X=(\widehat J_0,0,0,\ldots)$ and
$\Phi_{\delta,\varepsilon}=0$.
If $1\leqslant\dim X=N<\infty$, take $\Phi_{\delta,\varepsilon}=I_X$ and adjoin to
$\widehat J_0$ the linear coordinates $\ell_1,\ldots,\ell_N$ and the affine coordinates $1+\ell_j$.
At every $a\in X$, we have $(\ell_j(a),1+\ell_j(a))\ne(0,0)$,
so one may select from each pair a coordinate that is nonzero at $a$.
Let $e_1,\ldots,e_N\in X$ be the basis satisfying
$\ell_j(e_k)=\delta_{jk}$.
Choose $c_j\in\{0,1\}$ so that $p_j=\ell_j+c_j$ satisfies $p_j(a)\ne0$.
Then
\begin{equation*}
 B(t_1,\ldots,t_N)=\sum\limits_{j=1}^N(t_j-c_j)e_j,\qquad
 \Phi_{\delta,\varepsilon}(x)=x=B(p_1(x),\ldots,p_N(x)).
\end{equation*}

If there are countably infinitely many blocks, enumerate them as $E_n$, and choose a finite system of linear coordinates
$(\xi_{n,l})_{1\leqslant l\leqslant\dim E_n}$ on each block. For every $n\geqslant1$, adjoin to the fixed polynomial family
\begin{equation*}
 \xi_{n,l}\circ R_n,\qquad 1+\xi_{n,l}\circ R_n,
 \qquad 1+Q\circ T_n.
\end{equation*}
These scalar-valued polynomials form a countable family. Given $\delta,\varepsilon>0$, choose
\begin{equation*}
 \begin{gathered}
 0<r<\delta,\qquad L\in\mathbb N,\qquad
 L\geqslant\max\limits\left\{1,\left(\frac{2c_\sigma BT}{\varepsilon}\right)^2\right\},\\
 \Phi_{\delta,\varepsilon}=\Theta_{r,L}.
 \end{gathered}
\end{equation*}
By \eqref{eq:ufdd-ordered-tail-bounds},
\begin{equation*}
 \|D\Phi_{\delta,\varepsilon}\|_\infty
 \leqslant1+\frac{2c_\sigma BT}{\sqrt L}\leqslant1+\varepsilon.
\end{equation*}
Fix $a\in X$. Since $T_na\to0$, choose $N$ such that
$B\|T_Na\|<r2^{-L}$. Set $V=\{x\in X:B\|T_Nx\|<r2^{-L}\}$.
By continuity of $T_N$, the set $V$ is an open neighborhood of $a$.
For every $x\in V$ and $n\geqslant N$, the identities $T_n=T_nT_N$ and $\|T_n\|\leqslant1$ imply that
\begin{equation*}
 \mu(T_nx)\leqslant B\|T_nx\|
 \leqslant B\|T_Nx\|<r2^{-L}.
\end{equation*}
Thus, for every $x\in V$,
\begin{equation*}
 \Phi_{\delta,\varepsilon}(x)
 =\sum\limits_{n=1}^{N-1}
 \chi_L\left(\frac{\mu(T_nx)}{r}\right)R_nx,
\end{equation*}
Let $e_{n,l}$ be the basis of $E_n$ dual to $(\xi_{n,l})_l$.
For $n<N$, choose $c_{n,l}\in\{0,1\}$ so that
$p_{n,l}=\xi_{n,l}\circ R_n+c_{n,l}$ is nonzero at $a$, and set
$t_n=1+Q\circ T_n$. With $\beta_r$ as in~\eqref{eq:ufdd-tail-coefficient},
define
\begin{equation*}
 B((s_{n,l}),(t_n))
 =\sum\limits_{n=1}^{N-1}\beta_r(t_n-1)
           \sum\limits_{l=1}^{\dim E_n}(s_{n,l}-c_{n,l})e_{n,l}.
\end{equation*}
The function $\beta_r$ is smooth on all of $\mathbb R$, so $B$ is smooth.
For $x\in V$,
\begin{equation*}
 \Phi_{\delta,\varepsilon}(x)=B((p_{n,l}(x))_{n<N,l},(t_n(x))_{n<N}).
\end{equation*}
All coordinates are nonzero at $a$, since $t_n(a)\geqslant1$.
The polynomials used belong to the fixed family and are independent of $r,L$.

\smallskip\noindent\proofstep{}\label{step:uniform-identity:2}
Suppose that $|I|>\aleph_0$, and assume that the proposition holds for every space with strictly fewer than $|I|$ blocks.
Retain the well-ordering of $I$ and the associated mappings $\Pi_i,T_i,T_{i+}$ and $K_{i,p}$.
For each $i\in I$, set $X_i=\Pi_iX$.
The family $(E_j)_{j\in I_{\leqslant i}}$ is an unconditional finite-dimensional decomposition of $X_i$.
By~\eqref{eq:well-order-small-parts}, the number of its blocks satisfies $|I_{\leqslant i}|<|I|$.
With the inherited norm and separating polynomial, the induction hypothesis gives a fixed polynomial family
$(\mathcal P_{i,k}:X_i\to H_{i,k})_{k\geqslant1}$.
For $p\in\mathbb Q[t]$, $k\geqslant1$, and an integer $D\geqslant0$, set
\begin{equation}\label{eq:ufdd-degree-group}
 I_{p,k,D}
 =\{i\in I:\deg(\mathcal P_{i,k}\circ K_{i,p})\leqslant D\}.
\end{equation}
For each $i\in I_{p,k,D}$, write the complex homogeneous decomposition
\begin{equation*}
 \mathcal P_{i,k}^{\mathbb C}\circ K_{i,p}^{\mathbb C}
 =\sum\limits_{\nu=0}^{D}P_{i,p,k}^{(\nu)},
 \qquad
 a_{i,p,k,D}
 =\frac{1}{1+\sum\limits_{\nu=0}^{D}\|P_{i,p,k}^{(\nu)}\|}>0.
\end{equation*}
Here polynomial norms are supremum norms on the complex unit ball. Then, for every $z\in X_{\mathbb C}$,
\begin{equation}\label{eq:ufdd-normalized-child-feature}
 \left\|a_{i,p,k,D}
 \mathcal P_{i,k}^{\mathbb C}(K_{i,p}^{\mathbb C}(z))\right\|
 \leqslant(1+\|z\|)^D.
\end{equation}
Define
\begin{equation}\label{eq:ufdd-gated-polynomial-array}
 H_{p,k,D}(x)
 =\left(
 a_{i,p,k,D}q_i(x)\mathcal P_{i,k}(K_{i,p}(x))
 \right)_{i\in I_{p,k,D}}.
\end{equation}
Regard $H_{p,k,D}$ as taking values in the Hilbert direct sum of the corresponding spaces $H_{i,k}$.
By
\eqref{eq:ufdd-block-summability} and \eqref{eq:ufdd-normalized-child-feature},
\begin{equation*}
 \|H_{p,k,D}(x)\|
 \leqslant(1+\|x\|)^D\|(q_i(x))_i\|_2
 \leqslant\|A\|(1+\|x\|)^D\|x\|^d.
\end{equation*}
Each scalar coordinate has degree at most $D+d$, so Lemma~\ref{lem:ufdd-polynomial-arrays}
shows that $H_{p,k,D}$ is a continuous Hilbert-valued polynomial.
There are only countably many triples $(p,k,D)$, so
\begin{equation}\label{eq:ufdd-universal-feature-family}
 \widehat J_0\quad\text{and all }H_{p,k,D}
\end{equation}
can be arranged into a countable family $\mathcal E_X$ independent of $\delta,\varepsilon$.

\smallskip\noindent\proofstep{}\label{step:uniform-identity:3}
Given $\delta,\varepsilon>0$, choose
\begin{equation*}
 \begin{gathered}
 \eta=\frac{\min\limits\{\varepsilon,1\}}{8},\qquad r=\frac{\delta}{2},\\
 L\in\mathbb N,\qquad
 L\geqslant\max\limits\left\{1,\left(\frac{2c_\sigma BT}{\eta}\right)^2\right\}.
 \end{gathered}
\end{equation*}
Then
\begin{equation*}
 \|D\Theta_{r,L}\|_\infty
 \leqslant1+\frac{2c_\sigma BT}{\sqrt L}\leqslant1+\eta.
\end{equation*}
Apply Corollary~\ref{cor:ufdd-quantitative-partition} to the cover \eqref{eq:ufdd-tail-chart-cover}
to obtain $(\varphi_{n,j})$ and constants $D_n$.
Using these constants $D_n$, set
\begin{equation}\label{eq:ufdd-en-definition}
 e_n=\min\limits\left\{\frac{\delta}{4},
 \frac{\eta2^{-n}}{1+D_n}\right\}>0.
\end{equation}
By Lemma~\ref{lem:ufdd-fixed-chart-polynomials}, for every $n,m$ we may choose
$p_{n,m}\in\mathbb Q[t]$ such that, for every $i$, the following inequalities hold on $V_{i,m}$:
\begin{equation}\label{eq:ufdd-local-polynomial-estimates}
 \|K_{i,p_{n,m}}-\Theta_{r,L}\|<\frac{e_n}{2},\qquad
 \|DK_{i,p_{n,m}}-D\Theta_{r,L}\|<\eta.
\end{equation}
By the induction hypothesis on $X_i$, choose $G_{i,n}:X_i\to X_i$ such that
\begin{equation}\label{eq:ufdd-child-approximation}
 \|G_{i,n}-I_{X_i}\|_\infty<\frac{e_n}{2},\qquad
 \|DG_{i,n}\|_\infty\leqslant1+\eta,
\end{equation}
with the locally finite-dimensional range and polynomial coordinate representation required by the induction hypothesis.
For each $(n,j)$, fix a member of the cover containing $\operatorname{supp}\varphi_{n,j}$.
If the chosen member is $V_{i,m}$, define the globally smooth mapping
\begin{equation*}
 F_{n,j}=G_{i,n}\circ K_{i,p_{n,m}}.
\end{equation*}
If the chosen member is $V_0$, take $F_{n,j}=0$.
For $x\in V_{i,m}$, write $K=K_{i,p_{n,m}}$. Then
\begin{equation*}
 \begin{aligned}
 \|F_{n,j}(x)-\Theta_{r,L}(x)\|
 &\leqslant\|G_{i,n}(K(x))-K(x)\|+\|K(x)-\Theta_{r,L}(x)\|<e_n,\\
 DF_{n,j}(x)&=DG_{i,n}(K(x))\circ DK(x),\\
 \|DF_{n,j}(x)\|
 &\leqslant(1+\eta)(\|D\Theta_{r,L}(x)\|+\eta)
 \leqslant(1+\eta)(1+2\eta).
 \end{aligned}
\end{equation*}
On $V_0$, both $F_{n,j}$ and $\Theta_{r,L}$ vanish. Hence, on the
corresponding member of the cover,
\begin{equation}\label{eq:ufdd-local-composition-bounds}
 \|F_{n,j}-\Theta_{r,L}\|<e_n,\qquad
 \|DF_{n,j}\|\leqslant(1+\eta)(1+2\eta).
\end{equation}

\smallskip\noindent\proofstep{}\label{step:uniform-identity:4}
Inspired by~\cite[proof of Theorem~14, pp.~575--576]{HajekJohanis2010}, set
\begin{equation}\label{eq:ufdd-glued-approximation}
 \Phi(x)=\sum\limits_{n,j}\varphi_{n,j}(x)F_{n,j}(x).
\end{equation}
By Lemma~\ref{lem:partition-summation}, $\Phi\in C^\infty(X,X)$.
Fix $a\in X$, and set $\Lambda_a=\{(n,j):a\in\operatorname{supp}\varphi_{n,j}\}$.
Choose an open neighborhood $O_0$ of $a$ meeting only finitely many supports, and denote their index set by $M_a$.
Then $\Lambda_a\subset M_a$, and
\begin{equation*}
 O_1=O_0\setminus\bigcup\limits_{(n,j)\in M_a\setminus\Lambda_a}
       \operatorname{supp}\varphi_{n,j}
\end{equation*}
is an open neighborhood of $a$ on which $\varphi_{n,j}|_{O_1}=0$ whenever $(n,j)\notin\Lambda_a$.
Let $\Lambda_a^*$ consist of the indices in $\Lambda_a$ assigned to sets
$V_{i,m}$. For each $\lambda=(n,j)\in\Lambda_a^*$, write
$K_\lambda=K_{i,p_{n,m}}$. The induction hypothesis gives an open set
$W_\lambda\subset X_i$ and a finite set $F_\lambda\subset I_{\leqslant i}$
such that
\begin{equation*}
 K_\lambda(a)\in W_\lambda,\qquad
 G_{i,n}(W_\lambda)\subset\bigoplus\limits_{l\in F_\lambda}E_l.
\end{equation*}
Set
\begin{equation*}
 V=O_1\cap\bigcap\limits_{\lambda\in\Lambda_a^*}K_\lambda^{-1}(W_\lambda),
 \qquad I_a=\bigcup\limits_{\lambda\in\Lambda_a^*}F_\lambda.
\end{equation*}
Then $V$ is an open neighborhood of $a$, $I_a$ is finite, and, for every $x\in V$,
\begin{equation*}
 \Phi(x)=\sum\limits_{(n,j)\in\Lambda_a}\varphi_{n,j}(x)F_{n,j}(x)
       \in\bigoplus\limits_{l\in I_a}E_l.
\end{equation*}

By \eqref{eq:ufdd-local-composition-bounds},
\begin{equation*}
 \|\Phi(x)-x\|
 \leqslant\|\Theta_{r,L}(x)-x\|
 +\sum\limits_{n,j}\varphi_{n,j}(x)
             \|F_{n,j}(x)-\Theta_{r,L}(x)\|
 \leqslant\frac{3\delta}{4}<\delta.
\end{equation*}
Since $\sum\limits_{n,j}D\varphi_{n,j}(x)=0$, for $h\in X$,
\begin{equation*}
 \begin{aligned}
 D\Phi(x)[h]
 &=\sum\limits_{n,j}\varphi_{n,j}(x)DF_{n,j}(x)[h]\\
 &\quad+\sum\limits_{n,j}D\varphi_{n,j}(x)[h]
                      (F_{n,j}(x)-\Theta_{r,L}(x)).
 \end{aligned}
\end{equation*}
The support conditions and~\eqref{eq:ufdd-local-composition-bounds} now give
\begin{equation}\label{eq:ufdd-global-derivative-bound}
 \begin{aligned}
 \sup\limits_{x\in X}\|D\Phi(x)\|
 &\leqslant(1+\eta)(1+2\eta)+\sum\limits_nD_ne_n\\
 &\leqslant(1+\eta)(1+2\eta)+\eta
 \leqslant1+6\eta\leqslant1+\varepsilon.
 \end{aligned}
\end{equation}

\smallskip\noindent\proofstep{}\label{step:uniform-identity:5}
Retain $a,V,\Lambda_a,\Lambda_a^*$ from Step~\ref{step:uniform-identity:4}.
Fix $\lambda=(n,j)\in\Lambda_a^*$, and let $V_{i,m}$ be its assigned member of the cover.
Put
\begin{equation*}
	K_\lambda=K_{i,p_{n,m}},\qquad y_\lambda=K_\lambda(a).
\end{equation*}
Since $a\in V_{i,m}$, we have $q_i(a)>0$.

By the induction hypothesis, there exist $d_\lambda\geqslant0$, an open neighborhood
$W'_\lambda$ of $y_\lambda$, an open set
$\Omega^G_\lambda\subset\mathbb R^{d_\lambda}$, and
$g^G_\lambda\in C^\infty(\Omega^G_\lambda,X_i)$ such that
\begin{equation*}
	\begin{gathered}
		r_{\lambda,\nu}(y)
		=\bigl(\mathcal P_{i,k_{\lambda,\nu}}(y)\bigr)_{\gamma_{\lambda,\nu}},
		\qquad r_{\lambda,\nu}(y_\lambda)\ne0
		\quad\text{for every }1\leqslant\nu\leqslant d_\lambda,\\
		r_\lambda=(r_{\lambda,1},\ldots,r_{\lambda,d_\lambda}),\qquad
		r_\lambda(W'_\lambda)\subset\Omega^G_\lambda,\\
		G_{i,n}(y)=g^G_\lambda(r_\lambda(y))
        \quad\text{for every }y\in W'_\lambda.
	\end{gathered}
\end{equation*}
For each $1\leqslant\nu\leqslant d_\lambda$, choose
\begin{equation*}
	\begin{gathered}
		D_{\lambda,\nu}\geqslant
		\deg(\mathcal P_{i,k_{\lambda,\nu}}\circ K_\lambda),\qquad
		\alpha_{\lambda,\nu}
		=a_{i,p_{n,m},k_{\lambda,\nu},D_{\lambda,\nu}}>0,\\
		h_{\lambda,\nu}(x)
		=\bigl(H_{p_{n,m},k_{\lambda,\nu},D_{\lambda,\nu}}(x)
		\bigr)_{i,\gamma_{\lambda,\nu}}.
	\end{gathered}
\end{equation*}
By~\eqref{eq:ufdd-gated-polynomial-array}, for every $x\in X$ with $q_i(x)>0$,
\begin{equation}\label{eq:ufdd-child-feature-recovery}
	r_{\lambda,\nu}(K_\lambda(x))
	=\frac{h_{\lambda,\nu}(x)}{\alpha_{\lambda,\nu}q_i(x)}.
\end{equation}
In particular,
\begin{equation*}
	h_{\lambda,\nu}(a)
	=\alpha_{\lambda,\nu}q_i(a)r_{\lambda,\nu}(y_\lambda)\ne0.
\end{equation*}
Now define
\begin{equation*}
	\begin{gathered}
		h_\lambda=(h_{\lambda,1},\ldots,h_{\lambda,d_\lambda}),\qquad
		Q_\lambda=(h_\lambda,q_i),\\
		R_\lambda(s,t)
		=\left(\frac{s_\nu}{\alpha_{\lambda,\nu}t}\right)_{\nu=1}^{d_\lambda}
		\quad\text{for every }(s,t)\in\mathbb R^{d_\lambda}\times(0,\infty),\\
		\Omega^F_\lambda=R_\lambda^{-1}(\Omega^G_\lambda),\qquad
		g^F_\lambda=g^G_\lambda\circ R_\lambda
		\in C^\infty(\Omega^F_\lambda,X_i),\\
		U^F_\lambda
		=V\cap K_\lambda^{-1}(W'_\lambda)\cap q_i^{-1}((0,\infty)).
	\end{gathered}
\end{equation*}
Since $R_\lambda$ is smooth and~\eqref{eq:ufdd-child-feature-recovery} holds,
\begin{equation*}
	\begin{gathered}
		a\in U^F_\lambda,\qquad
		Q_\lambda(U^F_\lambda)\subset\Omega^F_\lambda,\\
		F_{n,j}(x)=g^F_\lambda(Q_\lambda(x))
		\quad\text{for every }x\in U^F_\lambda.
	\end{gathered}
\end{equation*}
By~\eqref{eq:ufdd-baseline-embedding} and
\eqref{eq:ufdd-universal-feature-family}, every component of $Q_\lambda$
is a scalar coordinate of $\mathcal E_X$ and is nonzero at $a$.
For $\lambda=(n,j)\in\Lambda_a\setminus\Lambda_a^*$, we have $F_{n,j}=0$.

On the other hand, Corollary~\ref{cor:ufdd-quantitative-partition}
gives the same type of local coordinate representation for each
$\varphi_{n,j}$. Consequently, for every $\lambda=(n,j)\in\Lambda_a$
there exist an integer $N_\lambda\geqslant0$, an open neighborhood
$U_\lambda$ of $a$, an open set
$\Omega_\lambda\subset\mathbb R^{N_\lambda}$, scalar coordinates
$p_{\lambda,1},\ldots,p_{\lambda,N_\lambda}$ of $\mathcal E_X$,
and $g_\lambda\in C^\infty(\Omega_\lambda,\mathbb R\times X)$ such that
\begin{equation*}
	\begin{gathered}
		P_\lambda=(p_{\lambda,1},\ldots,p_{\lambda,N_\lambda}),\qquad
		P_\lambda(U_\lambda)\subset\Omega_\lambda,\\
		p_{\lambda,\nu}(a)\ne0
		\quad\text{for every }1\leqslant\nu\leqslant N_\lambda,\\
		\bigl(\varphi_{n,j}(x),F_{n,j}(x)\bigr)
		=g_\lambda(P_\lambda(x))\quad\text{for every }x\in U_\lambda.
	\end{gathered}
\end{equation*}

Since $\Lambda_a$ is finite, set
\begin{equation*}
	\begin{gathered}
		N=\sum\limits_{\lambda\in\Lambda_a}N_\lambda,\qquad
		P=(P_\lambda)_{\lambda\in\Lambda_a}=(p_1,\ldots,p_N),\\
		\Omega=\prod\limits_{\lambda\in\Lambda_a}\Omega_\lambda
		\subset\mathbb R^N,\qquad
		V'=V\cap\bigcap\limits_{\lambda\in\Lambda_a}U_\lambda.
	\end{gathered}
\end{equation*}
For $\lambda=(n,j)\in\Lambda_a$, let
$\pi_\lambda:\mathbb R^N\to\mathbb R^{N_\lambda}$
be the corresponding coordinate projection, and define
\begin{equation*}
	(b_{n,j},c_{n,j})
	=g_\lambda\circ(\pi_\lambda|_\Omega)
	\in C^\infty(\Omega,\mathbb R\times X).
\end{equation*}
Then $a\in V'$, $P(V')\subset\Omega$, and $p_\nu(a)\ne0$ for every $1\leqslant\nu\leqslant N$.
Moreover, for every $(n,j)\in\Lambda_a$ and $x\in V'$,
\begin{equation*}
	\varphi_{n,j}(x)=b_{n,j}(P(x)),\qquad
	F_{n,j}(x)=c_{n,j}(P(x)).
\end{equation*}
Therefore, defining $B$ for every $t\in\Omega$ by
\begin{equation*}
	B(t)=\sum\limits_{(n,j)\in\Lambda_a}b_{n,j}(t)c_{n,j}(t),
\end{equation*}
we obtain
\begin{equation*}
	B\in C^\infty(\Omega,X),\qquad
	\Phi(x)=B(p_1(x),\ldots,p_N(x))\quad\text{for every }x\in V'.
\end{equation*}
This is the representation required by
Definition~\ref{def:ufdd-polynomial-features} relative to
\eqref{eq:ufdd-universal-feature-family}.

The resulting $\Phi$ satisfies every requirement of the inductive assertion,
completing the cardinal induction.

\end{proof}

\subsection[Approximation on decomposition spaces]{\texorpdfstring{$C^1$}{C1} real-analytic approximation on decomposition spaces}\label{subsec:ufdd-entire-embedding}

\begin{lemma}\label{lem:ufdd-entire-compilation}
Let $X$ be an arbitrary real Banach space, and let $(H_n)_{n\geqslant1}$ be a family of real Hilbert spaces,
each equipped with fixed orthonormal coordinates. Let $\mathcal P_n:X\to H_n$ be continuous polynomials.
Suppose that there exist $n_0\geqslant1$, a real Hilbert space $H_0$,
a bounded linear mapping $L_0:H_{n_0}\to H_0$, and a polynomial topological embedding $J_0:X\to H_0$ such that
\begin{equation*}
 J_0=L_0\circ\mathcal P_{n_0}.
\end{equation*}
Then there exists an entire real-analytic topological embedding
\begin{equation*}
 J:X\longrightarrow H=\bigoplus\limits_{n\geqslant1}^{\ell^2}H_n,
\end{equation*}
and a family of positive numbers $(w_n)_{n\geqslant1}$ such that $\pi_n\circ J=w_n\mathcal P_n$ for every $n$,
where $\pi_n:H\to H_n$ is the canonical projection onto the $n$th summand.
Moreover, let $U\subset X$ be open, let $Z$ be a real Banach space, and let $k\in\{1,\infty\}$.
Suppose that $u:U\to Z$ and that, for every $a\in U$, there exist
an integer $N\geqslant0$, an open neighborhood $V\subset U$ of $a$,
an open set $O\subset\mathbb R^N$, a mapping $b\in C^k(O,Z)$,
and scalar coordinates $p_1,\ldots,p_N$ of the polynomials
$\mathcal P_n$ such that $(p_1,\ldots,p_N)(V)\subset O$,
$p_j(a)\ne0$ for every $1\leqslant j\leqslant N$, and
\begin{equation*}
 u(x)=b(p_1(x),\ldots,p_N(x))
 \qquad\text{for every }x\in V.
\end{equation*}
Then $u\in\mathscr F_J^k(U,Z)$.
In particular, if $u:X\to Z$ has a representation as in
Definition~\ref{def:ufdd-polynomial-features} relative to
$(\mathcal P_n)_{n\geqslant1}$, then, near every $a\in X$, it can be
written as a smooth function of finitely many scalar coordinates of $J$
whose values at $a$ are nonzero.
\end{lemma}
\begin{proof}
Write $d_n=\deg\mathcal P_n$, and take the finite homogeneous expansion
\begin{equation*}
 \mathcal P_n^{\mathbb C}(z)=\sum\limits_{j=0}^{d_n}P_n^{(j)}(z),\qquad
 C_n=1+\sum\limits_{j=0}^{d_n}\|P_n^{(j)}\|,
\end{equation*}
where each $P_n^{(j)}$ is a continuous homogeneous polynomial of degree $j$, possibly zero.
Polynomial norms are supremum norms on the complex unit ball. Thus
\begin{equation*}
 \|\mathcal P_n^{\mathbb C}(z)\|
 \leqslant\sum\limits_{j=0}^{d_n}\|P_n^{(j)}\|
                  \max\limits\{1,\|z\|\}^j
 \leqslant C_n\max\limits\{1,\|z\|\}^{d_n}.
\end{equation*}
Set
\begin{equation}\label{eq:ufdd-entire-feature-embedding}
 w_n=\frac{2^{-n}}{C_n n^{d_n}},\qquad
 J^{\mathbb C}(z)=
 \bigl(w_n\mathcal P_n^{\mathbb C}(z)\bigr)_{n\geqslant1}.
\end{equation}
For $n\geqslant\max\limits\{1,R\}$,
\begin{equation*}
 \sup\limits_{\|z\|\leqslant R}
 \|w_n\mathcal P_n^{\mathbb C}(z)\|
 \leqslant2^{-n}\left(\frac{\max\limits\{1,R\}}{n}\right)^{d_n}
 \leqslant2^{-n}.
\end{equation*}
For every integer $N\geqslant\max\limits\{1,R\}$,
\begin{equation*}
 \sup\limits_{\|z\|\leqslant R}
 \left\|J^{\mathbb C}(z)
       -\sum\limits_{n=1}^N\iota_nw_n\mathcal P_n^{\mathbb C}(z)\right\|^2
 \leqslant\sum\limits_{n>N}4^{-n}=\frac{4^{-N}}{3},
\end{equation*}
where $\iota_n:H_{n,\mathbb C}\to H_{\mathbb C}$ is the canonical inclusion.
Thus the finite sums converge uniformly on every complex ball.
By the locally uniform limit theorem for Banach-valued holomorphic mappings
(\cite[Proposition~9.13, p.~72]{Mujica1986}), $J^{\mathbb C}$ is entire.
Write $J=J^{\mathbb C}|_X$, and let $\pi_{n_0}:H\to H_{n_0}$ be the coordinate projection.
Set $L=w_{n_0}^{-1}L_0\pi_{n_0}$. Then $L$ is bounded and linear, and $J_0=L\circ J$.
Thus $J$ is injective and $J^{-1}=J_0^{-1}\circ L|_{J(X)}$ is continuous.
For every scalar coordinate polynomial $p=(\mathcal P_n)_\gamma$, we have
\begin{equation*}
 p(x)=\frac{(J(x))_{n,\gamma}}{w_n},\qquad
 p(a)\ne0\quad\Longleftrightarrow\quad(J(a))_{n,\gamma}\ne0.
\end{equation*}
Write $p_j=(\mathcal P_{n_j})_{\gamma_j}$ and
$S=\{(n_j,\gamma_j):1\leqslant j\leqslant N\}$.
Define a linear mapping $L_S:\mathbb R^S\to\mathbb R^N$ by
\begin{equation*}
 (L_St)_j=\frac{t_{n_j,\gamma_j}}{w_{n_j}},\qquad
 O_S=L_S^{-1}(\operatorname{dom}b),\qquad B=b\circ L_S|_{O_S}.
\end{equation*}
Then $B\in C^k(O_S,Z)$, $S\subset\operatorname{supp}J(a)$, and
$u(x)=B((J(x))_S)$ on some open neighborhood $V_a\subset U$ of $a$.
Choose an open $O_a\subset H$ such that $J(a)\in O_a$ and
$J^{-1}(O_a)\subset V_a$, and set
\begin{equation*}
 \begin{gathered}
 W=O_a\cap\{y:y_S\in O_S\}
       \cap\bigcap\limits_{\gamma\in S}\{y:y_\gamma\ne0\},\\
 \widetilde u(y)=B(y_S)\quad\text{for every }y\in W.
 \end{gathered}
\end{equation*}
For each $y\in W$, the set $S$ lies in $\operatorname{supp}y$.
The outer mapping $(s,t)\mapsto B(s)$ on $O_S\times\mathbb R$ gives
$\widetilde u\in\mathscr R_H^k(W,Z)$. Finally,
\begin{equation*}
 J^{-1}(W)\subset U,\qquad
 \widetilde u(J(x))=B((J(x))_S)=u(x)
 \quad\text{for every }x\in J^{-1}(W).
\end{equation*}
This is the local extension representation required by
Definition~\ref{def:liftable-class}.
\end{proof}

\begin{proof}[Proof of Theorem~\ref{thm:ufdd}]
If $X=\{0\}$, take the zero mapping.
If $X\ne\{0\}$, apply Lemma~\ref{lem:ufdd-entire-compilation} to
the fixed polynomial family in Proposition~\ref{prop:ufdd-strengthened-approximation} to obtain a common embedding $J$.
Set $\Phi_n=\Phi_{2^{-n},1}$. Then~\eqref{eq:ufdd-final-static-bounds} holds.
The local range condition involving finitely many blocks follows from Proposition~\ref{prop:ufdd-strengthened-approximation}.
Lemma~\ref{lem:ufdd-entire-compilation} rewrites the representation in nonzero polynomial coordinates as~\eqref{eq:ufdd-final-active-lift}
and gives $\Phi_n\in\mathscr F_J^\infty(X,X)$.
\end{proof}

If each $X_1,\ldots,X_N$ admits a separating polynomial and an unconditional finite-dimensional decomposition,
choose on $X_i$ a positive homogeneous separating polynomial $Q_i$ of degree $d_i$, and let $d$ be a common multiple of the $d_i$. Then
\begin{equation}\label{eq:finite-sum-separator}
 Q(x_1,\ldots,x_N)=\sum\limits_{i=1}^N Q_i(x_i)^{\frac{d}{d_i}}
\end{equation}
separates the unit sphere for the maximum direct sum norm after normalizing
each $Q_i$ so that $Q_i(x_i)\geqslant\|x_i\|^{d_i}$. Indeed,
\begin{equation*}
 Q(x_1,\ldots,x_N)\geqslant\sum\limits_{i=1}^N\|x_i\|^d
 \geqslant\max\limits_{1\leqslant i\leqslant N}\|x_i\|^d,
 \qquad Q(0)=0.
\end{equation*}
For a decomposition $(E_{i,\gamma})_{\gamma\in I_i}$ of $X_i$ with projections
$R_{i,\gamma}$, use the index set
$\{(i,\gamma):1\leqslant i\leqslant N,\ \gamma\in I_i\}$ and the projections
\begin{equation*}
 \widehat R_{i,\gamma}(x_1,\ldots,x_N)
 =(0,\ldots,0,R_{i,\gamma}x_i,0,\ldots,0).
\end{equation*}
The sum of these projections applied to $x$ is $x$, unconditionally,
and each range is finite-dimensional.

\begin{proof}[Proof of Theorem~\ref{thm:main-complete}]
Let $X$ be a complemented closed subspace of a space $E$ admitting a separating polynomial and an unconditional finite-dimensional decomposition.
By Theorem~\ref{thm:ufdd} and Remark~\ref{rem:equivalent-norm-bounds},
the space $E$, with its original norm, satisfies the hypotheses of Theorem~\ref{thm:criterion}. Thus every Banach-valued $C^1$ mapping on an arbitrary open subset of $E$
admits fine $C^1$ real-analytic approximation.
Proposition~\ref{prop:permanence} transfers this property to $X$, giving~\eqref{eq:main-approximation}.
\end{proof}

\section[\texorpdfstring{$C^1$}{C1} Approximation on Separable Domains]{Fine \texorpdfstring{$C^1$}{C1} Real-Analytic Approximation on Separable Domains}\label{sec:target-geometry}

\subsection[Approximation through local polynomial representations]{Real-analytic approximation of mappings with local polynomial composition representations}\label{subsec:target-polynomial-approximation}

The definition is inspired by
\cite[Chapter~5, Definition~78, Fact~79, and Lemma~81, pp.~279--280]{HajekJohanis2014}.

\begin{definition}\label{def:target-polynomial-class}
Let $X,Z$ be real Banach spaces, and let $U\subset X$ be open.
Suppose that $u\in C^\infty(U,Z)$ and that there is a countable family of continuous real polynomials
$(p_n:X\to\mathbb R)_{n\geqslant1}$ such that, for every $a\in U$, one can choose finitely many indices
$n_1,\ldots,n_N$, an open neighborhood $a\in V\subset U$, an open set $O\subset\mathbb R^N$,
and $b\in C^\infty(O,Z)$ such that $(p_{n_1},\ldots,p_{n_N})(V)\subset O$ and, for every $x\in V$,
\begin{equation}\label{eq:target-polynomial-class}
 u(x)=b\bigl(p_{n_1}(x),\ldots,p_{n_N}(x)\bigr).
\end{equation}
Then we write $u\in\mathscr S(U,Z)$.
The case $N=0$ is allowed and means that $u|_V$ is constant.
In the scalar case, write $\mathscr S(U)=\mathscr S(U,\mathbb R)$.
\end{definition}

\begin{lemma}\label{lem:target-class-closure}
The class $\mathscr S$ is stable under restriction to open subsets, finite products, continuous multilinear operations,
smooth postcomposition, and countable locally finite sums.
If $V\subset U$ is open, $\psi\in\mathscr S(U)$, $v\in\mathscr S(V,Z)$,
and $\operatorname{supp}_U\psi\subset V$, then the extension of $\psi v$ by zero belongs to $\mathscr S(U,Z)$.
Moreover, for every continuous affine mapping $L:E\to X$ and
$u\in\mathscr S(U,Z)$, we have $u\circ L\in\mathscr S(L^{-1}(U),Z)$.
\end{lemma}
\begin{proof}
Restriction to an open subset uses the same polynomial family and the
intersections of the original representation neighborhoods with that subset.
For finitely many mappings $u_j\in\mathscr S(U,Z_j)$, fix $a\in U$ and write
\begin{equation*}
 u_j|_{V_j}=b_j\circ P_j|_{V_j},\qquad
 P_j:X\to\mathbb R^{N_j},\qquad P_j(V_j)\subset O_j,
\end{equation*}
where the components of $P_j$ belong to the countable polynomial family for
$u_j$. Set
\begin{equation*}
 \begin{gathered}
 P=(P_1,\ldots,P_m),\qquad V=\bigcap\limits_{j=1}^mV_j,\qquad
 O=\prod\limits_{j=1}^mO_j,\\
 B(t_1,\ldots,t_m)=(b_1(t_1),\ldots,b_m(t_m)).
 \end{gathered}
\end{equation*}
Then $P(V)\subset O$, $B\in C^\infty(O,\prod\limits_j Z_j)$, and
$(u_1,\ldots,u_m)|_V=B\circ P|_V$.
For a continuous multilinear map $M:Z_1\times\cdots\times Z_m\to Z$, use
\begin{equation*}
 B_M(t_1,\ldots,t_m)=M(b_1(t_1),\ldots,b_m(t_m)),\qquad
 M(u_1,\ldots,u_m)|_V=B_M\circ P|_V.
\end{equation*}
For smooth postcomposition by $\Phi:W\to Z'$, where $W$ is open and
$u(U)\subset W$, a representation $u=b\circ P$ on $V$ gives
\begin{equation*}
 O'=O\cap b^{-1}(W),\qquad P(V)\subset O',\qquad
 (\Phi\circ u)|_V=(\Phi\circ b|_{O'})\circ P|_V.
\end{equation*}
These constructions use a finite union of countable polynomial families.

Let $(u_i)_{i\geqslant1}$ be a locally finite family in $\mathscr S(U,Z)$,
and denote its polynomial families by $\mathcal P_i$.
The set $\mathcal P=\bigcup\limits_{i\geqslant1}\mathcal P_i$ is countable.
For each $a\in U$, choose an open $V_a\ni a$ and a finite $I_a$ such that
$u_i|_{V_a}=0$ for $i\notin I_a$.
The preceding product representation gives, after restricting $V_a$,
\begin{equation*}
 \begin{gathered}
 u_i(x)=b_i(P_i(x))\quad\text{for every }i\in I_a,\\
 B_\Sigma((t_i)_{i\in I_a})=\sum\limits_{i\in I_a}b_i(t_i),\\
 \sum\limits_{i\geqslant1}u_i(x)=B_\Sigma((P_i(x))_{i\in I_a}).
 \end{gathered}
\end{equation*}
Every scalar component of each $P_i$ belongs to the fixed
family $\mathcal P$, which is independent of $a$.

For extension by zero, define
\begin{equation*}
 w(x)=\begin{cases}\psi(x)v(x),&x\in V,\\0,&x\in U\setminus V.\end{cases}
\end{equation*}
On $V$, use the local product representation just proved; on
$U\setminus\operatorname{supp}_U\psi$, use $w=0$.
These two open sets cover $U$, and the polynomials belong to the union of
the families for $\psi$ and $v$. Hence $w\in\mathscr S(U,Z)$.

Finally, let $L:E\to X$ be continuous and affine.
If $u=b\circ(p_{n_1},\ldots,p_{n_N})$ on $V\subset U$, then
\begin{equation*}
 (u\circ L)|_{L^{-1}(V)}
 =b\circ(p_{n_1}\circ L,\ldots,p_{n_N}\circ L)|_{L^{-1}(V)}.
\end{equation*}
The sets $L^{-1}(V)$ are open and the polynomials $p_n\circ L$ form a
countable family, proving the last assertion.
\end{proof}

\begin{proposition}\label{prop:target-polynomial-analytic}
Let $X$ be a real Banach space admitting a separating polynomial, and let $U\subset X$ be open.
For every real Banach space $Z$, each $u\in\mathscr S(U,Z)$ admits fine $C^1$ real-analytic approximation.
\end{proposition}
\begin{proof}
By Theorem~\ref{thm:general-polynomial-hilbert-embedding}, choose a polynomial topological embedding
$J_0:X\to H_0$. Let $(p_n)$ be the countable polynomial family used in the local representations of $u$,
and regard $J_0$ together with all $p_n,1+p_n$ as a countable family of Hilbert-valued polynomials.
Lemma~\ref{lem:ufdd-entire-compilation} gives a common entire real-analytic topological embedding
$J:X\to H$ such that each listed polynomial $P$ can be written
as $P=c\,\pi\circ J$, where $c>0$ and $\pi$ is the canonical
projection onto the summand corresponding to $P$.

For a representation $u=b\circ(p_{n_1},\ldots,p_{n_N})$ on $V\ni a$,
with $b\in C^\infty(O,Z)$, define
\begin{equation*}
 \begin{gathered}
 c_j=\begin{cases}1,&p_{n_j}(a)=0,\\0,&p_{n_j}(a)\ne0,\end{cases}
 \qquad \widehat p_j=p_{n_j}+c_j,\qquad c=(c_1,\ldots,c_N),\\
 \widehat O=O+c,\qquad \widehat b(t)=b(t-c)
 \quad\text{for every }t\in\widehat O.
 \end{gathered}
\end{equation*}
Then $\widehat p_j(a)\ne0$, $(\widehat p_1,\ldots,\widehat p_N)(V)\subset\widehat O$,
and $u|_V=\widehat b\circ(\widehat p_1,\ldots,\widehat p_N)|_V$.
The last assertion of Lemma~\ref{lem:ufdd-entire-compilation} therefore gives
\begin{equation}\label{eq:target-polynomial-liftability}
 u\in\mathscr F_J^\infty(U,Z).
\end{equation}
Proposition~\ref{prop:analytic-liftable} gives fine $C^1$ real-analytic approximation of $u$.
\end{proof}

\subsection[Approximation of scalar Lipschitz functions]{Approximation of scalar Lipschitz functions in \texorpdfstring{$\mathscr S(X)$}{S(X)}}\label{subsec:target-scalar-approximation}

For the remainder of this section, fix a separable real Banach space $X$ admitting a separating polynomial.
All assertions hold for an empty domain. If $X=\{0\}$, the original mapping itself is the required real-analytic approximation.
Henceforth assume that $X\ne\{0\}$ and that the domain is nonempty.
By Lemma~\ref{lem:homogeneous-separator}, fix an even integer $d\geqslant2$, a continuous symmetric $d$-linear form
$\mathcal Q$, and a constant $C_Q\geqslant1$ such that
\begin{equation}\label{eq:target-fixed-separator}
 Q(x)=\mathcal Q(x^d),\qquad
 \|x\|^d\leqslant Q(x)\leqslant C_Q\|x\|^d.
\end{equation}

\begin{lemma}\label{lem:target-countable-partition}
For every open set $U\subset X$ and every open cover of $U$, there is a subordinate countable locally finite partition of unity
$(\psi_i)_{i\geqslant1}\subset\mathscr S(U)$ satisfying
\begin{equation*}
 0\leqslant\psi_i\leqslant1,
 \qquad M_i:=\sup\limits_{x\in U}\|D\psi_i(x)\|<\infty.
\end{equation*}
The relative closed support of each member is contained in its assigned member of the cover.
\end{lemma}
\begin{proof}
Given $a\in U$ and a member $V_a$ of the cover containing $a$, choose $r_a>0$ such that
$\overline B_X(a,r_a)\subset V_a$.
Set $\zeta(t)=1-\sigma(2t-1)$, using~\eqref{ns:eq:smooth-step-function}.
Then $\zeta\in C^\infty(\mathbb R,[0,1])$, $\zeta=1$ on
$(-\infty,\frac{1}{2}]$, and $\zeta=0$ on $[1,\infty)$.
Choose $0<t_a<r_a^d$ and set $\beta_a(x)=\zeta(\frac{Q(x-a)}{t_a})$.
By~\eqref{eq:target-fixed-separator}, $\beta_a=1$ on
$B_X(a,(\frac{t_a}{2C_Q})^{\frac{1}{d}})$, and
\begin{equation*}
 \operatorname{supp}\beta_a
 \subset\{x:Q(x-a)\leqslant t_a\}
 \subset\overline B_X(a,t_a^{\frac{1}{d}})\subset V_a.
\end{equation*}
By the chain rule,
\begin{equation}\label{eq:target-partition-cutoff-derivative}
 \begin{gathered}
 D\beta_a(x)=\frac{1}{t_a}\zeta'\left(\frac{Q(x-a)}{t_a}\right)DQ(x-a),\\
 \|D\beta_a\|_\infty
 \leqslant d\|\mathcal Q\|\,\|\zeta'\|_\infty\,t_a^{-\frac{1}{d}}<\infty,
 \end{gathered}
\end{equation}
where we have used $D\beta_a(x)\ne0\Longrightarrow\|x-a\|\leqslant t_a^{\frac{1}{d}}$.

The open neighborhoods on which the functions $\beta_a$ equal one cover $U$.
By separability, choose a countable subcover, and denote the corresponding functions by $\beta_i$.
If the subcover is finite, repeat one of its members to extend both the functions and their assigned cover members to sequences.
Apply Lemma~\ref{lem:successive-partition} to $(\beta_i)$,
and use~\eqref{eq:successive-partition} to obtain a countable locally finite partition of unity $(\psi_i)$.
Since $\operatorname{supp}_U\psi_i\subset\operatorname{supp}\beta_i\subset V_{a_i}$, this partition is subordinate to the given cover.
By~\eqref{eq:successive-partition-derivative} and~\eqref{eq:target-partition-cutoff-derivative},
$\|D\psi_i\|_\infty\leqslant\sum\limits_{j\leqslant i}\|D\beta_j\|_\infty<\infty$.
Each member uses only finitely many translated polynomials $Q(x-a_j)$, and the family of all polynomials used is countable.
Hence $\psi_i\in\mathscr S(U)$.
\end{proof}

\begin{lemma}\label{lem:target-scalar-lipschitz}
There exists $K_X\geqslant1$ such that, for every $L$-Lipschitz function $F:X\to\mathbb R$ and every
$\eta>0$, there is $v\in\mathscr S(X)$ satisfying
\begin{equation}\label{eq:target-scalar-lipschitz-result}
 \sup\limits_{x\in X}|v(x)-F(x)|<\eta,
 \qquad
 \sup\limits_{x\in X}\|Dv(x)\|\leqslant K_XL.
\end{equation}
The constant $K_X$ is independent of $F,L,\eta$.
\end{lemma}
\begin{proof}
\smallskip\noindent\firstproofstep{}\label{step:polynomial-smoothing:1}
Inspired by~\cite[Lemma~2, pp.~131--133]{AFK2012},
set $q(x)=(1+Q(x))^{\frac{1}{d}}-1$.
By~\eqref{eq:target-fixed-separator}, we have $q\geqslant0$, $q(0)=0$, and
\begin{equation*}
 \begin{aligned}
 Dq(x)[h]&=(1+Q(x))^{\frac{1}{d}-1}\mathcal Q(x^{d-1},h),\\
 \|Dq(x)\|&\leqslant
 \frac{\|\mathcal Q\|\|x\|^{d-1}}{(1+\|x\|^d)^{\frac{d-1}{d}}}
 \leqslant\|\mathcal Q\|.
 \end{aligned}
\end{equation*}
Thus $q$ is globally Lipschitz, and
\begin{equation}\label{eq:target-q-small-ball}
 q(x)<4\quad\Longrightarrow\quad\|x\|<5.
\end{equation}

Choose a dense sequence $(a_n)$ in $X$, and set $q_n(x)=q(80(x-a_n))$.
Set $\alpha(t)=1-\sigma(t-3)$ and $\beta(t)=\sigma(2t-3)$.
Then $\alpha,\beta\in C^\infty(\mathbb R,[0,1])$ and
\begin{equation*}
 \begin{aligned}
 &\alpha=1\text{ on }(-\infty,3],\qquad
 \alpha=0\text{ on }[4,\infty),\\
 &\beta=0\text{ on }(-\infty,\frac{3}{2}],\qquad
 \beta=1\text{ on }[2,\infty).
 \end{aligned}
\end{equation*}
Inspired by~\cite[proof of the unnumbered theorem, p.~142]{Johanis2003},
define
\begin{equation*}
 \begin{gathered}
 \kappa(t)=c_\kappa\beta_0(1-16t^2),\qquad
 c_\kappa=\left(\int\limits_{-\frac{1}{4}}^{\frac{1}{4}}
                      \beta_0(1-16t^2)\,\mathrm{d}t\right)^{-1},\\
 \kappa_N(s_1,\ldots,s_N)=\prod\limits_{j=1}^N\kappa(s_j).
 \end{gathered}
\end{equation*}
Here $\beta_0$ is defined in~\eqref{ns:eq:smooth-step-function}.
Then $\kappa_N\in C_c^\infty(\mathbb R^N)$ satisfies
\begin{equation*}
 \kappa_N\geqslant0,\qquad
 \operatorname{supp}\kappa_N\subset[-\frac{1}{4},\frac{1}{4}]^N,\qquad
 \int\limits_{\mathbb R^N}\kappa_N(s)\,\mathrm{d}s=1.
\end{equation*}
Define
\begin{equation*}
 m_N(t)=\int\limits_{\mathbb R^N}
       \min\limits_{1\leqslant j\leqslant N}(t_j-s_j)\kappa_N(s)\,\mathrm{d}s.
\end{equation*}
Write $m(t)=\min\limits_{1\leqslant j\leqslant N}t_j$.
The change of variables $u=t-s$ gives, for every multi-index $\gamma$,
\begin{equation*}
 m_N(t)=\int\limits_{\mathbb R^N}m(u)\kappa_N(t-u)\,\mathrm{d}u,\qquad
 \partial^\gamma m_N(t)
 =\int\limits_{\mathbb R^N}m(u)\partial^\gamma\kappa_N(t-u)\,\mathrm{d}u.
\end{equation*}
For $t$ in a bounded neighborhood, the integrands have a common compact
support and $m$ is bounded there. Thus differentiation is justified at every
order and $m_N\in C^\infty(\mathbb R^N)$.
Since $|m(t)-m(t')|\leqslant\|t-t'\|_\infty$,
\begin{equation*}
 \begin{aligned}
 |m_N(t)-m(t)|
 &\leqslant\int\limits_{\mathbb R^N}\|s\|_\infty\kappa_N(s)\,\mathrm{d}s
 \leqslant\frac{1}{4},\\
 |m_N(t)-m_N(t')|
 &\leqslant\int\limits_{\mathbb R^N}
             |m(t-s)-m(t'-s)|\kappa_N(s)\,\mathrm{d}s
 \leqslant\|t-t'\|_\infty.
 \end{aligned}
\end{equation*}
Consequently,
\begin{equation}\label{eq:target-smooth-minimum}
 \left|m_N(t)-\min\limits_{1\leqslant j\leqslant N}t_j\right|\leqslant\frac{1}{4},
 \qquad
 \operatorname{Lip}_{\|\cdot\|_\infty}(m_N)\leqslant1.
\end{equation}
Inspired by~\cite[Lemma~1 and the subsequent construction, pp.~127--128]{Fry2004},
\cite[Section~2, Fact~2 and Theorem~3, pp.~565--569]{HajekJohanis2010},
and~\cite[Chapter~7, proof of Theorem~64, pp.~427--428]{HajekJohanis2014}, define
\begin{equation}\label{eq:target-sup-partition}
 \varphi_1=\alpha(q_1),\qquad
 \varphi_n=\alpha(q_n)\,
 \beta\bigl(m_{n-1}(q_1,\ldots,q_{n-1})\bigr)
 \quad\text{for every }n\geqslant2.
\end{equation}
Each $\varphi_n$ is a smooth function of finitely many polynomials $Q(80(x-a_j))$.

Fix $x\in X$, and choose the least index $n$ such that $q_n(x)<3$.
If $n>1$, then
\begin{equation*}
 \min\limits_{j<n}q_j(x)\geqslant3,\qquad
 m_{n-1}(q_1(x),\ldots,q_{n-1}(x))\geqslant\frac{11}{4},
\end{equation*}
so $\varphi_n(x)=1$. If $n=1$, then $q_1(x)<3$, and hence $\varphi_1(x)=\alpha(q_1(x))=1$.
Thus $\sup\limits_{n\geqslant1}\varphi_n(x)=1$.
We may also choose $k$ such that $q_k(x)<1$. By continuity of $q_k$, the set $V=\{y\in X:q_k(y)<1\}$ is an open neighborhood of $x$.
For $n>k$, \eqref{eq:target-smooth-minimum} gives
$m_{n-1}(q_1,\ldots,q_{n-1})<\frac{5}{4}$ on $V$, so $\varphi_n|_V=0$.
Thus $(\varphi_n)$ is locally finite. For $n\geqslant2$, set
$Q_{<n}(x)=(q_1(x),\ldots,q_{n-1}(x))$, with the maximum norm on its range.
Then $\|DQ_{<n}(x)\|\leqslant80\|\mathcal Q\|$ and
\begin{equation*}
 \begin{aligned}
 D\varphi_n(x)[h]
 &=\alpha'(q_n(x))Dq_n(x)[h]\,\beta(m_{n-1}(Q_{<n}(x)))\\
 &\quad+\alpha(q_n(x))\beta'(m_{n-1}(Q_{<n}(x)))
             Dm_{n-1}(Q_{<n}(x))[DQ_{<n}(x)[h]].
 \end{aligned}
\end{equation*}
The norm of $Dm_{n-1}$ as a functional on the maximum-norm space is at most
one by~\eqref{eq:target-smooth-minimum}. Therefore
\begin{equation*}
 \|D\varphi_n(x)\|
 \leqslant80\|\mathcal Q\|(\|\alpha'\|_\infty+\|\beta'\|_\infty).
\end{equation*}
For $n=1$, use $D\varphi_1=\alpha'(q_1)Dq_1$.
Together with~\eqref{eq:target-q-small-ball}, this proves
\begin{equation}\label{eq:target-sup-partition-bounds}
 \varphi_n(x)\ne0\quad\Longrightarrow\quad
 \|x-a_n\|<\frac{1}{16},
 \qquad \operatorname{Lip}(\varphi_n)\leqslant M,
\end{equation}
where $M=80\|\mathcal Q\|(\|\alpha'\|_\infty+\|\beta'\|_\infty)$ is independent of $n$.

\smallskip\noindent\proofstep{}\label{step:polynomial-smoothing:2}
Inspired by~\cite[proof of $(i)\Rightarrow(iii)$ in Theorem~3, p.~567]{HajekJohanis2010},
choose $\chi_0\in C^\infty(\mathbb R,[0,1])$ vanishing on $(-\infty,\frac{1}{4}]$
and equal to one on $[1,\infty)$, and set $h=\chi_0\circ Q\in\mathscr S(X)$.
Then $h$ vanishes on a neighborhood of the origin and equals one on $X\setminus B_X(0,1)$.
Its derivative can be nonzero only on $\{Q\leqslant1\}\subset\overline B_X(0,1)$.
In fact,
\begin{equation*}
 Dh(x)[v]=\chi_0'(Q(x))d\mathcal Q(x^{d-1},v),\qquad
 \|Dh\|_\infty\leqslant d\|\mathcal Q\|\|\chi_0'\|_\infty.
\end{equation*}
Thus $C_h:=\operatorname{Lip}(h)\leqslant\|Dh\|_\infty<\infty$.
Fix $t>1$, and define
\begin{equation}\label{eq:target-c0-embedding}
 \Phi(x)=\left(t^m\varphi_j(t^{-m}x)h(t^{-m}x)\right)_{(m,j)\in\mathbb Z\times\mathbb N}.
\end{equation}
Fix $x\in X$ and $\delta>0$. Since $0\leqslant\varphi_j,h\leqslant1$,
we have $|\Phi_{m,j}(x)|\leqslant t^m$, so all these coordinates have absolute value less than $\delta$ when $m$ is sufficiently negative.
When $m$ is sufficiently large, $t^{-m}x$ lies in a neighborhood on which $h$ vanishes, so, for every $j$,
$\Phi_{m,j}(x)=0$. For each of the finitely many remaining $m$, local finiteness of $(\varphi_j)$
implies that only finitely many coordinates are nonzero at $t^{-m}x$. Therefore,
\begin{equation*}
 \#\{(m,j):|\Phi_{m,j}(x)|\geqslant\delta\}<\infty,
\end{equation*}
that is, $\Phi(x)\in c_0(\mathbb Z\times\mathbb N)$.

For every $x,y\in X$ and $(m,j)\in\mathbb Z\times\mathbb N$, the product estimate gives
\begin{equation*}
 \begin{aligned}
 |\Phi_{m,j}(x)-\Phi_{m,j}(y)|
 &\leqslant t^m(M+C_h)\|t^{-m}x-t^{-m}y\|\\
 &=(M+C_h)\|x-y\|.
 \end{aligned}
\end{equation*}
Taking the supremum over all coordinates gives the upper Lipschitz bound.
If $x\ne y$, choose an integer $m$ such that $2t^m\leqslant\|x-y\|<2t^{m+1}$.
Interchanging $x,y$ if necessary, assume that $\|x\|\geqslant\frac{1}{2}\|x-y\|\geqslant t^m$, so
$h(t^{-m}x)=1$.
Step~\ref{step:polynomial-smoothing:1} gives an index $j$ such that $\varphi_j(t^{-m}x)=1$.
If also $\varphi_j(t^{-m}y)\ne0$, then~\eqref{eq:target-sup-partition-bounds}
places both $t^{-m}x$ and $t^{-m}y$ in $B_X(a_j,\frac{1}{16})$, and hence
\begin{equation*}
 2\leqslant t^{-m}\|x-y\|<\frac{1}{8},
\end{equation*}
a contradiction. Thus $\Phi_{m,j}(y)=0$, whereas $\Phi_{m,j}(x)=t^m$. Consequently,
\begin{equation*}
 \|\Phi(x)-\Phi(y)\|_\infty\geqslant t^m>\frac{\|x-y\|}{2t}.
\end{equation*}
We have therefore proved that, for all $x,y\in X$,
\begin{equation}\label{eq:target-c0-bilipschitz}
 \begin{gathered}
 a\|x-y\|\leqslant\|\Phi(x)-\Phi(y)\|_\infty\leqslant b\|x-y\|,\\
 a=\frac{1}{2t},\qquad b=M+C_h.
 \end{gathered}
\end{equation}
By Lemma~\ref{lem:target-class-closure},
$\Phi_{m,j}=t^m(\varphi_j\circ t^{-m}I_X)(h\circ t^{-m}I_X)$ belongs to $\mathscr S(X)$.
Since the index set $\mathbb Z\times\mathbb N$ is countable, the union of the polynomial families used by all coordinates remains countable.

\smallskip\noindent\proofstep{}\label{step:polynomial-smoothing:3}
If $L=0$, take $v=F$. Assume henceforth that $L>0$, and identify the countable coordinate set $\mathbb Z\times\mathbb N$ with $\mathbb N$.
Define $H(\Phi(x))=F(x)$ on $\Phi(X)\subset c_0$. Then
$\operatorname{Lip}(H)\leqslant \frac{L}{a}$.
The McShane extension~\cite[Theorem~1, p.~838]{McShane1934} is defined, for every $y\in c_0$, by
\begin{equation*}
 \widetilde H(y)=\inf\limits_{x\in X}
       \left(F(x)+\frac{L}{a}\|y-\Phi(x)\|_\infty\right).
\end{equation*}
For a fixed $x_0\in X$, the Lipschitz bound on $H$ and the triangle inequality give
\begin{equation*}
 F(x_0)-\frac{L}{a}\|y-\Phi(x_0)\|_\infty
 \leqslant\widetilde H(y)
 \leqslant F(x_0)+\frac{L}{a}\|y-\Phi(x_0)\|_\infty.
\end{equation*}
Thus $\widetilde H$ is finite. Comparing the defining infima at $y,z$ and
then interchanging $y,z$ gives
\begin{equation*}
 |\widetilde H(y)-\widetilde H(z)|\leqslant\frac{L}{a}\|y-z\|_\infty,
 \qquad \widetilde H(\Phi(x))=F(x).
\end{equation*}
By~\cite[Theorem~G, p.~564]{HajekJohanis2010}, there is a mapping
$G\in C^\infty(c_0,\mathbb R)$ locally depending on only finitely many coordinates such that
\begin{equation*}
 \|G-\widetilde H\|_\infty\leqslant\frac{\eta}{2},
 \qquad\operatorname{Lip}(G)\leqslant\frac{L}{a}.
\end{equation*}
Set $v=G\circ\Phi$. Fix $x_0\in X$, and choose an open neighborhood $W\subset c_0$ of $\Phi(x_0)$,
a finite coordinate set $S$, an open set $O\subset\mathbb R^S$, and $B\in C^\infty(O,\mathbb R)$
such that $y_S\in O$ and $G(y)=B(y_S)$ for every $y\in W$.
Since $\Phi$ is continuous, $V=\Phi^{-1}(W)$ is an open neighborhood of $x_0$, and, for every $x\in V$,
\begin{equation*}
 v(x)=B\bigl((\Phi_\gamma(x))_{\gamma\in S}\bigr).
\end{equation*}
For each $\gamma\in S$, write locally
$\Phi_\gamma=b_\gamma\circ P_\gamma$, where the components of $P_\gamma$
belong to the countable family from Step~\ref{step:polynomial-smoothing:2}.
Set
\begin{equation*}
 P=(P_\gamma)_{\gamma\in S},\qquad
 C((t_\gamma)_{\gamma\in S})=(b_\gamma(t_\gamma))_{\gamma\in S},\qquad
 \widetilde B=B\circ C.
\end{equation*}
The domain of $\widetilde B$ is the open subset of
$\prod\limits_{\gamma\in S}\operatorname{dom}b_\gamma$ on which $C$ takes
values in $O$. Restricting $V$ to the intersection of the finitely many
representation neighborhoods gives $v|_V=\widetilde B\circ P|_V$.
All polynomial coordinates belong to the same fixed countable family, so
$v\in\mathscr S(X)$. 
Moreover,
\begin{equation*}
 \|v-F\|_\infty\leqslant\frac{\eta}{2}<\eta,
 \qquad \|Dv\|_\infty\leqslant\operatorname{Lip}(v)\leqslant\frac{b}{a}L.
\end{equation*}
Taking $K_X=\max\limits\{1,\frac{b}{a}\}$ gives~\eqref{eq:target-scalar-lipschitz-result}.
\end{proof}

\subsection[Scalar and finite-dimensional-range approximation]{\texorpdfstring{$C^1$}{C1} approximation for scalar mappings and mappings with locally finite-dimensional range}\label{subsec:target-finite-range-approximation}

\begin{proposition}\label{prop:target-affine-gluing}
Suppose that a real Banach space $Y$ satisfies the following condition: there exists $K\geqslant1$ such that, for every nonempty open ball
$B\subset X$, every $L>0$, every $r\in C^1(B,Y)$, and every $\delta>0$, if
$\sup\limits_B\|Dr\|\leqslant L$, then there exists $v\in\mathscr S(B,Y)$ satisfying
\begin{equation}\label{eq:target-local-approximation-property}
 \sup\limits_B\|v-r\|<\delta,
 \qquad\sup\limits_B\|Dv\|\leqslant KL.
\end{equation}
Then, for every open set $U\subset X$, each $f\in C^1(U,Y)$ admits approximation in the fine $C^1$ topology by mappings of the form
\begin{equation}\label{eq:target-affine-class}
 h(x)=b(x)+A(x)x,
 \qquad b\in\mathscr S(U,Y),\quad
 A\in\mathscr S(U,\mathcal L(X,Y))
\end{equation}
and consequently also admits fine $C^1$ real-analytic approximation.
\end{proposition}
\begin{proof}
Let $\varepsilon\in C(U,(0,\infty))$.
For every $a\in U$, choose an open ball $B_a\subset U$ centered at $a$ such that, for every $x\in B_a$,
\begin{equation}\label{eq:target-local-oscillation}
 \varepsilon(x)>\frac{\varepsilon(a)}{2},
 \qquad\|Df(x)-Df(a)\|<\frac{\varepsilon(a)}{32(K+1)}.
\end{equation}
By Lemma~\ref{lem:target-countable-partition}, take a subordinate partition of unity $(\psi_i)$,
and denote the corresponding balls and centers by $B_i,a_i$.
In Lemma~\ref{lem:affine-difference-gluing}, take $V_i=B_i$ and $\varphi_i=\psi_i$,
retain the local affine mappings $T_i$ defined there, write
$r_i=f|_{B_i}-T_i|_{B_i}$, and set
\begin{equation*}
 \begin{gathered}
 M_i=\|D\psi_i\|_\infty,\qquad
 L_i=\frac{\varepsilon(a_i)}{32(K+1)},\qquad
 \tau_i=\frac{2^{-i-5}\varepsilon(a_i)}{1+M_i}.
 \end{gathered}
\end{equation*}
By~\eqref{eq:target-local-oscillation}, $\|Dr_i\|_{\infty,B_i}\leqslant L_i$.
By~\eqref{eq:target-local-approximation-property}, choose $v_i\in\mathscr S(B_i,Y)$ such that
$\|v_i-r_i\|_{\infty,B_i}<\tau_i$, $\|Dv_i\|_{\infty,B_i}\leqslant KL_i$.
Define $A,b,h$ by~\eqref{eq:affine-difference-gluing-definition}, and set $e_i=v_i-r_i$.
Constant mappings belong to $\mathscr S$, and $\operatorname{supp}_U\psi_i\subset B_i$.
The extension-by-zero and countable locally finite summation properties in Lemma~\ref{lem:target-class-closure}
show that $b,A$ belong to the respective classes $\mathscr S$ in~\eqref{eq:target-affine-class}.
For $x\in U$, set $I_x=\{i:\psi_i(x)>0\}$.
If $i\in I_x$, then~\eqref{eq:target-local-oscillation} gives
$\varepsilon(a_i)<2\varepsilon(x)$, $\|e_i(x)\|<\tau_i$, and
$\|De_i(x)\|\leqslant(K+1)L_i$.
Consequently, \eqref{eq:partition-difference-derivative} gives
\begin{equation*}
 \begin{gathered}
 \|h(x)-f(x)\|
 \leqslant\sum\limits_{i\in I_x}\psi_i(x)\|e_i(x)\|
 <\frac{\varepsilon(x)}{16},\\
 \sum\limits_{i\in I_x}\psi_i(x)\|De_i(x)\|<\frac{\varepsilon(x)}{16},\\
 \sum\limits_{i\in I_x}\|D\psi_i(x)\|\,\|e_i(x)\|
 <\varepsilon(x)\sum\limits_{i\geqslant1}2^{-i-4}
 =\frac{\varepsilon(x)}{16}.
 \end{gathered}
\end{equation*}
Thus $\max\limits_{j=0,1}\|D^jh-D^jf\|<\frac{\varepsilon}{8}$.
Apply Proposition~\ref{prop:target-polynomial-analytic} to $(b,A)$,
and then use Lemma~\ref{lem:affine-analytic-approximation} to choose
$g\in C^\omega(U,Y)$ such that, for every $x\in U$,
\begin{equation*}
 \max\limits_{j=0,1}\|D^jg(x)-D^jh(x)\|
 <\frac{\varepsilon(x)}{2}.
\end{equation*}
This proves the assertion.
\end{proof}

\begin{corollary}\label{cor:target-scalar-fine}
For every open set $U\subset X$ and every finite-dimensional real Banach space $Z$,
the class $\mathscr S(U,Z)$ is dense in $C^1(U,Z)$ for the fine $C^1$ topology.
\end{corollary}
\begin{proof}
Let $f\in C^1(U,Z)$ and $\varepsilon\in C(U,(0,\infty))$.
First take $Z=\mathbb R$, and verify the hypothesis of Proposition~\ref{prop:target-affine-gluing}.
If $r\in C^1(B,\mathbb R)$ and $\sup\limits_B\|Dr\|\leqslant L$,
then for $x,y\in B$,
\begin{equation*}
 r(y)-r(x)=\int\limits_0^1Dr(x+t(y-x))[y-x] \,\mathrm{d}t,\qquad
 |r(y)-r(x)|\leqslant L\|y-x\|.
\end{equation*}
The McShane extension~\cite[Theorem~1, p.~838]{McShane1934} is given, for every $x\in X$, by
\begin{equation*}
 R(x)=\inf\limits_{z\in B}\bigl(r(z)+L\|x-z\|\bigr).
\end{equation*}
It satisfies $R|_B=r$ and $\operatorname{Lip}(R)\leqslant L$.
Lemma~\ref{lem:target-scalar-lipschitz} therefore establishes
\eqref{eq:target-local-approximation-property} with $K=K_X$.

In the scalar case, each $Df(a_i)$ is a continuous real linear polynomial, and
\begin{equation*}
 A(x)x=\sum\limits_i\psi_i(x)Df(a_i)[x].
\end{equation*}
Adjoin the countably many linear polynomials $x\mapsto Df(a_i)[x]$ to the polynomial family used.
Lemma~\ref{lem:target-class-closure} gives $x\mapsto A(x)x\in\mathscr S(U)$,
and hence $h\in\mathscr S(U)$.
Thus $\mathscr S(U)$ is dense in $C^1(U,\mathbb R)$ for the fine $C^1$ topology.

If $Z=\{0\}$, take the zero mapping. Otherwise, let $\dim Z=N\geqslant1$,
choose a basis $z_1,\ldots,z_N$, and write
$f=\sum\limits_{j=1}^Nf_jz_j$.
For each $f_j$, choose $g_j\in\mathscr S(U)$ such that
\begin{equation*}
 \max\limits_{k=0,1}\|D^kg_j(x)-D^kf_j(x)\|
 <\frac{\varepsilon(x)}{1+\sum\limits_{j=1}^N\|z_j\|}.
\end{equation*}
Set $g=\sum\limits_{j=1}^Ng_jz_j\in\mathscr S(U,Z)$. Then
\begin{equation*}
\max\limits_{k=0,1}\|D^kg(x)-D^kf(x)\|
 <\frac{\sum\limits_{j=1}^N\|z_j\|}{1+\sum\limits_{j=1}^N\|z_j\|}
 \varepsilon(x)<\varepsilon(x).\qedhere
\end{equation*}
\end{proof}

\begin{lemma}\label{lem:target-locally-finite-range}
Let $U\subset X$ be open, let $Y$ be a real Banach space, and let $h\in C^1(U,Y)$. Suppose that every point has an open neighborhood on which the range of $h$ lies in a finite-dimensional subspace.
Then $h$ admits fine $C^1$ approximation by mappings in $\mathscr S(U,Y)$.
\end{lemma}
\begin{proof}
Choose an open cover $(V_\alpha)$ of $U$ such that, for every
$\alpha$, the set $h(V_\alpha)$ is contained in a finite-dimensional
linear subspace of $Y$.
By Lemma~\ref{lem:target-countable-partition}, choose a countable
locally finite partition of unity $(\psi_i)$ subordinate to this cover.
Its supports lie in open sets $V_i$ with $h(V_i)\subset Y_i$, where each $Y_i\subset Y$ is finite-dimensional.
For $x\in V_i$ and $v\in X$, closedness of $Y_i$ gives
\begin{equation*}
 Dh(x)[v]=\lim\limits_{t\to0}\frac{h(x+tv)-h(x)}{t}\in Y_i.
\end{equation*}
Thus $h|_{V_i}\in C^1(V_i,Y_i)$ with the inherited norm.
Given $\eta\in C(U,(0,\infty))$, use Corollary~\ref{cor:target-scalar-fine}
to choose $v_i\in\mathscr S(V_i,Y_i)$ such that, for every $x\in V_i$,
\begin{equation*}
 \max\limits_{j=0,1}\|D^jv_i(x)-D^jh(x)\|
 <\frac{2^{-i-3}\eta(x)}{1+\|D\psi_i(x)\|}.
\end{equation*}
Inspired by~\cite[proof of Theorem~14, pp.~575--576]{HajekJohanis2010},
set $v=\sum\limits_i\psi_iv_i$, extending each product by zero outside $V_i$.
Since $\operatorname{supp}_U\psi_i\subset V_i$, Lemma~\ref{lem:target-class-closure} gives
$v\in\mathscr S(U,Y)$.
Write $e_i=v_i-h|_{V_i}$. For $I_x=\{i:\psi_i(x)>0\}$,
\eqref{eq:partition-difference-derivative} gives
\begin{equation*}
 \begin{aligned}
 \|v(x)-h(x)\|&\leqslant\sum\limits_{i\in I_x}\psi_i(x)\|e_i(x)\|,\\
 \|Dv(x)-Dh(x)\|
 &\leqslant\sum\limits_{i\in I_x}
       \bigl(\psi_i(x)\|De_i(x)\|+\|D\psi_i(x)\|\|e_i(x)\|\bigr)\\
 &\leqslant\eta(x)\sum\limits_{i\in I_x}2^{-i-3}
             \frac{\psi_i(x)+\|D\psi_i(x)\|}{1+\|D\psi_i(x)\|}
 \leqslant\frac{\eta(x)}{8}.
 \end{aligned}
\end{equation*}
The value estimate is bounded by the same sum, so
$\max\limits_{j=0,1}\|D^jv(x)-D^jh(x)\|\leqslant\frac{\eta(x)}{8}<\eta(x)$.
\end{proof}

\subsection[Approximation for targets with decompositions]{Fine \texorpdfstring{$C^1$}{C1} approximation for targets with decompositions}\label{subsec:target-first-order-result}

\begin{lemma}\label{lem:target-tail-cutoff}
Suppose that $Y$ admits a countable unconditional finite-dimensional decomposition and an equivalent norm that is $C^1$ smooth away from the origin.
Then, for a suitable equivalent norm, there exists $C_Y\geqslant1$ such that, for every $r>0$, there is
$\Theta_r\in C^1(Y,Y)$ satisfying
\begin{equation}\label{eq:target-tail-cutoff-bounds}
 \sup\limits_{y\in Y}\|\Theta_r(y)-y\|\leqslant r,
 \qquad
 \sup\limits_{y\in Y}\|D\Theta_r(y)\|\leqslant C_Y.
\end{equation}
Near every point, each $\Theta_r$ takes values in the sum of finitely many decomposition blocks.
\end{lemma}
\begin{proof}
	If $Y=\{0\}$, take $\Theta_r=0$.
	Assume henceforth that $Y\ne\{0\}$.
	Denote the equivalent norm from Lemma~\ref{lem:unconditional-renorming} by
	$\|\cdot\|$, and use it throughout this proof.
	Let $(Y_n)_{n\geqslant1}$ be an unconditional finite-dimensional decomposition of $Y$,
	and let $R_n:Y\to Y_n$ be the corresponding block projections.
	If the decomposition is finite, append zero blocks after the last term.

	Define the tail projections
	$T_n=I_Y-\sum\limits_{j<n}R_j
	 =\sum\limits_{j\geqslant n}R_j$.
	Equation~\eqref{eq:unconditional-decomposition-multiplier} gives $\|T_n\|\leqslant1$,
	and~\eqref{eq:unconditional-decomposition-definition} gives $T_ny\to0$.

	Let $\nu$ be the equivalent norm in the hypothesis, which is $C^1$ smooth away from the origin.
	By equivalence of the norms, there exist constants $0<\alpha\leqslant\beta$ such that, for every $y\in Y$,
	\begin{equation*}
	\alpha\|y\|\leqslant\nu(y)
	\leqslant\beta\|y\|.
	\end{equation*}
	Set $\mu=\frac{1}{\alpha}\nu$ and $B=\frac{\beta}{\alpha}\geqslant1$.
	Then $\mu$ is still an equivalent norm that is $C^1$ smooth away from the origin, and
	\begin{equation*}
	\|y\|\leqslant\mu(y)\leqslant B\|y\|.
	\end{equation*}
	By the triangle inequality for the norm, for $y\ne0$ and $h\in Y$,
	\begin{equation*}
	|D\mu(y)[h]|
	=\lim\limits_{t\to0}
	\frac{|\mu(y+th)-\mu(y)|}{|t|}
	\leqslant\mu(h)\leqslant B\|h\|.
	\end{equation*}
	Hence $\|D\mu(y)\|\leqslant B$.
Inspired by~\cite[Lemma~20, pp.~580--581]{HajekJohanis2010},
choose $\zeta\in C^\infty(\mathbb R,[0,1])$ such that
$\zeta=0$ on $(-\infty,\frac{1}{2}]$ and $\zeta=1$ on $[1,\infty)$, and define
\begin{equation}\label{eq:target-tail-cutoff-definition}
 \Theta_r(y)=\sum\limits_{n\geqslant1}\zeta\left(\frac{\mu(T_ny)}{r}\right)R_ny.
\end{equation}
Fix $y_0\in Y$. Since $T_ny_0\to0$, choose $N$ such that $\|T_Ny_0\|<\frac{r}{4B}$.
For every $y\in B_Y(y_0,\frac{r}{4B})$ and $n\geqslant N$, the identity $T_n=T_nT_N$ gives
\begin{equation*}
 \mu(T_ny)\leqslant B\|T_ny\|
 \leqslant B\bigl(\|T_Ny_0\|+\|y-y_0\|\bigr)<\frac{r}{2}.
\end{equation*}
Thus~\eqref{eq:target-tail-cutoff-definition} contains only the finitely many terms with $n<N$ on this ball.
Each coefficient is a $C^1$ composition where $T_ny\ne0$. Where $T_ny=0$, it vanishes throughout
$\{\|T_nz\|<\frac{r}{2B}\}$. Therefore $\Theta_r\in C^1(Y,Y)$ and
\begin{equation*}
 \Theta_r\bigl(B_Y(y_0,\frac{r}{4B})\bigr)
 \subset\bigoplus\limits_{n=1}^{N-1}R_nY.
\end{equation*}
Fix $y\in Y$, and write $\lambda_n=\zeta(\frac{\mu(T_ny)}{r})$.
Since $T_ny\to0$, there is a least index $n_0$ such that $\lambda_{n_0}\ne1$, and
$\mu(T_{n_0}y)<r$. For $n\geqslant n_0$, we have $R_nT_{n_0}y=R_ny$, so
\begin{equation}\label{eq:target-tail-cutoff-value-multiplier}
 y-\Theta_r(y)
 =\sum\limits_{n\geqslant n_0}(1-\lambda_n)R_n(T_{n_0}y).
\end{equation}
Since $0\leqslant1-\lambda_n\leqslant1$, \eqref{eq:unconditional-decomposition-multiplier} gives
\begin{equation*}
 \|y-\Theta_r(y)\|
 \leqslant\left(\sup\limits_{n\geqslant n_0}|1-\lambda_n|\right)\|T_{n_0}y\|
 \leqslant\|T_{n_0}y\|
 \leqslant\mu(T_{n_0}y)<r.
\end{equation*}
Termwise differentiation gives
\begin{equation}\label{eq:target-tail-cutoff-derivative}
 \begin{gathered}
 D\Theta_r(y)[h]=\sum\limits_n\zeta\left(\frac{\mu(T_ny)}{r}\right)R_nh
       +\sum\limits_n c_nR_ny,\\
 c_n=\frac{1}{r}\zeta'\left(\frac{\mu(T_ny)}{r}\right)D\mu(T_ny)[T_nh],
 \end{gathered}
\end{equation}
with $c_n=0$ when $T_ny=0$.
By~\eqref{eq:unconditional-decomposition-multiplier}, the first term satisfies
\begin{equation*}
 \left\|\sum\limits_n\zeta\left(\frac{\mu(T_ny)}{r}\right)R_nh\right\|
 \leqslant\left(\sup\limits_n\left|\zeta\left(\frac{\mu(T_ny)}{r}\right)\right|\right)\|h\|
 \leqslant\|h\|.
\end{equation*}
If some $c_n\ne0$, take the least such index $n_1$. Then
$\|T_{n_1}y\|\leqslant\mu(T_{n_1}y)<r$.
Using $\|D\mu(T_ny)\|\leqslant B$ for $T_ny\ne0$, $\|T_n\|\leqslant1$,
and the convention $c_n=0$ when $T_ny=0$, we obtain
\begin{equation*}
 |c_n|\leqslant\frac{\|\zeta'\|_\infty B\|h\|}{r}.
\end{equation*}
Since also $c_n=0$ for $n<n_1$,
\begin{equation}\label{eq:target-tail-cutoff-derivative-multiplier}
 \sum\limits_n c_nR_ny
 =\sum\limits_{n\geqslant n_1}c_nR_n(T_{n_1}y).
\end{equation}
Applying~\eqref{eq:unconditional-decomposition-multiplier} again gives
\begin{equation*}
 \left\|\sum\limits_n c_nR_ny\right\|
 \leqslant\left(\sup\limits_{n\geqslant n_1}|c_n|\right)\|T_{n_1}y\|
 \leqslant\|\zeta'\|_\infty B\|h\|.
\end{equation*}
Thus, for every $y,h\in Y$,
\begin{equation*}
 \|D\Theta_r(y)[h]\|
 \leqslant\|h\|+\|\zeta'\|_\infty B\|h\|
 =(1+\|\zeta'\|_\infty B)\|h\|.
\end{equation*}
We may therefore take $C_Y=1+\|\zeta'\|_\infty B$.
\end{proof}

\begin{lemma}\label{lem:target-local-approximation}
Suppose that $Y$ admits a countable unconditional finite-dimensional
decomposition and an equivalent norm that is $C^1$ smooth away from the
origin. Then the condition in Proposition~\ref{prop:target-affine-gluing} holds.
In the norm of Lemma~\ref{lem:target-tail-cutoff}, one may take $K=C_Y+1$.
\end{lemma}
\begin{proof}
Let $B\subset X$ be a nonempty open ball, let $r\in C^1(B,Y)$, and let $L>0$ and $\delta>0$. Suppose that
$\sup\limits_B\|Dr\|\leqslant L$. Set $u=\Theta_{\frac{\delta}{4}}\circ r$.
The chain rule gives, for every $x\in B$,
\begin{equation*}
 Du(x)=D\Theta_{\frac{\delta}{4}}(r(x))\circ Dr(x).
\end{equation*}
Thus $u\in C^1(B,Y)$ and~\eqref{eq:target-tail-cutoff-bounds} gives
\begin{equation*}
 \|u-r\|_{\infty,B}\leqslant\frac{\delta}{4},
 \qquad \|Du\|_{\infty,B}\leqslant C_YL.
\end{equation*}
By Lemma~\ref{lem:target-tail-cutoff}, choose an open neighborhood $W_a$ of $r(a)$ and a finite-dimensional subspace $Y_a$
such that $\Theta_{\frac{\delta}{4}}(W_a)\subset Y_a$.
Then $r^{-1}(W_a)$ is an open neighborhood of $a$, and $u(r^{-1}(W_a))\subset Y_a$. Hence $u$ has locally finite-dimensional range.
Set $\delta_0=\min\limits\{\frac{\delta}{4},\frac{L}{2}\}$. By Lemma~\ref{lem:target-locally-finite-range},
choose $v\in\mathscr S(B,Y)$ such that, for every $x\in B$,
\begin{equation*}
 \max\limits\{\|v(x)-u(x)\|,\|Dv(x)-Du(x)\|\}<\delta_0.
\end{equation*}
Therefore,
\begin{equation*}
 \begin{aligned}
 \|v-r\|_{\infty,B}
 &\leqslant\delta_0+\frac{\delta}{4}\leqslant\frac{\delta}{2}<\delta,\\
 \|Dv\|_{\infty,B}&\leqslant\delta_0+C_YL
 \leqslant\left(C_Y+\frac{1}{2}\right)L\leqslant KL.
 \end{aligned}
\end{equation*}
This proves \eqref{eq:target-local-approximation-property}.
\end{proof}

\begin{proof}[Proof of Theorem~\ref{thm:main-target}]
If $U=\varnothing$, there is nothing to prove.
Let $(Y_i)_{i\in I}$ be an unconditional finite-dimensional decomposition
of $Y$, and let $R_i$ be the corresponding block projections, which are
bounded by Lemma~\ref{lem:unconditional-renorming}.
For each $y\in Y$, choose finite sets $F_n\subset I$ such that
$\sum\limits_{i\in F_n}R_i y\to y$.
If $i\notin\bigcup\limits_{n\geqslant1}F_n$, applying $R_i$ to this limit
gives $R_i y=0$.
Thus, for every $y\in Y$, the set
$\{i\in I:R_i y\ne0\}$ is at most countable.
Choose a countable dense subset $D\subset U$ and set
\begin{equation*}
 I_0=\bigcup\limits_{a\in D}\{i\in I:R_i f(a)\ne0\},
 \qquad
 Y_0=\overline{\operatorname{span}}\bigcup\limits_{i\in I_0}Y_i.
\end{equation*}
The set $I_0$ is at most countable, and continuity of $f$ gives
$f(U)\subset Y_0$.
If $Y_0=\{0\}$, take $g=f$.
Otherwise, equip $Y_0$ with the norm inherited from $Y$.
The family $(Y_i)_{i\in I_0}$ is an unconditional finite-dimensional
decomposition of $Y_0$, and the restriction to $Y_0$ of the given
equivalent norm on $Y$ is $C^1$ smooth away from the origin.
For every $x\in U$ and $h\in X$, closedness of $Y_0$ gives
\begin{equation*}
 Df(x)[h]=\lim\limits_{t\to0}
 \frac{f(x+th)-f(x)}{t}\in Y_0.
\end{equation*}
The inclusion $\mathcal L(X,Y_0)\to\mathcal L(X,Y)$ is isometric, so the
Fr\'echet remainder estimates and the continuity of $Df$ are unchanged.
Hence $f\in C^1(U,Y_0)$.
Composing an approximation with the inclusion $Y_0\to Y$ preserves
real analyticity and both norms in~\eqref{eq:main-target-approximation}.
It therefore suffices to work with $Y_0$ as the target; replace $Y$ by
$Y_0$ for the remainder of the proof.

Denote the norm chosen in Lemma~\ref{lem:target-tail-cutoff} by $\|\cdot\|$,
denote the original norm by $\|\cdot\|_0$, and choose $a>0$ such that $a\|y\|_0\leqslant\|y\|$.
By Lemma~\ref{lem:target-local-approximation} and Proposition~\ref{prop:target-affine-gluing},
there exists $g\in C^\omega(U,Y)$ such that, for every $x\in U$,
\begin{equation*}
 \max\limits_{j=0,1}\|D^jg(x)-D^jf(x)\|
 <a\varepsilon(x).
\end{equation*}
For $T\in\mathcal L(X,Y)$, the range norm inequality gives
\begin{equation*}
 \|T\|_0=\sup\limits_{\|v\|_X\leqslant1}\|Tv\|_0
       \leqslant\frac{1}{a}\sup\limits_{\|v\|_X\leqslant1}\|Tv\|
       =\frac{1}{a}\|T\|.
\end{equation*}
Consequently,
\begin{equation*}
 \max\limits_{j=0,1}\|D^jg(x)-D^jf(x)\|_0
 \leqslant\frac{1}{a}\max\limits_{j=0,1}\|D^jg(x)-D^jf(x)\|
 <\varepsilon(x),
\end{equation*}
which is~\eqref{eq:main-target-approximation}.
Equivalent norms do not affect the $C^1$ or real-analytic properties, which completes the proof.
\end{proof}
\section[Approximation on Spaces with Property (K)]{Fine \texorpdfstring{$C^0$}{C0} Real-Analytic Approximation on Spaces with Property \texorpdfstring{$(K)$}{(K)}}\label{sec:lipschitz-upgrade}
Fix an arbitrary real Banach space $E$, and write
\begin{equation*}
 Z=\mathbb R\oplus E,\qquad
 \|(t,z)\|_Z=\max\limits\{|t|,\|z\|\}.
\end{equation*}
\refstepcounter{analyticcondition}\label{eq:lipschitz-analytic-property}
\noindent\textbf{Condition $(\mathrm L_E)$.}
For every open set $W\subset Z$, real Banach space $Y$,
Lipschitz mapping $F:W\to Y$, and $\eta>0$, there is a mapping $G\in C^\omega(W,Y)$ such that
\begin{equation*}
 \sup\limits_{w\in W}\|G(w)-F(w)\|<\eta .
\end{equation*}

\begin{theorem}
\label{thm:lipschitz-to-fine}
Suppose that $E$ satisfies \eqref{eq:lipschitz-analytic-property}. Then, for every subset
$A\subset E$, real Banach space $Y$,
$f\in C(A,Y)$, and $\varepsilon\in C(A,(0,\infty))$,
there are an open set $U\subset E$ containing $A$ and a mapping $g\in C^\omega(U,Y)$
such that, for every $x\in A$,
\begin{equation}\label{eq:lipschitz-upgrade-conclusion}
 \|g(x)-f(x)\|<\varepsilon(x).
\end{equation}
\end{theorem}

\begin{lemma}\label{lem:analytic-integral}
Let $E,F$ be real Banach spaces and $b\in C^\omega(\mathbb R\oplus E,F)$. Then the mapping $I_b$ defined, for every $t>0$ and $x\in E$, by
\begin{equation*}
 I_b(t,x)=\int\limits_0^t\frac{b(s,sx)}{(1+s)^2}\,\mathrm{d}s
\end{equation*}
is real analytic.
\end{lemma}
\begin{proof}
The change of variables $s=tu$ gives
\begin{equation*}
I_b(t,x)
=
t\int\limits_0^1
\frac{b(tu,tux)}{(1+tu)^2}\,\mathrm{d}u.
\end{equation*}
It therefore suffices to show that the right-hand side is real analytic in $(t,x)$.
Continuous Banach-valued mappings on compact intervals are integrable; see
\cite[Proposition~6.5, pp.~42--43]{Mujica1986}.

Fix $(t_0,x_0)\in(0,\infty)\times E$.
Choose complexifications $E_{\mathbb C},F_{\mathbb C}$ of $E,F$,
and write
\begin{equation*}
Z_{\mathbb C}
=
\mathbb C\oplus E_{\mathbb C}.
\end{equation*}
For $u\in[0,1]$, set
$
p_u=(ut_0,ut_0x_0).
$

By the complex extension theorem for real-analytic mappings
\cite[Chapter~1, Theorem~171, pp.~75--76]{HajekJohanis2014}, for every $u\in[0,1]$
there are a complex open neighborhood $V_u$ of $p_u$ in $Z_{\mathbb C}$
and a holomorphic mapping
$
B_u:V_u\longrightarrow F_{\mathbb C},
$
such that, for every $w\in V_u\cap(\mathbb R\oplus E)$,
\begin{equation*}
B_u(w)=b(w).
\end{equation*}
Consider the mapping
$
P(v,\tau,z)=(v\tau,v\tau z),
$
where
$
v\in[0,1],
(\tau,z)\in\mathbb C\times E_{\mathbb C}.
$
The mapping $P$ is continuous and, for each fixed $v$,
is holomorphic in $(\tau,z)$.

Since
$
P(u,t_0,x_0)=p_u\in V_u,
$
we can choose an open interval $J_u\subset\mathbb R$ containing $u$
and an open neighborhood of $(t_0,x_0)$
in $\mathbb C\times E_{\mathbb C}$,
denoted by $\Omega_u$, such that
$
P(v,\tau,z)\in V_u
$
for all
$
v\in J_u\cap[0,1],
(\tau,z)\in\Omega_u
$.

Shrinking $\Omega_u$ further, we may assume that every point in it satisfies
$
\operatorname{Re}\tau>\frac{t_0}{2}.
$
Thus, for $v\in[0,1]$,
\begin{equation*}
\operatorname{Re}(1+v\tau)
=
1+v\operatorname{Re}\tau>0,
\end{equation*}
so that $1+v\tau$ does not vanish. By compactness of $[0,1]$, there are finitely many points
$u_1,\ldots,u_m\in[0,1]$
such that
$
[0,1]\subset\bigcup\limits_{j=1}^m J_{u_j}.
$
Choose a partition
\begin{equation*}
0=a_0<a_1<\cdots<a_N=1,
\end{equation*}
such that each closed interval $[a_{k-1},a_k]$
is contained in some $J_{u_{j(k)}}$.

Set
$
\Omega=\bigcap\limits_{j=1}^m\Omega_{u_j}.
$
Then $\Omega$ is a complex open neighborhood of $(t_0,x_0)$.
For each $k=1,\ldots,N$, define
\begin{equation*}
\Psi_k(v,\tau,z)
=
\frac{\tau}{(1+v\tau)^2}
B_{u_{j(k)}}(v\tau,v\tau z),
\end{equation*}
where
$
v\in[a_{k-1},a_k],
(\tau,z)\in\Omega.
$

The continuity of $P$, the holomorphy of $B_{u_{j(k)}}$, and the fact that $1+v\tau\ne0$ imply that
$\Psi_k$ is jointly continuous on
$[a_{k-1},a_k]\times\Omega$
and that, for each fixed $v$,
the mapping
$
(\tau,z)\longmapsto\Psi_k(v,\tau,z)
$
is holomorphic on $\Omega$.

We next show that
\begin{equation*}
\Phi_k(\tau,z)
=
\int\limits_{a_{k-1}}^{a_k}
\Psi_k(v,\tau,z)\,\mathrm{d}v
\end{equation*}
is holomorphic on $\Omega$.

For $q,q_0\in\Omega$,
\begin{equation*}
 \|\Phi_k(q)-\Phi_k(q_0)\|
 \leqslant(a_k-a_{k-1})
       \sup\limits_{v\in[a_{k-1},a_k]}
           \|\Psi_k(v,q)-\Psi_k(v,q_0)\|.
\end{equation*}
Joint continuity and a finite cover of the compact set
$[a_{k-1},a_k]\times\{q_0\}$ show that the supremum tends to zero as
$q\to q_0$. Hence $\Phi_k$ is continuous on $\Omega$.

Take any $q=(\tau,z)\in\Omega$
and any direction
$h\in\mathbb C\times E_{\mathbb C}$.
Since $\Omega$ is open,
there is a number $\rho>0$ such that
$q+\lambda h\in\Omega$ for every $\lambda$ with $|\lambda|\leqslant\rho$.
For each fixed
$v\in[a_{k-1},a_k]$,
the function
$
\lambda\longmapsto
\Psi_k(v,q+\lambda h)
$
is holomorphic on an open neighborhood of the closed disc
$\{|\lambda|\leqslant\rho\}$.
By the Cauchy integral formula for Banach-valued holomorphic functions
\cite[Theorem~7.1, pp.~45--46]{Mujica1986},
for $|\zeta|<\rho$ we have
\begin{equation*}
\Psi_k(v,q+\zeta h)
=
\frac{1}{2\pi i}
\int\limits_{|\lambda|=\rho}
\frac{\Psi_k(v,q+\lambda h)}
{\lambda-\zeta}\,\mathrm{d}\lambda.
\end{equation*}

Since $\Psi_k$ is continuous on the compact set
$
[a_{k-1},a_k]
\times
\{q+\lambda h:|\lambda|=\rho\}
$,
the integrand is uniformly bounded.
For each $\ell\in F_{\mathbb C}^*$,
\cite[Proposition~6.4(a), p.~42]{Mujica1986} allows us to move $\ell$ inside the Bochner integral.
The scalar Fubini theorem and the fact that $F_{\mathbb C}^*$ separates points then allow us to interchange the two Bochner integrals,
giving
\begin{equation*}
\begin{aligned}
	\Phi_k(q+\zeta h)
	&=
	\int\limits_{a_{k-1}}^{a_k}
	\Psi_k(v,q+\zeta h)\,\mathrm{d}v\\
	&=
	\frac{1}{2\pi i}
	\int\limits_{|\lambda|=\rho}
	\frac{1}{\lambda-\zeta}
	\left(
	\int\limits_{a_{k-1}}^{a_k}
	\Psi_k(v,q+\lambda h)\,\mathrm{d}v
	\right)\,\mathrm{d}\lambda\\
	&=
	\frac{1}{2\pi i}
	\int\limits_{|\lambda|=\rho}
	\frac{\Phi_k(q+\lambda h)}
	{\lambda-\zeta}\,\mathrm{d}\lambda.
\end{aligned}
\end{equation*}
For $|\zeta|\leqslant r<\rho$ and $|\lambda|=\rho$,
\begin{equation*}
 \left|\frac{1}{\lambda-\zeta}\right|\leqslant\frac{1}{\rho-r},\qquad
 \left|\frac{1}{(\lambda-\zeta)^2}\right|\leqslant\frac{1}{(\rho-r)^2}.
\end{equation*}
Since $\Phi_k$ is bounded on the integration circle, differentiation under
the integral sign gives
\begin{equation*}
 \frac{\mathrm{d}}{\mathrm{d}\zeta}\Phi_k(q+\zeta h)
 =\frac{1}{2\pi i}\int\limits_{|\lambda|=\rho}
       \frac{\Phi_k(q+\lambda h)}{(\lambda-\zeta)^2}\,\mathrm{d}\lambda.
\end{equation*}
Thus $\zeta\mapsto\Phi_k(q+\zeta h)$ is holomorphic for $|\zeta|<\rho$.

Since $q$ and $h$ were arbitrary,
$\Phi_k$ is G\^ateaux holomorphic.
Moreover, $\Phi_k$ is continuous and hence locally bounded.
The criterion for Fr\'echet holomorphy of locally bounded G\^ateaux holomorphic mappings
\cite[Proposition~8.6 and Theorem~8.7, p.~61]{Mujica1986}
therefore gives
$
\Phi_k\in\mathcal H(\Omega,F_{\mathbb C}).
$

Set
\begin{equation*}
\Phi(\tau,z)=\sum\limits_{k=1}^N\Phi_k(\tau,z).
\end{equation*}
A finite sum of holomorphic mappings is holomorphic,
so
$
\Phi\in\mathcal H(\Omega,F_{\mathbb C}).
$ Finally, for real variables
$(t,x)\in\Omega\cap((0,\infty)\times E)$,
each local holomorphic extension agrees with $b$ on the real subspace, and hence
we have
\begin{equation*}
\begin{aligned}
	\Phi(t,x)
	&=
	\sum\limits_{k=1}^N
	\int\limits_{a_{k-1}}^{a_k}
	\frac{t\,b(vt,vtx)}{(1+vt)^2}\,\mathrm{d}v=
	t\int\limits_0^1
	\frac{b(vt,vtx)}{(1+vt)^2}\,\mathrm{d}v=
	I_b(t,x).
\end{aligned}
\end{equation*}

Thus $I_b$ has the holomorphic extension $\Phi$ on a neighborhood of $(t_0,x_0)$,
and is therefore real analytic on the corresponding real open neighborhood by
\cite[Chapter~1, Fact~156, p.~70]{HajekJohanis2014}.
Since $(t_0,x_0)$ was arbitrary,
we obtain
$
I_b\in C^\omega((0,\infty)\times E,F).
$
\end{proof}

\begin{lemma}\label{lem:analytic-majorant}
If $E$ satisfies \eqref{eq:lipschitz-analytic-property}, then for every nonempty open set $U\subset E$ and
$v\in C(U,\mathbb R)$ there is a mapping $a\in C^\omega(U,(0,\infty))$ such that, for every $x\in U$,
\begin{equation}\label{eq:analytic-majorant}
 a(x)>\max\limits\{1,v(x)\}.
\end{equation}
\end{lemma}

\begin{proof}
\firstproofstep{}\label{step:radial-domain:1}
Set $b(x)=\max\limits\{1,v(x)\}$ and let
\begin{equation*}
 \begin{gathered}
 V=\left\{(t,z)\in Z:t>0,\ \frac{z}{t}\in U,
                  \ t>b\!\left(\frac{z}{t}\right)\right\},\\
 d(w)=\operatorname{dist}(w,Z\setminus V).
 \end{gathered}
\end{equation*}
The conditions defining $V$ are strict inequalities between continuous
functions on $\{(t,z):t>0,\ \frac{z}{t}\in U\}$, so $V$ is open.
If $w=(t,z)\in V$ and $\lambda\geqslant1$, then
\begin{equation*}
 \frac{\lambda z}{\lambda t}=\frac{z}{t}\in U,\qquad
 \lambda t\geqslant t>b\left(\frac{z}{t}\right)
                    =b\left(\frac{\lambda z}{\lambda t}\right).
\end{equation*}
Hence $\lambda V\subset V$. For $C=Z\setminus V$, this implies
$C\subset\lambda C$, and therefore
\begin{equation*}
 d(\lambda w)=\inf\limits_{v\in C}\|\lambda w-v\|
 \geqslant\inf\limits_{v\in\lambda C}\|\lambda w-v\|
 =\lambda\inf\limits_{v\in C}\|w-v\|=\lambda d(w).
\end{equation*}
If $t\leqslant0$, then $(t,z)\in C$ and $d(t,z)=0$; if $t>0$, then
$(0,z)\in C$ and $d(t,z)\leqslant\|(t,z)-(0,z)\|_Z=t$.
Also, $|d(w)-d(w')|\leqslant\|w-w'\|_Z$. Thus, for $\lambda\geqslant1$,
\begin{equation}\label{eq:radial-distance}
 \begin{gathered}
 \lambda V\subset V,\qquad
 Z\setminus V\subset\lambda(Z\setminus V),\\
 d(\lambda w)\geqslant\lambda d(w),\qquad
 0\leqslant d(t,z)\leqslant\max\limits\{t,0\},\qquad
 \operatorname{Lip}(d)\leqslant1.
 \end{gathered}
\end{equation}
By \eqref{eq:lipschitz-analytic-property}, choose $h\in C^\omega(Z,\mathbb R)$ such that, for every $w\in Z$,
\begin{equation}\label{eq:distance-analytic-approximation}
 |h(w)-d(w)|<1.
\end{equation}

\proofstep{}\label{step:radial-domain:2}

Inspired by
\cite[proof of the first lemma in Section~6(A), p.~483]{EellsElworthy1970} and
\cite[proof of Lemma~A, p.~179]{Elworthy1972}, we define
\begin{equation}\label{eq:radial-integral}
 H(t,x)=\int\limits_0^t\frac{e^{h(s,sx)}}{(1+s)^2}\,\mathrm{d}s.
\end{equation}
Applying Lemma~\ref{lem:analytic-integral} to $e^h$ and then using the fundamental theorem of calculus, we obtain
\begin{equation}\label{eq:radial-integral-derivative}
 H\in C^\omega((0,\infty)\times E,\mathbb R),\qquad
 \partial_tH(t,x)=\frac{e^{h(t,tx)}}{(1+t)^2}>0.
\end{equation}
Fix $x\in U$. If $0\leqslant s\leqslant b(x)$, then $(s,sx)\notin V$, and hence
\begin{equation}\label{eq:radial-integral-before-threshold}
 H(b(x),x)<e\int\limits_0^{b(x)}\frac{1}{(1+s)^2}\,\mathrm{d}s<e<3.
\end{equation}
On the other hand, set
\begin{equation*}
 T_x=b(x)+1,\qquad c_x=\frac{d(T_x,T_xx)}{T_x}>0.
\end{equation*}
By \eqref{eq:radial-distance} and \eqref{eq:distance-analytic-approximation},
\begin{equation*}
 s\geqslant T_x\quad\Longrightarrow\quad
 d(s,sx)\geqslant c_xs,\qquad h(s,sx)>c_xs-1.
\end{equation*}
Therefore, as $t\to\infty$,
\begin{equation}\label{eq:radial-integral-growth}
 H(t,x)\geqslant\int\limits_{T_x}^t\frac{e^{c_xs-1}}{(1+s)^2}\,\mathrm{d}s
 \longrightarrow\infty.
\end{equation}

\proofstep{}\label{step:radial-domain:3}
By \eqref{eq:radial-integral-derivative}--\eqref{eq:radial-integral-growth},
for each $x\in U$ there is a unique $a(x)>b(x)$ such that
\begin{equation*}
 H(a(x),x)=3.
\end{equation*}
Set $\mathcal G(t,x)=H(t,x)-3$. For $x_0\in U$,
\begin{equation*}
 \mathcal G(a(x_0),x_0)=0,\qquad
 \partial_t\mathcal G(a(x_0),x_0)
 =\frac{e^{h(a(x_0),a(x_0)x_0)}}{(1+a(x_0))^2}>0.
\end{equation*}
The analytic implicit function theorem
\cite[Theorem~15.3(a), p.~151]{Deimling1985} gives an open $V_{x_0}\subset U$
and $a_{x_0}\in C^\omega(V_{x_0},(0,\infty))$ with
$H(a_{x_0}(x),x)=3$. Strict monotonicity of $H(\,\cdot,x)$ gives
$a_{x_0}(x)=a(x)$ for every $x\in V_{x_0}$, and hence
\begin{equation*}
 a_{x_0}|_{V_{x_0}\cap V_{x_1}}=a_{x_1}|_{V_{x_0}\cap V_{x_1}}.
\end{equation*}
Thus $a\in C^\omega(U,(0,\infty))$, $a>b$, and, for every $x\in U$ and $v\in E$,
\begin{equation*}
 Da(x)[v]=-
 \frac{D_xH(a(x),x)[v]}{\partial_tH(a(x),x)}.
\end{equation*}
In particular, $a$ satisfies~\eqref{eq:analytic-majorant}.
\end{proof}

\begin{lemma}
\label{lem:locally-lipschitz-radius}
Let $E,Y$ be arbitrary real Banach spaces, let $U\subset E$ be open, and let $f\in C(U,Y)$ and $\eta>0$.
Then there are a locally Lipschitz mapping $f_0:U\to Y$ and
$r\in C(U,(0,\frac{1}{2}])$ such that, for every $x\in U$ and $y\in B_E(x,r(x))$,
\begin{equation}\label{eq:locally-lipschitz-radius}
 \begin{gathered}
 \|f_0(x)-f(x)\|<\eta,\qquad B_E(x,r(x))\subset U,\\
 \|f_0(y)-f_0(x)\|\leqslant\frac{\|y-x\|}{r(x)}.
 \end{gathered}
\end{equation}
\end{lemma}

\begin{proof}
The assertion is immediate if $U=\varnothing$. Henceforth assume $U\ne\varnothing$.

\firstproofstep{}\label{step:local-lipschitz:1}
Apply Lemma~\ref{lem:arbitrary-subset-neighborhood} with $A=U$ and $\varepsilon\equiv4\eta$,
and restrict the resulting mapping to $U$. This gives a locally Lipschitz mapping $f_0:U\to Y$ satisfying
$\|f_0(x)-f(x)\|<\eta$.

\proofstep{}\label{step:local-lipschitz:2}
Choose an open cover $(O_\alpha)$ of $U$ and constants $L_\alpha\geqslant1$ such that
$\operatorname{Lip}(f_0|_{O_\alpha})\leqslant L_\alpha$. Define
\begin{equation*}
 d_\alpha(x)=\min\limits\left\{\frac{1}{L_\alpha},
                 \operatorname{dist}(x,E\setminus O_\alpha)\right\},
 \qquad
 r(x)=\frac{1}{2}\sup\limits_\alpha d_\alpha(x).
\end{equation*}
Each $d_\alpha$ is $1$-Lipschitz and satisfies $0\leqslant d_\alpha\leqslant1$, so
\begin{equation*}
 |r(x)-r(y)|\leqslant\frac{1}{2}\sup\limits_\alpha
 |d_\alpha(x)-d_\alpha(y)|\leqslant\frac{1}{2}\|x-y\|.
\end{equation*}
The covering property ensures that $r(x)>0$ on $U$. Fix $x\in U$ and choose $\alpha$ such that
$d_\alpha(x)>r(x)$. Then
\begin{equation*}
 0<r(x)\leqslant\frac{1}{2},\qquad
 B_E(x,r(x))\subset O_\alpha\subset U,\qquad
 L_\alpha<\frac{1}{r(x)}.
\end{equation*}
Thus, for $y\in B_E(x,r(x))$,
$\|f_0(y)-f_0(x)\|\leqslant L_\alpha\|y-x\|\leqslant\frac{\|y-x\|}{r(x)}$.
\end{proof}

\begin{lemma}\label{lem:radial-lipschitz}
Let $U\subset E$ be a nonempty open set, let $Y$ be a real Banach space, and suppose that $f_0\in C(U,Y)$ and
$r\in C(U,(0,\frac{1}{2}])$ satisfy the last two conditions in \eqref{eq:locally-lipschitz-radius}. Set
\begin{equation}\label{eq:radial-lipschitz-definition}
 \begin{gathered}
 v(x)=\frac{4(1+\|x\|)(1+\|f_0(x)\|)}{r(x)},\\
 W=\left\{(t,z)\in Z:t>0,\ \frac{z}{t}\in U,
                         \ t>v\!\left(\frac{z}{t}\right)\right\},\qquad
 F(t,z)=f_0\!\left(\frac{z}{t}\right).
 \end{gathered}
\end{equation}
Then $W$ is open and $\operatorname{Lip}(F)\leqslant1$.
\end{lemma}

\begin{proof}
The mapping $(t,z)\mapsto\frac{z}{t}$ is continuous on $\{t>0\}$.
Since $U$ is open and $v$ is continuous, $W$ is open. Take
$w=(t,tx),w'=(s,sy)\in W$; interchanging the two points if necessary, we may assume that
$\|f_0(x)\|\geqslant\|f_0(y)\|$. Write
\begin{equation*}
 \begin{gathered}
 D=\|w-w'\|_Z=\max\limits\{|t-s|,\|tx-sy\|\},\\
 T=\frac{t r(x)}{2(1+\|x\|)}>2(1+\|f_0(x)\|).
 \end{gathered}
\end{equation*}
If $D\geqslant T$, then
\begin{equation*}
 \|F(w)-F(w')\|\leqslant2\|f_0(x)\|<T\leqslant D.
\end{equation*}
If $0<D<T$, then $T\leqslant\frac{t}{4}$, and hence
\begin{equation*}
 \begin{gathered}
 s\geqslant t-D>\frac{t}{2},\qquad
 s(x-y)=(tx-sy)+(s-t)x,\\
 \|x-y\|\leqslant\frac{D(1+\|x\|)}{s}
       <\frac{2D(1+\|x\|)}{t}<r(x).
 \end{gathered}
\end{equation*}
Consequently,
\begin{equation*}
 \|F(w)-F(w')\|\leqslant\frac{\|x-y\|}{r(x)}
 <\frac{2D(1+\|x\|)}{t r(x)}
 <\frac{D}{2(1+\|f_0(x)\|)}\leqslant D.
\end{equation*}
\end{proof}

\begin{proof}[Proof of Theorem~\ref{thm:lipschitz-to-fine}]
If $A=\varnothing$, take $U=E$ and $g=0$. Assume now that $A\ne\varnothing$.
We first prove the assertion when $A=U$ is open.

\firstproofstep{}\label{step:K-fine:1}
Given $\delta>0$, use Lemma~\ref{lem:locally-lipschitz-radius} to choose $f_0,r$ such that, for every $x\in U$,
\begin{equation*}
 \|f_0(x)-f(x)\|<\frac{\delta}{4},
\end{equation*}
and such that the last two conditions in \eqref{eq:locally-lipschitz-radius} hold.
Define $v,W,F$ by \eqref{eq:radial-lipschitz-definition}.
Lemma~\ref{lem:radial-lipschitz} gives $\operatorname{Lip}(F)\leqslant1$.
By Lemma~\ref{lem:analytic-majorant}, choose $a\in C^\omega(U,(0,\infty))$ such that, for every $x\in U$,
\begin{equation*}
 a(x)>v(x).
\end{equation*}
Then
\begin{equation*}
 J(x)=(a(x),a(x)x),\qquad
 J\in C^\omega(U,Z),\quad J(U)\subset W,\quad F\circ J=f_0.
\end{equation*}
By \eqref{eq:lipschitz-analytic-property}, choose $G\in C^\omega(W,Y)$ such that
$\sup\limits_{w\in W}\|G(w)-F(w)\|<\frac{\delta}{4}$. Set $g=G\circ J$. Then
\begin{equation}\label{eq:continuous-uniform-approximation}
 \sup\limits_{x\in U}\|g(x)-f(x)\|
 \leqslant\frac{\delta}{2}<\delta .
\end{equation}

\proofstep{}\label{step:K-fine:2}
Inspired by \cite[Chapter~7, Lemma~34, p.~405]{HajekJohanis2014},
we use Lemma~\ref{lem:analytic-majorant} to choose $a\in C^\omega(U,(0,\infty))$ such that
$a>\frac{1}{\varepsilon}$. Apply Step~\ref{step:K-fine:1} to the continuous mapping $af:U\to Y$ to obtain
$h\in C^\omega(U,Y)$ satisfying
\begin{equation*}
 \sup\limits_{x\in U}\|h(x)-a(x)f(x)\|<1.
\end{equation*}
Then $g=\frac{h}{a}\in C^\omega(U,Y)$ and, for every $x\in U$,
\begin{equation*}
 \|g(x)-f(x)\|<\frac{1}{a(x)}<\varepsilon(x).
\end{equation*}

\proofstep{}\label{step:K-fine:3}
Step~\ref{step:K-fine:2} proves that every continuous Banach-valued mapping on any open subset of $E$
admits fine $C^0$ real-analytic approximation.
For arbitrary $A\subset E$, Lemma~\ref{lem:arbitrary-subset-neighborhood}
gives an open set $U_0\supset A$, $F_0\in C(U_0,Y)$, and
$\delta\in C(U_0,(0,\infty))$ such that, for every $x\in A$,
\begin{equation*}
 \|F_0(x)-f(x)\|<\delta(x)<\frac{\varepsilon(x)}{4}.
\end{equation*}
Apply Step~\ref{step:K-fine:2} to $F_0,\delta$ on $U_0$ to obtain
$g\in C^\omega(U_0,Y)$ with $\|g-F_0\|<\delta$. Then, for every $x\in A$,
\begin{equation*}
 \|g(x)-f(x)\|
 \leqslant\|g(x)-F_0(x)\|+\|F_0(x)-f(x)\|
 <2\delta(x)<\frac{\varepsilon(x)}{2},
\end{equation*}
which gives~\eqref{eq:lipschitz-upgrade-conclusion}.
\end{proof}

\begin{proof}[Proof of Theorem~\ref{thm:main-K}]
For $d,\alpha$ in \eqref{eq:property-K-definition}, choose
$M>\sup\limits_{d(x)<\alpha}\|x\|$ and set $\widetilde d(x)=d(Mx)$. Then
\begin{equation*}
 \varnothing\ne\{\widetilde d<\alpha\}
 =M^{-1}\{d<\alpha\}\subset B_X(0,1),
\end{equation*}
Moreover, for $n\geqslant0$,
\begin{equation*}
 D^n\widetilde d(x)[h_1,\ldots,h_n]
 =M^nD^nd(Mx)[h_1,\ldots,h_n].
\end{equation*}
The Taylor series of $d$ at $Mx$ converges when $\|Mh\|<R_d$, so the Taylor
series of $\widetilde d$ at $x$ has radius at least $\frac{R_d}{M}$.
Thus the definition of property $(K)$ used here is equivalent to
the condition in \cite[p.~83, definition of property~$(K)$]{CepedelloHajek2001}.

Set $E=X$. Property $(K)$ is preserved under finite direct sums
\cite[p.~94, after Proposition~3]{CepedelloHajek2001},
and $t\mapsto t^2$ shows that $\mathbb R$ has property $(K)$.
Hence the separable space $Z=\mathbb R\oplus E$ has property $(K)$.
Given an open set $W\subset Z$, a real Banach space $Y$, a Lipschitz mapping $F:W\to Y$, and $\eta>0$,
the mapping $F$ is uniformly continuous on $W$.
Apply \cite[Theorem~1, p.~83]{CepedelloHajek2001} to $Z,W,F$ and the closed convex set $Y$
to obtain $G\in C^\omega(W,Y)$ satisfying
\begin{equation*}
 \sup\limits_{w\in W}\|G(w)-F(w)\|\leqslant\frac{\eta}{2}<\eta.
\end{equation*}
Therefore $E$ satisfies condition~\eqref{eq:lipschitz-analytic-property}.
Theorem~\ref{thm:lipschitz-to-fine} now yields \eqref{eq:main-K-approximation}.
\end{proof}
\section[Applications]{Applications}\label{sec:applications}

\subsection[Preservation of Lipschitz constants]{Real-analytic approximation preserving Lipschitz constants}\label{sec:lipschitz-whitney}
\begin{theorem}
\label{thm:lipschitz-analytic}
Suppose that $X$ admits a separating polynomial and an unconditional finite-dimensional decomposition, and equip $X$ with the multiplier norm of Lemma~\ref{lem:unconditional-renorming}.
For every real Banach space $Y$, $L$-Lipschitz mapping $f:X\to Y$,
and constant $\epsilon>0$, there is a mapping $g\in C^\omega(X,Y)$ such that
\begin{equation}\label{eq:lipschitz-analytic-result}
 \|g-f\|_\infty<\epsilon,
 \qquad\operatorname{Lip}(g)\leqslant L.
\end{equation}
If $L>0$, we may also require that, for every $x\in X$,
\begin{equation}\label{eq:lipschitz-analytic-strict-derivative}
 \|Dg(x)\|<L.
\end{equation}
In particular, for every positive integer $m$ and arbitrary $\Gamma$, the assertion holds for $X=\ell^{2m}(\Gamma)$
with its usual norm.
\end{theorem}

\begin{lemma}\label{lem:strict-derivative-analytic-transfer}
Let $X,Y$ be real Banach spaces, and suppose that every mapping in $C^1(X,Y)$ admits fine $C^1$ approximation by real-analytic mappings.
If $L>0$, $v\in C^1(X,Y)$, and $\|Dv(x)\|<L$ on $X$,
then for every $\epsilon>0$ there is a mapping $g\in C^\omega(X,Y)$ such that
\begin{equation*}
 \|g-v\|_\infty<\epsilon,\qquad
 \|Dg(x)\|<L,\qquad
 \operatorname{Lip}(g)\leqslant L.
\end{equation*}
\end{lemma}
\begin{proof}
The continuous function $x\mapsto L-\|Dv(x)\|$ is positive. Set
\begin{equation}\label{eq:lipschitz-exact-epsilon-definition}
 \varepsilon_*(x)=\min\limits\left\{\frac{\epsilon}{2},
                           \frac{L-\|Dv(x)\|}{2}\right\},
\end{equation}
and use the hypothesis to choose $g$ such that $\|g-v\|<\varepsilon_*$ and $\|Dg-Dv\|<\varepsilon_*$. Then
\begin{equation*}
 \|g-v\|_\infty\leqslant\frac{\epsilon}{2}<\epsilon,
 \qquad
 \|Dg(x)\|<\|Dv(x)\|+\frac{L-\|Dv(x)\|}{2}<L.
\end{equation*}
For $x,y\in X$, integration along the segment gives
\begin{equation*}
 \begin{aligned}
 g(y)-g(x)&=\int\limits_0^1Dg(x+t(y-x))[y-x] \,\mathrm{d}t,\\
 \|g(y)-g(x)\|
 &\leqslant\int\limits_0^1\|Dg(x+t(y-x))\|\|y-x\|\,\mathrm{d}t
 \leqslant L\|y-x\|.
 \end{aligned}
\end{equation*}
Hence $\operatorname{Lip}(g)\leqslant L$.
\end{proof}

\begin{lemma}
	\label{lem:lipschitz-local-smooth}
Suppose that $X$ admits a separating polynomial and an unconditional finite-dimensional decomposition, and equip it with
the multiplier norm of Lemma~\ref{lem:unconditional-renorming}.
Let $Y$ be an arbitrary real Banach space.
Given a $K$-Lipschitz mapping $F:X\to Y$ and $\eta>0$,
there is a mapping $v\in C^\infty(X,Y)$ such that the following estimates hold for every $x\in X$:
	\begin{equation}\label{eq:lipschitz-smooth-result}
		\|v-F\|_\infty<\eta,\qquad
		\|Dv(x)\|\leqslant
		\operatorname{Lip}\bigl(F,B_X(x,\eta)\bigr)+\eta.
	\end{equation}
\end{lemma}

\begin{proof}
If $K=0$, then $F$ is constant and we may take $v=F$.
Henceforth assume $K>0$.

	\firstproofstep{}\label{step:lipschitz-local:1}
Inspired by \cite[proof of Theorem~14, pp.~575--576]{HajekJohanis2010}, we set
	\begin{equation}\label{eq:lipschitz-delta-theta-definition}
		\delta=\theta=\frac{\eta}{8(1+K)}.
	\end{equation}
Let $(E_i)_{i\in I}$ be the decomposition used to define the
multiplier norm, and let $R_i$ be its associated coordinate projections.
By Proposition~\ref{prop:ufdd-strengthened-approximation},
there is a mapping $\Phi\in C^\infty(X,X)$ satisfying
	\begin{equation}\label{eq:lipschitz-Phi-estimates}
		\|\Phi-I_X\|_\infty<\delta,\qquad
		\|D\Phi\|_\infty\leqslant1+\theta,
	\end{equation}
and there are an open cover $(V_\alpha)$ and subspaces
$M_\alpha$, each a sum of finitely many original blocks, such that
	$
	\Phi(V_\alpha)\subset M_\alpha.
	$
Choose a finite set $F_\alpha\subset I$ such that
$M_\alpha=\bigoplus\limits_{i\in F_\alpha}E_i$, and set
$P_\alpha=\sum\limits_{i\in F_\alpha}R_i:X\to M_\alpha$.
By \eqref{eq:unconditional-decomposition-multiplier},
	$\|P_\alpha\|\leqslant1$.

Apply Corollary~\ref{cor:ufdd-quantitative-partition}
to the open cover $(V_\alpha)$ to obtain a nonnegative, locally finite
$C^\infty$ partition of unity
	$(\varphi_{n,j})_{n\geqslant1,\,j\in A_n}$,
together with indices $\alpha(n,j)$ and constants $D_n<\infty$, such that
	\begin{equation}\label{eq:lipschitz-partition-bounds}
		\operatorname{supp}\varphi_{n,j}
		\subset V_{\alpha(n,j)},\qquad
		\sum\limits_{j\in A_n}\|D\varphi_{n,j}(x)\|
		\leqslant D_n
	\end{equation}
for every $n\geqslant1$ and $x\in X$.

	\proofstep{}\label{step:lipschitz-local:2}
For each $n\geqslant1$, choose
	\begin{equation}\label{eq:lipschitz-an-definition}
		a_n=\min\limits\left\{
		\frac{\eta}{8},
		\frac{\eta2^{-n}}{8(1+D_n)}
		\right\}.
	\end{equation}
Then
	\begin{equation*}
	\sup\limits_na_n\leqslant\frac{\eta}{8},
	\qquad
	\sum\limits_{n\geqslant1}D_na_n
	\leqslant\frac{\eta}{8}
	\sum\limits_{n\geqslant1}2^{-n}
	=\frac{\eta}{8}.
	\end{equation*}

Fix $(n,j)$ and write
	$M=M_{\alpha(n,j)}$ and $P=P_{\alpha(n,j)}$.
First suppose that $M\ne\{0\}$.
Fix a Lebesgue measure on the finite-dimensional space $M$,
and choose a nonnegative function $\kappa_{n,j}\in C_c^\infty(M)$
such that
	\begin{equation}\label{eq:lipschitz-local-kernel-support}
		\int\limits_M\kappa_{n,j}(w)\,\mathrm{d}w=1,
		\qquad
		\operatorname{supp}\kappa_{n,j}
		\subset B_M\left(
		0,\min\limits\left\{
		\frac{\eta}{8},
		\frac{a_n}{2(1+K)}
		\right\}\right).
	\end{equation}
For every $z\in M$, define
	\begin{equation}\label{eq:lipschitz-local-convolution-definition}
		h_{n,j}(z)
		=\int\limits_M F(z-w)\kappa_{n,j}(w)\,\mathrm{d}w.
	\end{equation}
Since $F|_M$ is continuous and $\kappa_{n,j}$ has compact support,
the integrand is continuous and compactly supported. Thus the Bochner integral in \eqref{eq:lipschitz-local-convolution-definition}
exists.

A change of variables gives
	\begin{equation*}
	h_{n,j}(z)
	=\int\limits_M F(u)\kappa_{n,j}(z-u)\,\mathrm{d}u.
	\end{equation*}
Fix $z_0\in M$ and choose a bounded open ball $W$ centered at $z_0$.
For $z\in W$, the supports of the integrands are all contained in
the same compact set
$\overline W-\operatorname{supp}\kappa_{n,j}$.
The mapping $F$ is bounded on this compact set,
and all derivatives of the kernel are continuous and bounded.
The dominated convergence theorem for Bochner integrals
therefore permits repeated differentiation under the integral sign, giving, for every $z\in W$, every integer $k\geqslant1$, and every $v_1,\ldots,v_k\in M$,
	\begin{equation*}
	\begin{aligned}
		&D^kh_{n,j}(z)[v_1,\ldots,v_k]=
		\int\limits_M F(u)
		D^k\kappa_{n,j}(z-u)[v_1,\ldots,v_k]\,\mathrm{d}u.
	\end{aligned}
	\end{equation*}
Another application of dominated convergence shows that all these derivatives are continuous, and hence
	$h_{n,j}\in C^\infty(M,Y)$.

By \eqref{eq:lipschitz-local-kernel-support},
	\begin{equation*}
	\begin{aligned}
		\|h_{n,j}(z)-F(z)\|
		&\leqslant
		\int\limits_M
		\|F(z-w)-F(z)\|
		\kappa_{n,j}(w)\,\mathrm{d}w\\
		&\leqslant
		K\int\limits_M\|w\|\kappa_{n,j}(w)\,\mathrm{d}w\leqslant\frac{Ka_n}{2(1+K)}<a_n.
	\end{aligned}
	\end{equation*}
Consequently,
	\begin{equation}\label{eq:lipschitz-finite-convolution}
		\|h_{n,j}-F|_M\|_\infty<a_n.
	\end{equation}

If $M=\{0\}$, define
	$h_{n,j}(0)=F(0)$.
In this case $h_{n,j}$ is smooth,
and \eqref{eq:lipschitz-finite-convolution} remains valid.

	\proofstep{}\label{step:lipschitz-local:3}
Fix
$x\in\operatorname{supp}\varphi_{n,j}$ and set
	$
	k_x=\operatorname{Lip}\bigl(F,B_X(x,\eta)\bigr).
	$
By Step~\ref{step:lipschitz-local:1},
$x\in V_{\alpha(n,j)}$, so
	$\Phi(x)\in M$.

Suppose that $M\ne\{0\}$.
For any
	$z,z'\in B_M(\Phi(x),\frac{\eta}{4})$
and $w\in\operatorname{supp}\kappa_{n,j}$,
\eqref{eq:lipschitz-Phi-estimates} and
\eqref{eq:lipschitz-local-kernel-support} give
	\begin{equation*}
	\begin{aligned}
		\|z-w-x\|
		&\leqslant
		\|z-\Phi(x)\|
		+\|\Phi(x)-x\|+\|w\|\\
		&<\frac{\eta}{4}
		+\frac{\eta}{8}
		+\frac{\eta}{8}<\eta.
	\end{aligned}
	\end{equation*}
Similarly, $\|z'-w-x\|<\eta$.
Thus $z-w,z'-w\in B_X(x,\eta)$, and hence
	\begin{equation*}
	\|F(z-w)-F(z'-w)\|
	\leqslant k_x\|z-z'\|.
	\end{equation*}
Substituting this inequality into \eqref{eq:lipschitz-local-convolution-definition}, we obtain
	\begin{equation*}
	\|h_{n,j}(z)-h_{n,j}(z')\|
	\leqslant k_x\|z-z'\|.
	\end{equation*}
That is,
	\begin{equation}\label{eq:lipschitz-local-convolution-lipschitz}
	\operatorname{Lip}\left(
	h_{n,j},
	B_M\left(\Phi(x),\frac{\eta}{4}\right)
	\right)\leqslant k_x.
	\end{equation}
Since $h_{n,j}$ is differentiable, for every $e\in M$
we may apply \eqref{eq:lipschitz-local-convolution-lipschitz} to
$\Phi(x)+te$ and $\Phi(x)$ and then let $t\to0$ to obtain
	\begin{equation*}
	\|Dh_{n,j}(\Phi(x))e\|
	\leqslant k_x\|e\|.
	\end{equation*}
Hence
	\begin{equation}\label{eq:lipschitz-local-convolution-derivative}
		\|Dh_{n,j}(\Phi(x))\|\leqslant k_x.
	\end{equation}
If $M=\{0\}$, then
$Dh_{n,j}=0$, so \eqref{eq:lipschitz-local-convolution-derivative} holds.

	\proofstep{}\label{step:lipschitz-local:4}
Define
	\begin{equation}\label{eq:lipschitz-local-gluing-definition}
		v_{n,j}
		=h_{n,j}\circ P_{\alpha(n,j)}\circ\Phi,
		\qquad
		v=\sum\limits_{n,j}\varphi_{n,j}v_{n,j}.
	\end{equation}
By the chain rule,
each $v_{n,j}$ belongs to $C^\infty(X,Y)$.
The partition of unity is locally finite, so
Lemma~\ref{lem:partition-summation} gives
	$v\in C^\infty(X,Y)$.

Fix $x\in\operatorname{supp}\varphi_{n,j}$,
and continue to write $M=M_{\alpha(n,j)}$ and
	$P=P_{\alpha(n,j)}$.
Since $\Phi(V_{\alpha(n,j)})\subset M$
and $P|_M=I_M$, we have, for every $y\in V_{\alpha(n,j)}$,
	\begin{equation*}
	P\Phi(y)=\Phi(y).
	\end{equation*}
Differentiating on this open set gives
	\begin{equation*}
	PD\Phi(x)=D\Phi(x).
	\end{equation*}
Thus
	\begin{equation*}
	v_{n,j}(x)=h_{n,j}(\Phi(x)),
	\end{equation*}
and the chain rule yields
	\begin{equation*}
	\begin{aligned}
		Dv_{n,j}(x)
		&=Dh_{n,j}(\Phi(x))\circ PD\Phi(x)\\
		&=Dh_{n,j}(\Phi(x))\circ D\Phi(x).
	\end{aligned}
	\end{equation*}
Combining \eqref{eq:lipschitz-finite-convolution},
\eqref{eq:lipschitz-local-convolution-derivative}, and
\eqref{eq:lipschitz-Phi-estimates}, we obtain
	\begin{equation}\label{eq:lipschitz-local-component-bounds}
		\begin{gathered}
			\|v_{n,j}(x)-F(\Phi(x))\|<a_n,\\
			\|Dv_{n,j}(x)\|
			\leqslant(1+\theta)k_x.
		\end{gathered}
	\end{equation}

Since $\varphi_{n,j}\geqslant0$ and
	$\sum\limits_{n,j}\varphi_{n,j}=1$,
for every $x\in X$ we have
	\begin{equation*}
	\begin{aligned}
		\|v(x)-F(\Phi(x))\|
		&\leqslant
		\sum\limits_{n,j}\varphi_{n,j}(x)
		\|v_{n,j}(x)-F(\Phi(x))\|\\
		&\leqslant
		\sum\limits_{n,j}\varphi_{n,j}(x)a_n\leqslant\sup\limits_na_n.
	\end{aligned}
	\end{equation*}
Since $F$ is $K$-Lipschitz, \eqref{eq:lipschitz-Phi-estimates} also gives

	\begin{equation*}
	\|F(\Phi(x))-F(x)\|
	\leqslant K\|\Phi(x)-x\|
	\leqslant K\delta.
	\end{equation*}
Therefore
	\begin{equation}\label{eq:lipschitz-local-value-bound}
		\begin{aligned}
			\|v-F\|_\infty
			&\leqslant\sup\limits_na_n+K\delta\leqslant\frac{\eta}{8}
			+\frac{K\eta}{8(1+K)}
			<\frac{\eta}{4}<\eta.
		\end{aligned}
	\end{equation}

Since $\sum\limits_{n,j}D\varphi_{n,j}(x)=0$, the product rule gives,
for $h\in X$,
\begin{equation*}
 \begin{aligned}
 Dv(x)[h]
 &=\sum\limits_{n,j}\varphi_{n,j}(x)Dv_{n,j}(x)[h]\\
 &\quad+\sum\limits_{n,j}D\varphi_{n,j}(x)[h]
                  (v_{n,j}(x)-F(\Phi(x))).
 \end{aligned}
\end{equation*}
 Using
\eqref{eq:lipschitz-local-component-bounds} and
\eqref{eq:lipschitz-partition-bounds}, we obtain
	\begin{equation*}
	\begin{aligned}
		\|Dv(x)\|
		&\leqslant(1+\theta)k_x
		+\sum\limits_{n\geqslant1}a_n
		\sum\limits_{j\in A_n}\|D\varphi_{n,j}(x)\|\\
		&\leqslant(1+\theta)k_x+\sum\limits_{n\geqslant1}D_na_n.
	\end{aligned}
	\end{equation*}
Note that $0\leqslant k_x\leqslant K$.
By \eqref{eq:lipschitz-delta-theta-definition},
\eqref{eq:lipschitz-an-definition}, and
\eqref{eq:lipschitz-partition-bounds},
	\begin{equation}\label{eq:lipschitz-local-derivative-bound}
	\begin{aligned}
		\|Dv(x)\|&\leqslant k_x+\theta K+\sum\limits_{n\geqslant1}D_na_n\\
		&\leqslant k_x+\frac{K\eta}{8(1+K)}+\frac{\eta}{8}
		<k_x+\frac{\eta}{4}<k_x+\eta.
	\end{aligned}
	\end{equation}
Together with \eqref{eq:lipschitz-local-value-bound},
this proves \eqref{eq:lipschitz-smooth-result}.

\end{proof}

\begin{lemma}
\label{lem:radial-strict-contraction}
Let $E$ be an arbitrary real normed space, and let $0<a<1$. For every $x\in E$, define
\begin{equation}\label{eq:radial-strict-map}
 R_a(x)=\left(1-\frac{a}{\sqrt{1+\|x\|^2}}\right)x.
\end{equation}
Then
\begin{equation}\label{eq:radial-strict-bounds}
 \|R_a-I_E\|_\infty\leqslant a,\qquad
 \operatorname{Lip}(R_a)\leqslant1,
\end{equation}
and, for every $R>0$,
\begin{equation}\label{eq:radial-strict-ball-bound}
 \operatorname{Lip}\bigl(R_a,\overline B_E(0,R)\bigr)
 \leqslant1-\frac{a}{(1+R^2)^{\frac{3}{2}}}<1.
\end{equation}
\end{lemma}
\begin{proof}
For every $r\geqslant0$, set
\begin{equation*}
 s(r)=1-\frac{a}{\sqrt{1+r^2}},\qquad
 \rho(r)=rs(r).
\end{equation*}
Then $s$ is nonnegative and increasing, and
\begin{equation}\label{eq:radial-strict-scalar-derivative}
 \rho'(r)=1-\frac{a}{(1+r^2)^{\frac{3}{2}}},\qquad
 0<\rho'(r)<1.
\end{equation}
For $x,y\in E$, assume that $r=\|x\|\geqslant t=\|y\|$.
If $r>t$, the reverse triangle inequality $r-t\leqslant\|x-y\|$ gives
\begin{equation}\label{eq:radial-secant-bound}
 \begin{aligned}
 \|R_a(x)-R_a(y)\|
 &\leqslant s(r)\|x-y\|+(s(r)-s(t))t\\
 &\leqslant\left(s(r)+\frac{t(s(r)-s(t))}{r-t}\right)\|x-y\|\\
 &=\frac{\rho(r)-\rho(t)}{r-t}\|x-y\|.
 \end{aligned}
\end{equation}
If $r=t$, then
$\|R_a(x)-R_a(y)\|=s(r)\|x-y\|$, and $s(r)\leqslant\rho'(r)$.
By \eqref{eq:radial-strict-scalar-derivative} and the mean value theorem, for all $t,r$ with $0\leqslant t<r\leqslant R$,
\begin{equation*}
 \frac{\rho(r)-\rho(t)}{r-t}
 \leqslant\max\limits_{0\leqslant u\leqslant R}\rho'(u)
 =1-\frac{a}{(1+R^2)^{\frac{3}{2}}}.
\end{equation*}
Together with the estimate for $r=t$, this gives \eqref{eq:radial-strict-ball-bound}.
Letting $R\to\infty$ yields $\operatorname{Lip}(R_a)\leqslant1$. Moreover,
\begin{equation*}
 \|R_a(x)-x\|
 =\frac{a\|x\|}{\sqrt{1+\|x\|^2}}\leqslant a,
\end{equation*}
which proves \eqref{eq:radial-strict-bounds}.
\end{proof}

\begin{lemma}
\label{lem:lipschitz-strict-smooth}
Let $X,Y$ and their norms satisfy the assumptions on the spaces and norms in Lemma~\ref{lem:lipschitz-local-smooth}. Let $L>0$ and let $u:X\to Y$ be an $L$-Lipschitz mapping.
If, for every $R>0$,
\begin{equation}\label{eq:lipschitz-strict-local-hypothesis}
 \operatorname{Lip}\bigl(u,\overline B_X(0,R)\bigr)<L,
\end{equation}
then for every $\epsilon>0$ there is a mapping $v\in C^\infty(X,Y)$ such that, for every $x\in X$,
\begin{equation}\label{eq:lipschitz-strict-smooth-result}
 \|v-u\|_\infty<\epsilon,
 \qquad \|Dv(x)\|<L.
\end{equation}
\end{lemma}
\begin{proof}
\firstproofstep{}\label{step:lipschitz-strict:1}
By Lemma~\ref{lem:homogeneous-separator}, choose an even integer $d\geqslant2$, a homogeneous polynomial $Q$,
and a constant $B\geqslant1$ such that $Q(x)\geqslant\|x\|^d$ and
$\mu=Q^{\frac{1}{d}}$ is smooth away from the origin with $\|D\mu\|\leqslant B$. Set
\begin{equation}\label{eq:lipschitz-coercive-gauge}
 \nu=(1+Q)^{\frac{1}{d}}.
\end{equation}
Then $\nu\in C^\infty(X,\mathbb R)$, $\nu(x)\geqslant\max\limits\{1,\|x\|\}$,
$D\nu(0)=0$, and, for every $x\ne0$,
\begin{equation*}
 D\nu(x)=\left(\frac{\mu(x)}{(1+\mu(x)^d)^{\frac{1}{d}}}\right)^{d-1}D\mu(x).
\end{equation*}
Hence $\|D\nu\|_\infty\leqslant B$.
Choose a nonincreasing function $\beta\in C^\infty(\mathbb R,[0,1])$ such that
$\beta=1$ on $(-\infty,0]$ and $\beta=0$ on $[1,\infty)$, and set
\begin{equation}\label{eq:lipschitz-annular-partition}
 \begin{gathered}
 \chi_n=\beta(\nu-n),\qquad \psi_1=\chi_1,\\
 \psi_n=\chi_n-\chi_{n-1}\quad\text{for every }n\geqslant2.
 \end{gathered}
\end{equation}
Since $\chi_n$ increases with $n$, we have $\psi_n\geqslant0$.
The function $\nu$ is bounded near every point, so $\chi_n$ is identically one on such a neighborhood for all sufficiently large $n$.
For every $N\geqslant1$,
\begin{equation*}
 \sum\limits_{n=1}^N\psi_n=\chi_N,
\end{equation*}
and, for every $n\geqslant2$, $\psi_n=\chi_n-\chi_{n-1}=0$ wherever $\chi_n=\chi_{n-1}=1$.
Thus $(\psi_n)$ is locally finite and $\sum\limits_n\psi_n=1$.
The chain rule gives $D\psi_1=\beta'(\nu-1)D\nu$ and, for every $n\geqslant2$,
\begin{equation*}
 D\psi_n=(\beta'(\nu-n)-\beta'(\nu-n+1))D\nu.
\end{equation*}
Consequently,
\begin{equation}\label{eq:lipschitz-annular-layer-estimate}
 \operatorname{supp}\psi_n\subset\{\nu\leqslant n+1\}
 \subset\overline B_X(0,n+1),\qquad
 \|D\psi_n\|_\infty\leqslant C:=2B\|\beta'\|_\infty.
\end{equation}
For $n\geqslant2$, $\psi_n=0$ on $\{\nu\leqslant n-1\}\cup\{\nu\geqslant n+1\}$.
Consequently, at any point at most two consecutive functions $\psi_n$ are nonzero.

\proofstep{}\label{step:lipschitz-strict:2}
Set
\begin{equation}\label{eq:lipschitz-annular-gaps}
 m_n=L-\operatorname{Lip}(u,\overline B_X(0,n+3))>0.
\end{equation}
The sequence $(m_n)$ is nonincreasing. Choose
\begin{equation}\label{eq:lipschitz-eta-n-choice}
 0<\eta_n<\min\limits\left\{\frac{\epsilon}{2},1,
                          \frac{m_{n+1}}{8(C+1)}\right\},
\end{equation}
and use Lemma~\ref{lem:lipschitz-local-smooth} to obtain $v_n\in C^\infty(X,Y)$ such that
\begin{equation*}
 \|v_n-u\|_\infty<\eta_n,\qquad
 \|Dv_n(x)\|\leqslant\operatorname{Lip}(u,B_X(x,\eta_n))+\eta_n.
\end{equation*}
If $x\in\operatorname{supp}\psi_n$, then $B_X(x,\eta_n)\subset\overline B_X(0,n+3)$, and hence
\begin{equation}\label{eq:lipschitz-annular-local-derivative}
 \|Dv_n(x)\|\leqslant L-m_n+\eta_n
 \leqslant L-\frac{7}{8}m_n.
\end{equation}

\proofstep{}\label{step:lipschitz-strict:3}
Set $v=\sum\limits_n\psi_nv_n$. By Lemma~\ref{lem:partition-summation},
$v\in C^\infty(X,Y)$ and
\begin{equation*}
 \|v-u\|_\infty\leqslant\sup\limits_n\eta_n
 \leqslant\frac{\epsilon}{2}<\epsilon.
\end{equation*}
Fix $x$. There is an integer $p\geqslant1$
such that, whenever $n\notin\{p,p+1\}$,
$\psi_n(x)=0$ and $D\psi_n(x)=0$.
Take $u$ as the comparison mapping in \eqref{eq:partition-summation-bound}.
By \eqref{eq:lipschitz-annular-local-derivative} and $m_p\geqslant m_{p+1}$,
the two terms on the right-hand side satisfy, respectively,
\begin{equation*}
 \sum\limits_n\psi_n(x)\|Dv_n(x)\|
 \leqslant\sum\limits_n\psi_n(x)\left(L-\frac{7}{8}m_n\right)
 \leqslant L-\frac{7}{8}m_{p+1},
\end{equation*}
\begin{equation*}
 \sum\limits_n\|D\psi_n(x)\|\|v_n(x)-u(x)\|
 \leqslant C(\eta_p+\eta_{p+1})
 \leqslant\frac{C(m_{p+1}+m_{p+2})}{8(C+1)}
 \leqslant\frac{1}{4}m_{p+1}.
\end{equation*}
Therefore
\begin{equation}\label{eq:lipschitz-strict-final-bound}
 \|Dv(x)\|\leqslant L-\frac{5}{8}m_{p+1}<L,
\end{equation}
which proves the assertion.
\end{proof}

\begin{proof}[Proof of Theorem~\ref{thm:lipschitz-analytic}]
If $L=0$, take $g=f$. Assume now that $L>0$, and choose
\begin{equation*}
 0<a<\min\limits\left\{\frac{1}{2},\frac{\epsilon}{4(1+L)}\right\},
 \qquad u=f\circ R_a.
\end{equation*}
By Lemma~\ref{lem:radial-strict-contraction},
\begin{equation}\label{eq:lipschitz-radial-preprocessing}
 \|u-f\|_\infty\leqslant La<\frac{\epsilon}{4},\qquad
 \operatorname{Lip}(u)\leqslant L,
\end{equation}
and, for every $R>0$,
\begin{equation*}
 \operatorname{Lip}(u,\overline B_X(0,R))
 \leqslant L\left(1-\frac{a}{(1+R^2)^{\frac{3}{2}}}\right)<L.
\end{equation*}
Lemma~\ref{lem:lipschitz-strict-smooth} gives $v\in C^\infty(X,Y)$ satisfying
$\|v-u\|_\infty<\frac{\epsilon}{4}$ and $\|Dv(x)\|<L$.
By Theorem~\ref{thm:main-complete}, we may apply
Lemma~\ref{lem:strict-derivative-analytic-transfer} to $v$ and obtain $g\in C^\omega(X,Y)$ such that
\begin{equation*}
 \|g-v\|_\infty<\frac{\epsilon}{4},\qquad
 \|Dg(x)\|<L,\qquad\operatorname{Lip}(g)\leqslant L.
\end{equation*}
The three approximations satisfy
\begin{equation*}
 \|g-f\|_\infty
 \leqslant\|g-v\|_\infty+\|v-u\|_\infty+\|u-f\|_\infty
 <\frac{3\epsilon}{4}<\epsilon.
\end{equation*}
For $\ell^{2m}(\Gamma)$ and $a\in\ell^\infty(\Gamma)$,
\begin{equation*}
 \|(a_\gamma x_\gamma)_{\gamma\in\Gamma}\|_{2m}^{2m}
 =\sum\limits_{\gamma\in\Gamma}|a_\gamma|^{2m}|x_\gamma|^{2m}
 \leqslant\|a\|_\infty^{2m}\|x\|_{2m}^{2m}.
\end{equation*}
Thus the standard basis satisfies~\eqref{eq:unconditional-decomposition-multiplier}
in the usual norm, and $x\mapsto\|x\|^{2m}$ is a separating polynomial.
The assertion therefore holds in the usual norm.
\end{proof}

\begin{theorem}
\label{thm:scalar-exact-lipschitz}
Let $X$ be a real Banach space satisfying one of the following conditions:
\begin{enumerate}[label=\textup{(\roman*)}]
 \item\label{case:scalar-lipschitz-separable} $X$ is separable and admits a separating polynomial;
 \item\label{case:scalar-lipschitz-main} $X$ belongs to the class of domain spaces in Theorem~\ref{thm:main-complete}.
\end{enumerate}
Fix any given equivalent norm on $X$, and use it to compute Lipschitz constants and derivative norms.
For every real-valued $L$-Lipschitz function $f:X\to\mathbb R$
and every constant $\epsilon>0$, there is a function $g\in C^\omega(X,\mathbb R)$ such that
\begin{equation}\label{eq:scalar-exact-lipschitz-result}
 \|g-f\|_\infty<\epsilon,\qquad
 \operatorname{Lip}(g)\leqslant L.
\end{equation}
If $L>0$, we may also require that $\|Dg(x)\|<L$ for every $x\in X$.
\end{theorem}
\begin{proof}
If $L=0$, take $g=f$. Henceforth assume $L>0$.
In both cases, $X$ admits a separating polynomial and is therefore superreflexive; see
\cite[Chapter~5, Theorem~64, pp.~269--270]{HajekJohanis2014}.
By \cite[Theorem~11.37, p.~454]{Pisier2016}, $X$ admits an equivalent norm whose modulus of smoothness has power type $p$ for some $1<p\leqslant2$.
By \cite[Theorem~9, p.~2618]{Johanis2024}, we may then choose $v\in C^1(X,\mathbb R)$ such that
\begin{equation}\label{eq:scalar-johanis-smoothing}
 \|v-f\|_\infty\leqslant\frac{\epsilon}{4},\qquad
 \operatorname{Lip}(v)\leqslant L.
\end{equation}
Here the Lipschitz constant is computed with respect to the given norm.

Choose $0<a\leqslant\min\limits\{\frac{1}{2},\frac{\epsilon}{4}\}$,
and apply Lemma~\ref{lem:radial-strict-contraction} with $E=\mathbb R$ to obtain the mapping
\begin{equation}\label{eq:scalar-strict-range-compression}
 S_a(t)=t-\frac{at}{\sqrt{1+t^2}},\qquad w=S_a\circ v.
\end{equation}
Differentiation gives
\begin{equation*}
 S_a'(t)=1-\frac{a}{(1+t^2)^{\frac{3}{2}}},\qquad
 Dw(x)[h]=\left(1-\frac{a}{(1+v(x)^2)^{\frac{3}{2}}}\right)Dv(x)[h].
\end{equation*}
Together with~\eqref{eq:radial-strict-bounds}, these identities give
\begin{equation}\label{eq:scalar-strict-c1-model}
 \begin{gathered}
 \|w-f\|_\infty\leqslant a+\frac{\epsilon}{4}\leqslant\frac{\epsilon}{2},\\
 \|Dw(x)\|\leqslant
 L\left(1-\frac{a}{(1+v(x)^2)^{\frac{3}{2}}}\right)<L.
 \end{gathered}
\end{equation}
In case~\ref{case:scalar-lipschitz-separable}, apply Theorem~\ref{thm:main-target} with $Y=\mathbb R$;
in case~\ref{case:scalar-lipschitz-main}, use Theorem~\ref{thm:main-complete}.
In either case the hypotheses of Lemma~\ref{lem:strict-derivative-analytic-transfer} are satisfied,
so there is a real-analytic function $g$ such that $\|g-w\|_\infty<\frac{\epsilon}{4}$, $\|Dg(x)\|<L$, and
$\operatorname{Lip}(g)\leqslant L$. Thus
$\|g-f\|_\infty<\frac{3\epsilon}{4}<\epsilon$, as required.
\end{proof}

\begin{remark}\label{rem:scalar-lipschitz-question}
Theorem~\ref{thm:scalar-exact-lipschitz} answers the Lipschitz-constant question
of Azagra, Fry, and Keener \cite[third question in Section~1]{AFK2012}
for the stated domain classes. It gives
$\operatorname{Lip}(g)\leqslant\operatorname{Lip}(f)$ for the prescribed norm,
and hence the requested bounds with every $C>1$.
\end{remark}

\subsection[Whitney extension and analyticity]{Whitney extension and analyticity off closed sets}\label{sec:hilbert-whitney-applications}

\begin{proposition}\label{thm:zero-order-whitney}
Let $X$ be a real Banach space, and suppose that every continuous Banach-valued mapping on every open subset of $X$
admits fine $C^0$ real-analytic approximation.
For every real Banach space $Y$, closed set $A\subset X$, and $f\in C(A,Y)$,
there is a mapping $F\in C(X,Y)$ satisfying
\begin{equation}\label{eq:zero-order-whitney-extension}
 F|_A=f,\qquad F|_{X\setminus A}\in C^\omega(X\setminus A,Y).
\end{equation}
\end{proposition}
\begin{proof}
If $A=\varnothing$, take $F=0$; if $A=X$, take $F=f$.
Otherwise, the Dugundji extension theorem
\cite[Theorem~4.1, pp.~357--358]{Dugundji1951} gives $F_0\in C(X,Y)$ with $F_0|_A=f$.
On $U=X\setminus A$, choose a real-analytic mapping $g$ satisfying
$\|g(x)-F_0(x)\|<\min\limits\{1,\operatorname{dist}(x,A)\}$,
and define
\begin{equation*}
 F(x)=\begin{cases}g(x),&x\in X\setminus A,\\f(x),&x\in A.\end{cases}
\end{equation*}
For $a\in A$ and $x\in U$,
\begin{equation*}
 \|F(x)-F(a)\|<\|x-a\|+\|F_0(x)-F_0(a)\|\longrightarrow0
\end{equation*}
as $x\to a$.
Since $F=F_0$ on $A$, $F$ is continuous at every $a\in A$.
\end{proof}

Let $H,K$ be real Hilbert spaces, let $A\subset H$ be closed, and let $f:A\to K$ and
$T:A\to\mathcal L(H,K)$.
Following \cite[Definition~25, p.~13]{JohanisZajicek2025},
we say that $(f,T)$ satisfies the \emph{first-order Whitney condition} if $f$ is continuous,
$T$ is continuous in the operator norm, and, for every $a\in A$,
\begin{equation}\label{eq:hilbert-whitney-condition}
 \lim\limits_{r\downarrow0}
 \sup\limits_{\substack{x,y\in A\cap B_H(a,r)\\x\ne y}}
 \frac{\|f(x)-f(y)-T(a)(x-y)\|}{\|x-y\|}=0,
\end{equation}
where the supremum over the empty set is taken to be zero.

\begin{theorem}\label{thm:hilbert-whitney}
Let $H,K$ be real Hilbert spaces, and let $A\subset H$ be closed.
If $f:A\to K$ and $T:A\to\mathcal L(H,K)$ satisfy the first-order Whitney condition,
then there is a mapping $F\in C^1(H,K)$ satisfying
\begin{equation}\label{eq:hilbert-whitney-result}
 F|_A=f,\qquad DF|_A=T,\qquad
 F|_{H\setminus A}\in C^\omega(H\setminus A,K).
\end{equation}
If, in addition, $\operatorname{Lip}(f)\leqslant L$ and
$\sup\limits_{a\in A}\|T(a)\|\leqslant L$, then, for every $\lambda>1$, we may also require that
\begin{equation}\label{eq:hilbert-whitney-lipschitz}
 \operatorname{Lip}(F)\leqslant\lambda L.
\end{equation}
\end{theorem}
\begin{proof}
If $A=\varnothing$, take a constant mapping. If $A=H$ and $a\in H$, then
\begin{equation*}
 \begin{aligned}
 &\frac{\|f(a+h)-f(a)-T(a)h\|}{\|h\|}\\
 &\qquad\leqslant\sup\limits_{\substack{x,y\in A\cap B_H(a,2\|h\|)\\x\ne y}}
 \frac{\|f(x)-f(y)-T(a)(x-y)\|}{\|x-y\|}\longrightarrow0
 \end{aligned}
\end{equation*}
as $h\to0$, $h\ne0$. Thus $Df=T$, and continuity of $T$ gives
$f\in C^1(H,K)$; take $F=f$. Henceforth assume $\varnothing\ne A\ne H$.

\firstproofstep{}\label{step:hilbert-whitney:1}
We first verify the conditions in \cite[Definitions~8 and~16, pp.~7 and~9]{JohanisZajicek2025},
namely $(\mathrm{LE})$ and $(\mathrm{LA}_1)$.

By the Kirszbraun--Valentine extension theorem
\cite[Theorem~5 and the subsequent corollary, p.~92]{Valentine1945},
every $L$-Lipschitz mapping $f_0:B\to K$ defined on an arbitrary subset $B\subset H$
has an $L$-Lipschitz extension $\widetilde f_0:H\to K$.
Thus $(H,K)$ satisfies $(\mathrm{LE})$ with $C_E=1$.

To verify $(\mathrm{LA}_1)$, take an $L$-Lipschitz mapping
$v:B_H(0,1)\to K$ and $\eta>0$.
By $(\mathrm{LE})$, choose an extension $V:H\to K$
with $\operatorname{Lip}(V)\leqslant L$.
Theorem~\ref{thm:lipschitz-analytic} gives $W\in C^\omega(H,K)$ satisfying
\begin{equation*}
 \|W-V\|_\infty<\eta,\qquad \operatorname{Lip}(W)\leqslant L.
\end{equation*}
Set $w=W|_{B_H(0,1)}\in C^1(B_H(0,1),K)$. Then
$\|w-v\|_{\infty,B_H(0,1)}<\eta$ and $\operatorname{Lip}(w)\leqslant L$.
Hence $(H,K)$ satisfies $(\mathrm{LA}_1)$ with $C_A=1$.

\proofstep{}\label{step:hilbert-whitney:2}
Without the additional Lipschitz assumptions, use
\cite[Theorem~34, p.~23]{JohanisZajicek2025}
to choose $F_0\in C^1(H,K)$ satisfying
\begin{equation}\label{eq:hilbert-whitney-pre-extension}
 F_0|_A=f,\qquad DF_0|_A=T,
\end{equation}
and set $c=1$.

Under the additional Lipschitz assumptions of the theorem, the case $L=0$ is handled by a constant extension.
If $L>0$, given $\lambda>1$, choose $1<\lambda_0<\lambda$.
By \cite[Theorem~35, p.~24]{JohanisZajicek2025}
and the equalities $C_E=C_A=1$ established in Step~\ref{step:hilbert-whitney:1},
choose a mapping satisfying \eqref{eq:hilbert-whitney-pre-extension},
denoted by $F_0\in C^1(H,K)$, such that
\begin{equation*}
 \operatorname{Lip}(F_0)\leqslant\lambda_0L,\qquad
 c=\frac{(\lambda-\lambda_0)L}{2}>0.
\end{equation*}

In both cases, write $\Omega=H\setminus A$ and set
$e(x)=\min\limits\{\operatorname{dist}(x,A)^2,c\}$ on $\Omega$.
By Theorem~\ref{thm:main-complete}, choose $h\in C^\omega(\Omega,K)$ such that, for every $x\in\Omega$,
\begin{equation}\label{eq:hilbert-whitney-exterior-approximation}
 \max\limits_{j=0,1}\|D^jh(x)-D^jF_0(x)\|<e(x).
\end{equation}
Set $F=h$ on $\Omega$ and $F=F_0$ on $A$, and write $R=F-F_0$.
Since $R|_A=0$ and $e(x)\leqslant\operatorname{dist}(x,A)^2$,
for every $a\in A$ and $z\in H$ we have
\begin{equation*}
 \|R(a+z)-R(a)\|
 \leqslant\operatorname{dist}(a+z,A)^2
 \leqslant\|z\|^2.
\end{equation*}
In particular,
\begin{equation*}
 \frac{\|R(a+z)-R(a)-0[z]\|}{\|z\|}\leqslant\|z\|\longrightarrow0
\end{equation*}
as $z\to0$ with $z\ne0$, so $R$ is Fr\'echet differentiable at $a$ and $DR(a)=0$.
On $\Omega$, we have $DR=Dh-DF_0$, and
\eqref{eq:hilbert-whitney-exterior-approximation} gives
$\|DR(x)\|<\operatorname{dist}(x,A)^2$.
Consequently, for every $a\in A$ and $x\in H$,
\begin{equation*}
 \|DR(x)-DR(a)\|\leqslant\|x-a\|^2.
\end{equation*}
This proves that $DR$ is continuous at every point of $A$.
Moreover, $R|_\Omega=h-F_0\in C^1(\Omega,K)$.
Hence $F\in C^1(H,K)$, and \eqref{eq:hilbert-whitney-pre-extension}
gives $F|_A=f$ and $DF|_A=T$.

\proofstep{}\label{step:hilbert-whitney:3}
Under the additional Lipschitz assumptions with $L>0$, the inequality
$e\leqslant c=\frac{(\lambda-\lambda_0)L}{2}$
gives, for every $x\in\Omega$,
\begin{equation*}
 \begin{aligned}
 \|DF(x)\|
 &\leqslant\|DF_0(x)\|+\|Dh(x)-DF_0(x)\|\\
 &\leqslant\lambda_0L+\frac{(\lambda-\lambda_0)L}{2}
 =\frac{\lambda+\lambda_0}{2}L.
 \end{aligned}
\end{equation*}
On $A$, we have $DF=T$ and $\|T(a)\|\leqslant L$, so
\begin{equation*}
 \sup\limits_{z\in H}\|DF(z)\|
 \leqslant\frac{\lambda+\lambda_0}{2}L<\lambda L.
\end{equation*}
For any $x,y\in H$, the fundamental theorem of calculus gives
\begin{equation*}
 F(y)-F(x)
 =\int\limits_0^1DF(x+t(y-x))(y-x)\,\mathrm{d}t.
\end{equation*}
Thus
\begin{equation*}
 \|F(y)-F(x)\|
 \leqslant\sup\limits_{z\in H}\|DF(z)\|\,\|y-x\|
 \leqslant\lambda L\|y-x\|,
\end{equation*}
which proves \eqref{eq:hilbert-whitney-lipschitz}.
\end{proof}

\begin{remark}\label{rem:hilbert-whitney-question}
Taking $K=H$ in Theorem~\ref{thm:hilbert-whitney} answers the question
of Johanis and Zaj\'i\v{c}ek \cite[Section~2.6]{JohanisZajicek2025}
about basic and Lipschitz $C^1$ Whitney extension for a nonseparable
Hilbert pair $(H,H)$. The extension also has the prescribed derivative
on $A$ and is real analytic on $H\setminus A$.
\end{remark}

\subsection[Approximation with surjective derivatives]{Real-analytic approximation without critical points and with surjective derivatives}\label{sec:hilbert-geometric-applications}

For a real-valued $C^1$ function $f$, a point $x$ is called a critical point of $f$ if $Df(x)=0$.  Azagra and Cepedello Boiso proved that every continuous
mapping from a separable infinite-dimensional Hilbert
space into a finite-dimensional Euclidean space admits
uniform approximation by smooth mappings without
critical points \cite{AzagraCepedello2004}.

\begin{proposition}
\label{prop:geo-surjective-transfer}
Let $E,F$ be real Banach spaces with $F\ne\{0\}$, and let $U\subset E$ be open.
Suppose that mappings in $C^1(U,F)$ admit fine $C^1$ real-analytic approximation.
If, for given $f\in C(U,F)$ and $\varepsilon\in C(U,(0,\infty))$,
there is a mapping $h\in C^1(U,F)$ such that, for every $x\in U$,
\begin{equation}\label{eq:geo-smooth-surjective-input}
 \|h(x)-f(x)\|<\frac{\varepsilon(x)}{3},
 \qquad Dh(x):E\longrightarrow F\text{ is surjective},
\end{equation}
then there is a mapping $g\in C^\omega(U,F)$ such that, for every $x\in U$,
\begin{equation}\label{eq:geo-analytic-surjective-output}
 \|g(x)-f(x)\|<\varepsilon(x),
 \qquad Dg(x):E\longrightarrow F\text{ is surjective}.
\end{equation}
\end{proposition}

\begin{proof}
Write
	\begin{equation*}
	\mathcal O
	=\{T\in\mathcal L(E,F):T(E)=F\}.
	\end{equation*}
The set $\mathcal O$ is open in the operator norm topology; see also
\cite[proof of Proposition~1.8, p.~77]{AzagraDobrowolskiGarcia2019}.
For completeness, fix $T\in\mathcal O$. The open mapping theorem gives
$C>0$ such that, for every $y\in F$, there exists $x\in E$ with
\begin{equation*}
 Tx=y,\qquad \|x\|\leqslant C\|y\|.
\end{equation*}
Let $S\in\mathcal L(E,F)$ satisfy $q=C\|S-T\|<1$.
For $y\in F$, set $r_0=y$ and inductively choose $x_n\in E$ with
\begin{equation*}
 Tx_n=r_n,\qquad \|x_n\|\leqslant C\|r_n\|,\qquad
 r_{n+1}=r_n-Sx_n=(T-S)x_n.
\end{equation*}
Then $\|r_n\|\leqslant q^n\|y\|$ and
$\sum\limits_{n\geqslant0}\|x_n\|\leqslant\frac{C}{1-q}\|y\|$.
Thus $x=\sum\limits_{n\geqslant0}x_n\in E$ satisfies
\begin{equation*}
 Sx=\lim\limits_{N\to\infty}\sum\limits_{n=0}^NSx_n
    =\lim\limits_{N\to\infty}(r_0-r_{N+1})=y.
\end{equation*}
Therefore $B_{\mathcal L(E,F)}(T,\frac{1}{C})\subset\mathcal O$.
Since $F\ne\{0\}$, we have $0\notin\mathcal O$.
By \eqref{eq:geo-smooth-surjective-input},
$Dh(x)\in\mathcal O$ for every $x\in U$.
Therefore
	\begin{equation}\label{eq:geo-surjective-delta-definition}
		\delta(x)=\min\limits\left\{
		\frac{\varepsilon(x)}{3},
		\frac{1}{2}\operatorname{dist}
		\bigl(Dh(x),\mathcal L(E,F)\setminus\mathcal O\bigr)
		\right\}
	\end{equation}
is a finite, continuous, strictly positive function on $U$.

By hypothesis,
there is a mapping $g\in C^\omega(U,F)$ such that, for every $x\in U$,
	\begin{equation*}
	\|g(x)-h(x)\|<\delta(x),
	\qquad
	\|Dg(x)-Dh(x)\|<\delta(x).
	\end{equation*}

Set $\mathcal C=\mathcal L(E,F)\setminus\mathcal O$.
The distance inequality and~\eqref{eq:geo-surjective-delta-definition} give
\begin{equation*}
 \begin{aligned}
 \operatorname{dist}(Dg(x),\mathcal C)
 &\geqslant\operatorname{dist}(Dh(x),\mathcal C)-\|Dg(x)-Dh(x)\|\\
 &>\frac{1}{2}\operatorname{dist}(Dh(x),\mathcal C)>0.
 \end{aligned}
\end{equation*}
Hence $Dg(x)\in\mathcal O$, so $Dg(x):E\to F$ is surjective.
On the other hand, \eqref{eq:geo-smooth-surjective-input} gives
	\begin{equation*}
	\begin{aligned}
		\|g(x)-f(x)\|
		&\leqslant
		\|g(x)-h(x)\|+\|h(x)-f(x)\|<\frac{2\varepsilon(x)}{3}
		<\varepsilon(x).
	\end{aligned}
	\end{equation*}
Thus $g$ satisfies \eqref{eq:geo-analytic-surjective-output}.
\end{proof}

\begin{corollary}
\label{cor:geo-vector-surjective-examples}
In each of the following cases, for every continuous mapping $f:U\to F$ and every
$\varepsilon\in C(U,(0,\infty))$ there is
a real-analytic approximation satisfying \eqref{eq:geo-analytic-surjective-output}:
\begin{enumerate}[label=\textup{(\roman*)}]
\item\label{case:surjective-hilbert-full} $E$ is an arbitrary infinite-dimensional real Hilbert space, $U=E$,
and $F$ is a nonzero Hilbert quotient of $E$, meaning that there is a bounded linear surjection $E\to F$.
\item\label{case:surjective-hilbert-open} $E$ is a separable infinite-dimensional real Hilbert space, $U\subset E$ is an arbitrary open set,
and $F$ is a nonzero separable real Hilbert space.
\item\label{case:surjective-ell} $E=\ell^{2m}(\Gamma)$, where $m$ is a positive integer and $\Gamma$ is an arbitrary infinite set,
$U=E$, and $F$ is an arbitrary nonzero Banach quotient of $E$.
\item\label{case:surjective-Lp} $E=L^{2m}[0,1]$, where $m$ is a positive integer, $U=E$,
and $F$ is an arbitrary nonzero Banach quotient of $E$.
\item\label{case:surjective-separable-finite} $E$ is an arbitrary separable infinite-dimensional real Banach space admitting a separating polynomial,
$U=E$, and $F=\mathbb R^m$ for some $m\geqslant1$.
\end{enumerate}
\end{corollary}
\begin{proof}
In each case, we first construct a mapping
$h\in C^1(U,F)$ satisfying \eqref{eq:geo-smooth-surjective-input}.
In cases~\ref{case:surjective-hilbert-full} and~\ref{case:surjective-ell},
use \cite[Corollary~1.3, p.~6; Theorem~1.2, p.~5]{AzagraGarciaJimenez2024}, respectively.
For case~\ref{case:surjective-hilbert-open}, use
\cite[Theorem~1.1, p.~2; discussion preceding Theorem~1.5, p.~6]{AzagraDobrowolskiGarcia2019};
for case~\ref{case:surjective-Lp}, use \cite[Theorem~1.2, p.~2]{AzagraDobrowolskiGarcia2019}.
Choosing $\frac{\varepsilon}{4}\in C(U,(0,\infty))$ in each cited result gives
$\|h(x)-f(x)\|\leqslant\frac{\varepsilon(x)}{4}<\frac{\varepsilon(x)}{3}$
for every $x\in U$, as required in~\eqref{eq:geo-smooth-surjective-input}.
Fine $C^1$ real-analytic approximation in these four cases follows from Theorem~\ref{thm:main-complete},
so Proposition~\ref{prop:geo-surjective-transfer} applies.

In case~\ref{case:surjective-separable-finite},
Lemma~\ref{lem:target-countable-partition} gives smooth partitions of unity on $E$.
In \cite[Corollary~1.9, p.~8]{AzagraGarciaJimenez2024}, take $E$ itself as the separable complemented subspace
and choose the positive continuous function $\frac{\varepsilon}{4}$ to obtain the required $C^1$ mapping $h$ with finite-dimensional range.
Theorem~\ref{thm:main-target} supplies fine $C^1$ real-analytic approximation in this case,
and Proposition~\ref{prop:geo-surjective-transfer} yields the conclusion.
\end{proof}

\begin{remark}\label{rem:noncritical-question}
Taking $E=U=\ell^2$ and $F=\mathbb R$ in
Corollary~\ref{cor:geo-vector-surjective-examples} answers affirmatively
Jim\'enez-Sevilla's concluding question \cite[Question]{JimenezSevilla2008}
about uniform real-analytic approximation of continuous real-valued
functions on $\ell^2$ without critical points. 
\end{remark}

\subsection[Analytic equivalence and isotopy]{Real-analytic equivalence and isotopy of open sets}\label{sec:analytic-equivalence}

For topological spaces $U,V$, a family $(h_t:U\to V)_{t\in[0,1]}$
is a \emph{topological isotopy} if the mapping
\begin{equation*}
 [0,1]\times U\longrightarrow[0,1]\times V,\qquad
 (t,x)\longmapsto(t,h_t(x))
\end{equation*}
is a homeomorphism.

\begin{lemma}\label{lem:diffeomorphism-small-perturbation}
Let $E$ be a real Banach space, and let $V\subset E$ be a nonempty open set. For every $y\in V$, set
\begin{equation*}
 r_V(y)=\min\limits\{1,\operatorname{dist}(y,E\setminus V)\}.
\end{equation*}
If $q\in C^1(V,E)$ satisfies, for every $y\in V$,
\begin{equation}\label{eq:diffeomorphism-perturbation-hypotheses}
 \|q(y)\|<\frac{r_V(y)}{8},\qquad
 \|Dq(y)\|\leqslant\frac{1}{8},
\end{equation}
then $h_t=I_V+tq$ is a $C^1$ diffeomorphism of $V$ onto itself for every $t\in[0,1]$.
Moreover, $(t,y)\mapsto h_t^{-1}(y)$ is continuous, so $(h_t)$ is a topological isotopy.
\end{lemma}
\begin{proof}
Inspired by \cite[Section~2, Corollary~1 and Lemmas~2--3, pp.~326--327]{Moulis1971},
we extend $q$ by zero on $E\setminus V$ and denote the extension by $\overline q$.
If $x,y\in V$ and the line segment $[x,y]$ joining $x$ and $y$ is contained in $V$, the mean value inequality gives
$\|q(x)-q(y)\|\leqslant\frac{1}{8}\|x-y\|$.
If $[x,y]$ contains a point $z\in E\setminus V$, then
\begin{equation*}
 \begin{aligned}
 \|q(x)-q(y)\|
 &\leqslant\frac{1}{8}\bigl(\operatorname{dist}(x,E\setminus V)
                         +\operatorname{dist}(y,E\setminus V)\bigr)\\
 &\leqslant\frac{1}{8}\bigl(\|x-z\|+\|z-y\|\bigr)
 =\frac{1}{8}\|x-y\|.
 \end{aligned}
\end{equation*}
If $x\in V$ and $y\in E\setminus V$, then
\begin{equation*}
 \|\overline q(x)-\overline q(y)\|=\|q(x)\|
 \leqslant\frac{1}{8}\operatorname{dist}(x,E\setminus V)
 \leqslant\frac{1}{8}\|x-y\|.
\end{equation*}
Thus
\begin{equation*}
 \operatorname{Lip}(\overline q,E)\leqslant\frac{1}{8},\qquad
 \|\overline q\|_\infty\leqslant\frac{1}{8}.
\end{equation*}
For $t\in[0,1]$ and $z\in E$, set $C_{t,z}(y)=z-t\overline q(y)$.
Then
\begin{equation*}
 \|C_{t,z}(u)-C_{t,z}(v)\|\leqslant\frac{t}{8}\|u-v\|.
\end{equation*}
The contraction mapping theorem gives a unique fixed point $y$.
Explicitly, with $y_0=z$ and $y_{n+1}=z-t\overline q(y_n)$,
\begin{equation*}
 y=\lim\limits_{n\to\infty}y_n,\qquad
 \|y-y_n\|\leqslant\frac{(\frac{t}{8})^n}{1-\frac{t}{8}}\|y_1-y_0\|.
\end{equation*}
Thus $\overline h_t=I_E+t\overline q$ is a bijection of $E$ onto itself.
For $z\in E\setminus V$, we have $\overline h_t(z)=z$; uniqueness of its
preimage then gives
\begin{equation*}
 \overline h_t^{-1}(E\setminus V)=E\setminus V,\qquad
 \overline h_t(V)=V.
\end{equation*}
By \eqref{eq:diffeomorphism-perturbation-hypotheses},
$\|tDq(y)\|\leqslant\frac{1}{8}<1$,
so the Neumann series gives
\begin{equation*}
 \begin{gathered}
 (Dh_t(y))^{-1}=(I_E+tDq(y))^{-1}
   =\sum\limits_{n=0}^\infty(-tDq(y))^n,\\
 \|(Dh_t(y))^{-1}\|\leqslant\frac{1}{1-\frac{t}{8}}\leqslant\frac{8}{7}.
 \end{gathered}
\end{equation*}
The inverse function theorem gives $h_t^{-1}\in C^1(V,V)$ and
\begin{equation*}
 Dh_t^{-1}(h_t(y))=(I_E+tDq(y))^{-1}.
\end{equation*}

For $u=h_t^{-1}(y)$ and $v=h_s^{-1}(z)$, the identities
$u=y-t\overline q(u)$ and $v=z-s\overline q(v)$ give
\begin{equation*}
 \|u-v\|\leqslant\|y-z\|+\frac{1}{8}\|u-v\|+\frac{1}{8}|t-s|.
\end{equation*}
Hence
\begin{equation}\label{eq:diffeomorphism-inverse-joint-bound}
 \|h_t^{-1}(y)-h_s^{-1}(z)\|
 \leqslant\frac{8}{7}\|y-z\|+\frac{1}{7}|t-s|.
\end{equation}
The mapping $(t,x)\mapsto(t,h_t(x))$ is continuous and bijective, and \eqref{eq:diffeomorphism-inverse-joint-bound}
ensures that its inverse is continuous. Therefore $(h_t)$ is a topological isotopy.
\end{proof}

\begin{theorem}\label{thm:analytic-diffeomorphism}
Let $X$ be any real Banach space to which Theorem~\ref{thm:main-complete} applies,
let $U,V\subset X$ be open,
and let $f:U\to V$ be a $C^1$ diffeomorphism.
For every $\varepsilon\in C(U,(0,\infty))$, there is a real-analytic diffeomorphism
$g:U\to V$ such that, for every $x\in U$,
\begin{equation}\label{eq:analytic-diffeomorphism-fine-approximation}
 \|g(x)-f(x)\|<\varepsilon(x),\qquad
 \|Dg(x)-Df(x)\|<\varepsilon(x).
\end{equation}
Moreover, for every $t\in[0,1]$,
\begin{equation}\label{eq:analytic-diffeomorphism-isotopy}
 f_t=(1-t)f+tg
\end{equation}
is a $C^1$ diffeomorphism from $U$ onto the same set $V$, and $(f_t)$ is a topological isotopy.
In particular, the assertion holds for every real Hilbert space.
\end{theorem}

\begin{proof}
If $U=\varnothing$ or $X=\{0\}$, take $g=f$.
Henceforth assume $U\ne\varnothing$ and $X\ne\{0\}$.
For every $x\in U$, set
\begin{equation}\label{eq:analytic-diffeomorphism-delta-definition}
 \delta(x)=
 \min\limits\left\{
 \frac{\varepsilon(x)}{2},
 \frac{r_V(f(x))}{8},
 \frac{1}{8\|Df(x)^{-1}\|}
 \right\},
\end{equation}
where $r_V$ is defined in Lemma~\ref{lem:diffeomorphism-small-perturbation}.
Since $f$ is a $C^1$ diffeomorphism, each $Df(x)$ is a bounded linear isomorphism.
The continuity of $Df$ and of inversion on $GL(X)$ imply that
$\delta$ is a continuous positive function on $U$.

By Theorem~\ref{thm:main-complete}, choose $g\in C^\omega(U,X)$ such that, for every $x\in U$,
\begin{equation}\label{eq:analytic-diffeomorphism-proof-approx}
 \|g(x)-f(x)\|<\delta(x),\qquad
 \|Dg(x)-Df(x)\|<\delta(x).
\end{equation}
Define $q=g\circ f^{-1}-I_V$ on $V$. Then $q\in C^1(V,X)$.
For any $y\in V$, set $x=f^{-1}(y)$.
The chain rule gives
\begin{equation*}
 Dq(y)
 =Dg(x)Df(x)^{-1}-I_X
 =(Dg(x)-Df(x))Df(x)^{-1}.
\end{equation*}
Combining \eqref{eq:analytic-diffeomorphism-delta-definition} and
\eqref{eq:analytic-diffeomorphism-proof-approx}, we obtain
\begin{equation*}
 \begin{aligned}
 \|q(y)\|
 &=\|g(x)-f(x)\|<\delta(x)\leqslant\frac{r_V(y)}{8},\\
 \|Dq(y)\|
 &\leqslant\|Dg(x)-Df(x)\|\,\|Df(x)^{-1}\|\\
 &<\delta(x)\|Df(x)^{-1}\|\leqslant\frac{1}{8}.
 \end{aligned}
\end{equation*}
We may therefore apply Lemma~\ref{lem:diffeomorphism-small-perturbation}:
$h_t=I_V+tq$ is a $C^1$ diffeomorphism of $V$ onto itself,
and $(t,y)\mapsto h_t^{-1}(y)$ is jointly continuous.
Since $q\circ f=g-f$,
\begin{equation*}
 f_t=(1-t)f+tg=h_t\circ f,\qquad
 f_t^{-1}=f^{-1}\circ h_t^{-1}.
\end{equation*}
Thus each $f_t$ is a $C^1$ diffeomorphism from $U$ onto $V$.
Both $(t,x)\mapsto f_t(x)$ and $(t,y)\mapsto f_t^{-1}(y)$ are jointly continuous,
so $(f_t)$ is a topological isotopy from $f$ to $g=f_1$.

Since $g:U\to V$ is a $C^1$ diffeomorphism, $Dg(x)$ is invertible at every point.
By the analytic inverse function theorem
\cite[Theorem~15.3(b), pp.~151--152]{Deimling1985},
$g$ has a real-analytic local inverse near every point.
By injectivity of $g$, these local inverses agree with the global inverse $g^{-1}$.
Hence $g^{-1}\in C^\omega(V,U)$, and $g$ is a real-analytic diffeomorphism.
Finally, \eqref{eq:analytic-diffeomorphism-proof-approx}
and $\delta\leqslant\frac{\varepsilon}{2}$ give
\eqref{eq:analytic-diffeomorphism-fine-approximation}.
\end{proof}

\begin{remark}\label{rem:hilbert-analytic-equivalence-question}
For a separable infinite-dimensional real Hilbert space,
Theorem~\ref{thm:analytic-diffeomorphism} answers the subquestion in
Anderson, Chapman, and Schori's problem list
\cite[Problem~D1, p.~22]{AndersonChapmanSchori1971}
asking whether two $C^\infty$-diffeomorphic open subsets are
real-analytically equivalent. The theorem requires only a $C^1$
diffeomorphism and also applies to nonseparable Hilbert spaces.
\end{remark}

\begin{corollary}
\label{cor:geo-homeomorphism-analytic-isotopy}
Let $H$ be a separable infinite-dimensional real Hilbert space, and let $U\subset H$ be open.
Every homeomorphism $h:U\to U$ is topologically isotopic to a real-analytic diffeomorphism $g:U\to U$.
In particular, $h$ and $g$ are homotopic.
\end{corollary}
\begin{proof}
The assertion is immediate if $U=\varnothing$. Otherwise, $U$ is a separable infinite-dimensional
Hilbert manifold. By Burghelea and Henderson
\cite[Theorem~0.3, p.~1262]{BurgheleaHenderson1970},
$h$ is topologically isotopic to a $C^\infty$ diffeomorphism
$h_0:U\to U$.
Denote this isotopy by $(k_t)_{0\leqslant t\leqslant1}$,
where $k_0=h$ and $k_1=h_0$.

Theorem~\ref{thm:analytic-diffeomorphism}
gives a real-analytic diffeomorphism $g:U\to U$ such that the family defined, for every $t\in[0,1]$, by
	\begin{equation}\label{eq:geo-smooth-to-analytic-isotopy}
		h_t=(1-t)h_0+tg
	\end{equation}
is a topological isotopy from $h_0$ to $g$.

For $0\leqslant t\leqslant1$, define
	\begin{equation*}
	H_t(x)=
	\begin{cases}
		k_{2t}(x), & 0\leqslant t\leqslant\frac{1}{2},\\[2mm]
		h_{2t-1}(x), & \frac{1}{2}\leqslant t\leqslant1.
	\end{cases}
	\end{equation*}
Since $k_1=h_0$, both formulas defining $H_t$ agree with $h_0$
at $t=\frac{1}{2}$. Each $H_t$ is
a homeomorphism of $U$ onto itself, with $H_0=h$ and $H_1=g$.
The inverse is
\begin{equation*}
 H_t^{-1}(y)=
 \begin{cases}
  k_{2t}^{-1}(y),&0\leqslant t\leqslant\frac{1}{2},\\
  h_{2t-1}^{-1}(y),&\frac{1}{2}\leqslant t\leqslant1.
 \end{cases}
\end{equation*}
At $t=\frac{1}{2}$, the two values are $k_1^{-1}(y)=h_0^{-1}(y)$.
By the joint continuity of the two isotopies and the pasting lemma,
both $(t,x)\mapsto H_t(x)$ and $(t,y)\mapsto H_t^{-1}(y)$ are continuous.
Thus $(H_t)_{0\leqslant t\leqslant1}$
is a topological isotopy from $h$ to $g$.
\end{proof}

\begin{remark}\label{rem:hilbert-analytic-isotopy-question}
Corollary~\ref{cor:geo-homeomorphism-analytic-isotopy} answers the subquestion
in Anderson, Chapman, and Schori's problem list
\cite[Problem~D1, p.~22]{AndersonChapmanSchori1971}
asking whether every homeomorphism of an open subset of a separable
infinite-dimensional Hilbert space onto itself is homotopic to
a real-analytic diffeomorphism. The corollary gives a topological isotopy.
\end{remark}

\end{document}